\documentclass[12pt, reqno]{amsart}
\usepackage{amsmath}
\usepackage{amsxtra}
\usepackage{amstext}
\usepackage{amssymb}
\usepackage{amsthm}
\usepackage{latexsym}
\usepackage{dsfont} 
\usepackage{verbatim}
\usepackage{enumerate,xcolor}
\usepackage{graphicx}
\usepackage{caption}
\usepackage{setspace}

\usepackage{color}

\theoremstyle{plain}
\newtheorem{thm}{Theorem}[section]
\newtheorem{cor}[thm]{Corollary}
\newtheorem{lem}[thm]{Lemma}
\newtheorem{prop}[thm]{Proposition}

\theoremstyle{definition}
\newtheorem{defn}[thm]{Definition}
\newtheorem{exa}[thm]{Example}

\numberwithin{equation}{section}
\begin{document}
			\title[ Slice Regular Capacity and  Rational Approximation]{\mbox{Slice regular capacity and rational approximation}}
	\author{Thankarajan  Prasad}
		\address{Department of Mathematics, Nalanda University, Rajgir-803116, District- Nalanda , Bihar, India}.
		\email{{prasad@nalandauniv.edu.in},   {prasadvalapil@gmail.com}}
	    \author{Rincy Susan Raju}
	    \address{Department of Mathematics, University College (Affiliated to University of Kerala), Thiruvananthapuram, Kerala, India}
	    \email{{rincysr@universitycollege.ac.in}}
		
		\author{Enose Shine Lal}
		\address{Department of Mathematics, University College (Affiliated to University of Kerala), Thiruvananthapuram, Kerala, India}
		\email{{shinelal@universitycollege.ac.in}}

		\begin{abstract}
			Approximation theory in the complex plane is a well established area of mathematical research, with extensive contributions documented in the literature, notably through the classical results of  Weierstrass, Runge,  Mergelyan, and Vitushkin. Approximation theory in quaternionic settings, particularly rational approximation on compact sets, which studies when slice regular functions can be uniformly approximated by rational functions, is one of the recent areas of interest in quaternionic rational approximation theory (See \cite{runge quaternionic,quatern approx}). In this paper, we discuss Vitushkin approximation  theorem in  the quaternionic setting  by using the quaternionic analogue of analytic capacity  and also we establish Davies modified Vitushkin scheme in quaternionic setting.  We also discuss  quaternionic analogue of $\gamma_{+}$- capacity and Cauchy capacity and  the  remarkable  results of Tolsa \cite{tolsa} and Younsi \cite{younis}  related  to these capacities in complex plane setting.

		\end{abstract}
		\maketitle
		
		\tableofcontents
		
		\section{Introduction}
	The theory of rational approximation is the classical problem of representing holomorphic functions by simple functions, especially 	polynomials and rational functions.	
	The foundations of approximation theory in the complex plane $\mathbb{C}$ were laid by Weierstrass (1885) and Runge (1885), establishing that continuous holomorphic functions could be approximated by polynomials under suitable topological conditions\cite{Runge,Weierstrass}. Runge \cite{complex approx,conway, Runge}  proved that a function holomorphic on a compact set $K \subset \mathbb{C}$ can be uniformly approximated by rational functions with poles outside $K$, and if $\mathbb{C} \setminus K$ is connected it can be approximated by polynomials.  The Mergelyan’s theorem \cite{Mergelyan}, which provided the more significant result- any continuous function on a compact set $K$ that is holomorphic in its interior can be uniformly approximated by polynomials if and only if $\mathbb{C} \setminus K$ is connected. On more general compact sets,  Vitushkin \cite{vitushkin} developed a theory about when functions holomorphic in the interior of a compact set can be uniformly approximated on the whole set by rational functions whose poles lie outside the set in terms  of  continuous analytic capacity, which measures the size of a set from the perspective of bounded analytic functions. 
	
	The notion of analytic capacity arose from one of the problems of complex function theory which is about the characterization of removable singularities for bounded analytic functions. Painlev\'e\cite{painleve} asked for geometric characterization of compact subsets of complex space such that every bounded analytic function on its complement extends analytically across that compact set. In par with this,  Ahlfors\cite{ahlfors} introduced analytic capacity to solve function-theoretic extremal problems on finitely connected planar domains. Later,  Garabedian\cite{Garabedian} refined this framework.  Havinson\cite{havinsonmain} used analytic
	capacity to  studying problems of function theory which is partly an extension of the work of Ahlfors and Garabedian to infinitely connected domains. During 1950's and 1960's Vitushkin\cite{vitushkin} recognized that analytic capacity has some connections with rational approximation theory.   Vitushkin also introduced continuous analytic capacity  and proved that rational functions are uniformly dense in the space of continuous functions on compact set (holomorphic in the interior) if and only if continuous analytic capacity and analytic capacity satisfy specific localized geometric density ratios. Davie\cite{davie} developed  a modified approximation scheme which improves Vitushkin method by using squares of varying size which is more flexible.
	
	Vitushkin’s techniques is not a complete description of analytic capacity in metric or geometric terms. The main advancement in this study is due to Tolsa  \cite{tolsa} , who proved that analytic capacity is comparable to $\gamma_{+}$ capacity which is easier to comprehend since it is more suitable to real analysis tools. Tolsa  \cite{tolsa}  solved a very difficult problem raised by Vitushkin in 1967 which is about the semi additivity of Analytic capacity. Then Havinson \cite{havinson7} introduced Cauchy capacity and Malik Younsi \cite{younis} studied its relation with analytic Capacity.

	The non-commutativity of the quaternions and the rigid nature of Feuter regularity which does not contain the identity mapping the extension of approximation theory from complex variables  to the  quaternionic setting encounters  obstacle. But, with the introduction of slice regular functions  by Gentili and Struppa\cite{gentili} it became possible to formulate right-slice quaternionic polynomials  and investigate approximation properties.   Gal and Sabadini \cite{quatern approx} developed quaternionic analogues of namely approximation results of Runge type,
	Mergelyan type and Arakelian type. Since the Quaternionic Riemann mapping theorem is available only for some extended axially-symmetric sets which is denoted by $\Re(\mathbb{H})$, a full analogue of Mergelyan's theorem is not immediately available in the quaternionic setting.     Gal and Sabadini \cite{quatern approx}     obtained uniform polynomial approximation for several important classes of sets, including star-like sets and  $\Re(\mathbb{H})$ and  established uniform approximation result on the boundary of open sets in $\Re(\mathbb{H})$ by quaternionic polynomials. Also they obtain quantitative approximation in the case of compact balls with real centers and for Cassini cells. 
	
	The quaternionic version is not easy  to study since the quaternionic field $\mathbb H$ is non-commutative and  every non-real quaternion	is naturally associated with the conjugacy sphere thus, the geometry relevant to slice regular functions is
	different from the ordinary pointwise geometry of the complex plane. If $K\subset\mathbb H$ is axially symmetric, then its intersection with a complex slice, contains the  two-dimensional complex-slice information of $K$. The axial symmetry ensures that the geometry of $K$ is compatible with the spherical structure of  slice regularity. This motivates the introduction of a quaternionic analogue of analytic capacity, which we call the slice regular capacity( see section 9). 
	
	In this work we  extend Vitushkin theorem  to the quaternionic setting of slice regular functions. Let $K\subset\mathbb H$ be a compact axially symmetric set.  We denote the space of all continuous left slice functions on $K$ by $\mathcal C_L(K)$. Let $R_L(K)$ be the set of  functions on $K$ which can be approximated by left rational functions with poles off $K$ and $\mathcal A_L(K)$ be the set of  functions on $\mathcal C_L(K)$ and left slice regular on the interior of $K$. In quaternionic setting,  we are trying to solve the conditions under which $\mathcal R_L(K)$ is equal to $\mathcal C_L(K)$ and  $\mathcal A_L(K)$. As in the  classical complex case we are trying to solve  this problem in terms of the quaternionic counterpart of analytic capacity and continuous analytic capacity.
	
   This paper contains 3 parts. In 0th part, we lists   known results of quaternionic variables. We give the basic definitions , properties and results of slice regular functions. Also the quaternionic Laurent series and in correspondence with this the singularities and its classification has been explained. As our main aim is to study approximation on compact sets the approximation theory of complex variables are presented  in this part. We focused mainly on the Vitushkin approximation scheme and Davies modified version of the Vitushkin  approximation scheme.
   
   In  part I we introduce  slice regular capacity and a detailed study of its geometric properties. The extension of $\gamma_{+}$ and Cauchy capacity in quaternionic is also discussed.  The removable singularity and  some famous results by Tolsa\cite{tolsa} and Younsi\cite{younis} has been discussed in quaternionic setting. 
   
  In the 3rd part, we discuss the Vitushkin approximation theorem and also the quaternionic analogue of Davies modified version of Vitushkin approximation  scheme in quaternionic setting.

\part*{\MakeUppercase{Part 0:   Notions, Terminologies, and Preliminaries} }
\vspace{0.5em}
\vspace{0.5em}
In this part we list notions, terminologies and some preliminaries from the literature which will be used in the sequel.
\section{Slice Regular Functions}
In this section we mainly refer \cite{Q1,Q2,Q3,Q4}. 
Let $\mathbb H$ denote the quaternionic field which consists of the elements  of the form $p=p_0+p_1i+p_2j+p_3k,$ where $p_0,p_1,p_2,p_3\in\mathbb R,$ and the imaginary units $i,j,k$ satisfy
$i^2=j^2=k^2=-1,\ ij=k, jk=i, ki=j,$ and $ij=-ji, jk=-kj, ki=-ik.$
For $p\in\mathbb H$, its real and imaginary parts are defined by
$\operatorname{Re}(p)=p_0,$ and $\operatorname{Im}(p)=p_1i+p_2j+p_3k.$ 

The conjugate and modulus of $q$ are $\overline p=p_0-p_1i-p_2j-p_3k$ and
\[|p|=\sqrt{p\overline p} = \sqrt{p_0^2+p_1^2+p_2^2+p_3^2},\]
respectively.
Let $\mathbb S =\{p_0+p_1i+p_2j+p_3k:\ p_0,p_1,p_2,p_3\in\mathbb R,\ p_1^2+p_2^2+p_3^2=1\}$ be the unit sphere of purely imaginary quaternions.
Equivalently, $\mathbb S=\{I\in\mathbb H:I^2=-1\}.$

For each $I\in\mathbb S$, we consider the  complex plane 
$\mathbb C_I=\{x+Iy:x,y\in\mathbb R\}$ which  is isomorphic to the complex field $\mathbb C$ and every quaternion belongs to at least one such plane. Hence, 
\[\mathbb H=\bigcup_{I\in\mathbb S}\mathbb C_I.\]
 In particular, every non-real quaternion $p\in\mathbb H$ can be written as $p=u+Iv,$ where
$u=\operatorname{Re}(p),\ v=|\operatorname{Im}(p)|>0,\ I=\frac{\operatorname{Im}(p)}{|\operatorname{Im}(p)|}\in\mathbb S.$ 
By this representation the sphere associated with $q$ is defined by
\[[p]=\{u+Jv:J\in\mathbb S\}.\]
It consists of all quaternions having the same real part and
the same modulus of the imaginary part as $p$. Its center is
$u$ and its radius is $|v|$. When $v=0$, the sphere reduces
to the single point $u$.

Let $I,J\in\mathbb S$ satisfy $IJ=-JI$. Then every quaternion
$p\in\mathbb H$ admits a unique representation $p=z+wJ,\ z,w\in\mathbb C_I.$ Thus,
\[\mathbb H=\mathbb C_I+\mathbb C_IJ.\]

Nowadays, the theory underlying slice regular functions has become firmly established. In some open sets  the class of slice regular functions has  nice properties. 

For an open set $U \subseteq \mathbb{H}$,  $U$ is axially symmetric if $[q] \subset U$ for any $q \in U$ and a slice domain if $U \cap \mathbb{R} \neq \emptyset$ and if $U\cap \mathbb{C}_I$ is a domain in $\mathbb{C}_I$ for any $I \in \mathbb{S}$. The balls centered at real points  are axially symmetric slice domain.

\noindent Now, we define slice regular functions on axially symmetric open sets $U$. For more detailed study of slice regular functions can be seen in \cite{Q1,Q2,Q3,Q4}.

	Let $\mathcal{U} = \{(u, v) \in \mathbb{R}^2 : u + \mathbb{S}v \subset U\}$. Let  $f : U \to \mathbb{H}$ be a function  of the form
	$f(p) = f_0(u, v) + j f_1(u, v)$ for $p = u + jv \in U$	where $f_0, f_1 : \mathcal{U} \to \mathbb{H}$ that satisfy 
	\begin{equation*}\label{cont}
		f_0(u, -v) = f_0(u, v), \quad f_1(u, -v) = -f_1(u, v). \tag{2.4}
	\end{equation*}
	is called left slice regular function if  $f_0$ and $f_1$ satisfy the Cauchy-Riemann equations
	\begin{align*}
		\frac{\partial}{\partial u} f_0(u, v) - \frac{\partial}{\partial v} f_1(u, v) &= 0, \tag{2.5} \\
		\frac{\partial}{\partial v} f_0(u, v) + \frac{\partial}{\partial u} f_1(u, v) &= 0, \tag{2.6}
	\end{align*}

 In similar way we can define right slice regular functions.
The function which satisfy \eqref{cont} is called to be  left  slice function.
If $f$ is a left (or right)  slice function such that $f_0$ and $f_1$ are real-valued, then $f$ is called intrinsic. If the slice function $f$ is continuous it is known to be continuous slice function.  For more details \cite{Q1,Q2,Q3,Q4}.

Now, we give the characterization of Continuous Slice function by Colombo,  Gantner, and  Kimsey \cite{Q4}.
\begin{lem}\label{contslice}\cite{Q4}
	Let $\mathcal{U} \subset \mathbb{H}$ be axially symmetric and let
	\[
	\mathcal{U}^+ := \{(u, v) \in \mathbb{R} \times [0, +\infty) : u + \mathbb{S}v \subset \mathcal{U}\}.
	\]
	A function $f : \mathcal{U} \to \mathbb{H}$ is a continuous left slice function if and only if there exist continuous functions $f_0, f_1 : \mathcal{U}^+ \to \mathbb{H}$, where $f_1(u, v) = 0$ if $v = 0$, such that $	f(q) = f_0(u, v) + jf_1(u, v)$	for all $p = u + Jv \in \mathcal{U}$ with $v \geq 0$.
\end{lem}
The class of all left slice regular functions will be denoted by $\mathcal{SR}_L(U)$, right slice regular functions by $\mathcal{SR}_R(U)$ and class of all intrinsic functions by $\mathcal N(U)$ on $U$ respectively. The set $\mathcal{SR}_L(U)$ is a right linear space over $\mathbb{H}$.

For, an intrinsic function $f$ on an axially symmetric set $U$ the following are equivalent:
\begin{enumerate}
	\item[\rm(i)] $f$ is intrinsic.
	\item[\rm(ii)] $f(U \cap \mathbb{C}_I) \subset \mathbb{C}_I$ for every $I \in \mathbb{S}$.
	\item[\rm(iii)]  $f(\bar{p}) = \overline{f(p)}$ for all $p \in U$.
\end{enumerate}
All polynomials with quaternionic coefficients on the right are  left slice regular functions and quaternionic coefficients on the left are right slice regular functions.
Multiplication and composition of slice regular with intrinsic functions preserve the slice regularity but not true for arbitrary slice functions.

The slice derivative of left slice regular function $f$   denoted by $\partial_s f$ of $f$,  is 
	\[\partial_s(f)(p) = 
	\begin{cases}	\partial_I(f)(q) = \frac{1}{2}\left(\frac{\partial}{\partial u} - I \frac{\partial}{\partial v}\right)f_I(u + Iv), & \text{if } p = u + Jv, \ v \neq 0, \\[10pt]
		\dfrac{\partial f}{\partial u}(u), & \text{if } p = u \in \mathbb{R}.
	\end{cases}
	\]

The definition is well posed since, $\partial_I f = \partial_J f$ in $\Omega_I \cap \Omega_J$ for any choice of $I, J \in \mathbb{S}$.
\subsection{Properties of Slice Regular Functions}
\hfill

The following result(splitting lemma) states that if we restrict a slice regular function on an axially symmetric open set $U$ to one of the complex planes $C_J$, then the restricted function   is holomorphic. 
\begin{lem}\cite{splitting}\label{splitting lemma}
	Let $I, J \in \mathbb{S}$ with $I \perp J$ and if $f \in \mathcal{SR}_L(U)$, then the restriction $f_I = f|_{U \cap \mathbb{C}_I}$ satisfies
	\[f_I(z) = f_1(z) + f_2(z)J\]
	where $f_1, f_2$ are complex holomorphic functions on $U \cap \mathbb{C}_I$.
\end{lem}
By the Splitting Lemma, any slice regular function $f$ on an axially symmetric open set have a power series expansion for any real point on $U$, especially for a ball centered at a real origin  can be locally represented on this plane via a power series in $z$. 
But this  cannot be extended to balls centered at arbitrary points in the quaternionic space $\mathbb{H}$ since even $f(p) = (p - q)^n$ is not slice regular slice regular unless the center $q$ is  real.

	A quaternionic-function $f $ on a unit ball $\mathbb{B}$ is left slice regular if, and only if, it has a series expansion of the form
	\[	f(p) = \sum_{n \geq 0} p^n \frac{1}{n!} \frac{\partial^n f}{\partial x^n}(0)\]
	converging on $\mathbb{B}$. In particular if $f$ is slice regular, then it is infinitely differentiable on $\mathbb{B}$. In case of slice domain $U$ it can be extended for slice regular function $f$:

	For  $p_0 \in \mathbb{R}$ in $U$, the function $f$ can be represented by power series
	\[	f(p) = \sum_{n \geq 0} (p - p_0)^n \frac{1}{n!} \frac{\partial^n f}{\partial x^n}(p_0)\]
	on the ball $\mathbb B(p_0, R)$ where $R = R_{p_0}$ is the largest positive real number such that $\mathbb B(p_0, R)$ is contained in $U$.

 Using power series expansions, several analogues of well-known results from the theory of holomorphic functions have been established for slice regular functions. The identity principle holds  on all axially symmetric slice domains $U$.
 \begin{thm}\cite{poles}\label{identity principle}
 	 If, $f$ and $g$, regular functions on  $U$ coincide on a subset of $U\cap \mathbb{C}_I$ for some $I \in \mathbb{S}$,  having an accumulation point in $U\cap \mathbb{C}_I$, then $f = g$ in $U$.
 \end{thm}
 The following Representation formula is one of the important property of slice functions:
 \begin{thm}\cite{repnslice}\label{representation formula}
 	   A function $f \in \mathcal{SR_L}(U)$  on an axially symmetric set $U \subseteq \mathbb{H}$ is a left slice function on $U$ if and only if for every $p = u + Jv \in U$ we have
 	\begin{equation*}\label{extension}
 		f(p) = \frac{1}{2}\Big[f(\bar{z}) + f(z)\Big] + \frac{1}{2}JI\Big[f(\bar{z}) - f(z)\Big] 
 	\end{equation*}
 	with $z = u + Iv$ for $I \in \mathbb{S}$.
 \end{thm}
 As result of  Representation formula, every holomorphic function that is defined on a suitable open set in $\mathbb{C}_J$ has a slice regular extension.
 
 	Let $U_J$ be open subset of $ \mathbb{C}_j$ and symmetric  with respect to the real axis and $O = \bigcup_{x+Jy \in U_J, I\in \mathbb{H}} x+Iy$. Then,
 	\begin{enumerate}
 		\item Every quaternionic function $f$ on $U_J$ has a unique extension  to a left slice function on $O$ which is given by \eqref{extension}.
		\item If $f$ is left holomorphic, i.e., it satisfies \eqref{extension}, then its extension is left slice regular \cite{repnslice}.
 	\end{enumerate}
 
 The product  of two slice regular functions is not slice regular unless the factor on the appropriate side is intrinsic. But, there exists a regularized product that preserves slice regularity:
 
 	For, $f = f_0 + jf_1, g = g_0 + jg_1 \in \mathcal{SR}_L(U)$, then  their left slice regular product is
 	\[ 	f *_L g = (f_0g_0 - f_1g_1) + j(f_0g_1 + f_1g_0).\]
  
 	The slice regular product is associative and distributive, but  in general it is not commutative. If $f$ is intrinsic, then $f *_L g = fg$  the pointwise product $fg$ and
 	\begin{equation*}
 		f *_L g = fg = g *_L f. 
 	\end{equation*}
 
 	If $f(p) = \sum_{n=0}^{+\infty} p^n a_n$ and $g(p) = \sum_{n=0}^{+\infty} p^n b_n$ are two left slice regular power series, then their slice regular product is
 	\begin{equation*}
 		\left( \sum_{n=0}^{+\infty} q^n a_n \right) *_L \left( \sum_{n=0}^{+\infty} p^n b_n \right) = (f *_L g)(p) = \sum_{n=0}^{+\infty} p^n \sum_{k=0}^{n} a_k b_{n-k}. 
 	\end{equation*}
 	the usual product of formal power series with coefficients in a non-commutative ring.
 
 Similarly as the composition of slice regular functions is not, in general, slice regular, and so $f(p) = (p + a)^n$ is not a slice regular function, unless $a \in \mathbb{R}$. But, for any $a \in \mathbb{H}$ we can define
 \[ g(p) = (p + a)^{\star n} = \sum_{k=0}^{n} \binom{n}{k} p^{n-k} a^k\]
 which is  slice regular.
 
 	For $f = f_0 + jf_1 \in \mathcal{SR}_L(U)$, we define its slice regular conjugate $f^c = \overline{f_0} + j\overline{f_1}$ and its symmetrization $f^s = f *_L f^c = f^c *_L f$. 
 
 Now if, $f = f_0 + jf_1$  then,
 \[ f^s = |f_0|^2 - |f_1|^2 + j 2\text{Re}\,(f_0 \overline{f_1}).\]
 Hence it is an intrinsic function. It is $f^s(p) = 0$ if and only if $f(\tilde{p}) = 0$ for some $\tilde{p} \in [p]$. Also,
 \begin{equation*}
 	f^c(p) =  \overline{f(\overline{p})}, 
 	 \end{equation*}
 and $ 	f *_L g(p) = f(p) g\left( f(p)^{-1} p f(p) \right) \quad \text{if } f(p) \neq 0. $ 
 	The slice regular inverse of $f \in \mathcal{SR}_L(U)$ with $f \not\equiv 0$ which is denoted by $f^{-*_L}$, which satisfies $f^{-*_L} *_L f = f *_L f^{-*_L} = 1$, is given by
 		\[	f^{-*_L} = (f^s)^{-1} *_L f^c = (f^s)^{-1} f^c,\]
 		and it is defined on $U \setminus [Z(f)]$, where $Z(f) = \{s \in U : f(s) = 0\}$. See \cite{Q3,Q1,Q2,Q4,repnslice} for more details.
 
Considering the series expansion of the complex Cauchy kernel and  its closed form under the assumption that $s$ and $p$ are quaternions that do not commute, the left  slice regular Cauchy kernels has been defined:

	The left slice regular Cauchy kernel is defined as
	\[S_L^{-1}(s, p) = -(p^2 - 2\text{Re}(s)p + |s|^2)^{-1}(p - \bar{s}), \quad p \notin [s],
	\]

 The left slice regular functions satisfy Cauchy formulas with  left slice regular Cauchy kernel.
\begin{thm}\cite{Q1,Q2}
	Let $U \subset \mathbb{H}$ be a bounded axially symmetric set such that $U\cap \mathbb{C}_I$ is a Cauchy domain in $\mathbb{C}_I$ for every $J \in \mathbb{S}$,  and set $ds_I = ds(-I)$. If $f$ is a left slice regular function on a set that contains $\overline{U}$, then
	\begin{equation*}
		f(p) = \frac{1}{2\pi} \int_{\partial(U \cap \mathbb{C}_I)} S_L^{-1}(s, p) \, ds_I \, f(s), \quad \text{for every } p \in U. 
	\end{equation*}
	which  depend  neither on $U$ nor on the imaginary unit $J\in \mathbb{S}.$ Moreover,
	\[	f^{(n)}(p) = \frac{n!}{2\pi} \int_{\partial(U \cap \mathbb{C}_I)} (p^2 - 2s_0p + |s|^2)^{-n-1}(p - \bar{s})^{(n+1)*} \, ds_I \, f(s),\]
	where
	\begin{equation*}
		(p - \bar{s})^{n*} = \sum_{k=0}^{n} \frac{n!}{(n - k)!k!} p^{n-k} \bar{s}^k. 
	\end{equation*}
\end{thm}
Similarly for  unbounded  axially symmetric set $U \subset \mathbb{H}$ such that $U\cap \mathbb{C}_J$ is a Cauchy domain in $\mathbb{C}_J$ for every $J \in \mathbb{S}$, Cauchy formula hold.
 
 Let $J \in \mathbb{S}$. If $f \in \mathcal{SR}_L(\overline{U})$ and $f(\infty) = \lim_{|p|\to\infty} f(p)$ exists, then
\[f(q) = f(\infty) + \frac{1}{2\pi} \int_{\partial(U \cap \mathbb{C}_j)} S_L^{-1}(s, p) \, ds_j \, f(s) \quad \text{for every } p \in U.
\]

The quaternionic rational functions is defined this way:
\begin{defn}\cite{Q4}
	Let $P \in \mathcal{N}(\mathbb{H})$ and $Q \in \mathcal{SH}_L(\mathbb{H})$ be polynomials. A function $f$ is called left rational if it is of the form,
	\[ f(q) = P(q)^{-1} Q(q).\]
\end{defn}
It is necessary that $P$ is intrinsic otherwise $P^{-1}$ is not slice regular.  For any rational left slice regular function is $f(q) = P(q)^{-*_L} *_L Q(p)$ with left slice regular (but not necessarily intrinsic) polynomials, it can be represented in the above form $f(q) = \tilde{P}(q)^{-1}\tilde{Q}(q)$ with an intrinsic polynomial $\tilde{P}$ and a left slice regular polynomial $\tilde{Q}$. 

\begin{cor}\cite{Q3}
	Let $f \in \mathcal{SR}_L(U)$,  $f_I = f_1 + f_2J$ with holomorphic components $f_1, f_2 : U \cap \mathbb{C}_I \to \mathbb{C}_I$, with $J, I \in \mathbb{S}$, $J \perp I$  according to Lemma \ref{splitting lemma}. Then $f$ is left rational if and only if $F_1$ and $F_2$ are rational functions on $\mathbb{C}_J$.
\end{cor}
The slice regular function satisfy the quaternionic analogue of Schwarz Lemma:
\begin{thm}\cite{gentili}
	Let $f : \mathbb{B} \to \mathbb{B}$ be a regular function satisfying $f(0) = 0$. Then, $	|f(q)| \le |q|$ for all $q \in \mathbb{B}$
	and $|\partial_sf(0)| \le 1$.
\end{thm}
An analog of the complex Maximum Modulus Principle holds  for slice regular functions and  proved using the Splitting Lemma.
\begin{thm}\label{maximum modulus}\cite{poles}
	 Let $U$ be a slice domain and $f$ be slice regular quaternion valued function on $U$which has a relative maximum at $p$. Then $f$ is constant.
\end{thm}
\section{Laurent series and Singularities of regular Function}
Quaternionic  Laurent series and classification of singularities are presented in this section. For more details readers can refer \cite{poles, laurent}.
 
A regular Laurent series centered at 0 is an expression of the form:
\begin{thm}\cite{poles}
	For $\Omega = A(0, R_1, R_2) = \{q \in \mathbb{H} : R_1 < |q| < R_2\},$ a quaternionic spherical shell, where $0 \le R_1 < R_2 \le +\infty$ and if $f$ is a quaternionic-slice regular function on $\Omega $, then there exists a unique sequence of quaternionic coefficients $\{a_n\}_{n=-\infty}^{+\infty} \subset \mathbb{H}$ such that for all $p \in \mathbb{H}$ with $R_1 < |p| < R_2$, $		f(p) = \sum_{n=-\infty}^{+\infty} p^n a_n$.
	The series converge absolutely and uniformly on every compact subset of  $\Omega$.
\end{thm}
The  regular Laurent series centered at  general quaternion $p \in \mathbb{H}$ is also defined.  

For any $\{a_n\}_{n \in \mathbb{Z}}$ in $\mathbb{H}$,  the regular Laurent series centered at $q$ is
	\begin{equation*}
		\sum_{n \in \mathbb{Z}} (p - q)^{*n} a_n
	\end{equation*}

In order to study the convergence of Laurent series, the authors in \cite{poles} defined these distances:

For all $p, q \in \mathbb{H}$
	\begin{equation}
		\sigma(p,q) = \begin{cases} 
			
			|p - q| & \text{if } p, q \text{ lie on the same complex line } \mathbb{C}_I \\ 
			\omega(p,q) & \text{otherwise} 
		\end{cases}
	\end{equation}
	where $\omega(p,q) = \sqrt{[Re(p) - Re(q)]^2 + [|Im(p)| + |Im(q)|]^2}.$

 	For, $q \in \mathbb{C}_I \subset \mathbb{H}$ and $R \geq 0$, denote the $\sigma$-ball of radius $R$ centered at $p$ as
 	\begin{equation}
 		\Sigma(q, R) = \{p \in \mathbb{H} : \sigma(p,q) < R\}.
 	\end{equation}
 	Also, define $\Omega(q, R) = \{p \in \mathbb{H} : \omega(p,q) < R\}$.
 	For all $p, q \in \mathbb{H}$ 
 	\begin{equation}
 		\tau(p, q) = \begin{cases}
 			
 			&|p - q|,\  \text{if } p, q \text{ lie on the same complex line } \mathbb{C}_I \\
 			&\sqrt{[Re(q) - Re(p)]^2 + [|Im(q)| - |Im(p)|]^2},\ \text{otherwise}
 		\end{cases}
 	\end{equation}
The function $\tau$ is  not a distance, since $\tau(p,q) = 0$ for all $p$ lying in the same 2-sphere $x + y\mathbb{S}$ as $q$ except $\bar{q}$.

Chosen $q \in \mathbb{H}$ and $R \geq 0$, let  $T(q, R) = \{p \in \mathbb{H} : \tau(p,q) < R\}$. If $q = x + Iy$ then
	\[
	T(q, R) = \bigcup_{J \neq -I} \Delta_J(x + Jy, R).
	\]
The authors in \cite{poles} define the sets  of convergence of regular Laurent series. For $0 \le R_1 < R_2 \le +\infty$, 
		\begin{gather}
			\Sigma(q, R_1, R_2) = \{p \in \mathbb{H} : \tau(p, q) > R_1, \sigma(p, q) < R_2\}, \\
			\Omega(q, R_1, R_2) = \Omega(q, R_2) \setminus \overline{T(q, R_1)},
		\end{gather}
		where $\Omega(q, R_2) = \{p\in \mathbb{H} : \omega(p, q) < R_2\}$.
Here the interior of $\Sigma(q, R_1, R_2)$ is $\Omega(p, R_1, R_2)$ and
	\[
	\overline{\Sigma(q, R_1, R_2)} = \{p \in \mathbb{H} : \tau(p, q) \ge R_1, \sigma(p, q) \le R_2\}.
	\] 
	The convergence of Laurent series is explained as follows:	
	\begin{thm}\cite{poles}
		Let $\{a_n\}_{n \in \mathbb{Z}}$  be any sequence in $\mathbb{H}$ and $R_1, R_2 \in [0, +\infty]$ be such that $$R_1 = \limsup_{m \rightarrow +\infty} |a_{-m}|^{1/m},\  1/R_2 = \limsup_{n \rightarrow +\infty} |a_n|^{1/n}.$$ For all $q \in \mathbb{H}$ the series
		\begin{equation}
			f(p) = \sum_{n \in \mathbb{Z}} (p - q)^{*n} a_n
		\end{equation}
		converges absolutely and uniformly on the compact subsets of $\Sigma(q, R_1, R_2)$, and it does not converge at any point of $\mathbb{H} \setminus \overline{\Sigma(q, R_2)}$  and $T(q, R_1)$. If $\Omega(q, R_1, R_2) \neq \emptyset$ then the sum of the series defines a regular function $f : \Omega(q, R_1, R_2) \rightarrow \mathbb{H}$.
	\end{thm}
	For the expansion property define,  if $q \in \mathbb C_I$ and if $A_I(q, R_1, R_2)$ denotes the annulus $\{z \in \mathbb C_I : R_1 < |z - q| < R_2\}$ then,
	$	\Sigma(q, R_1, R_2) = A_I(q, R_1, R_2) \cup \Omega(q, R_1, R_2).$
	In this annulus region the existence of Laurent series is given as follows.
	\begin{thm}\cite{poles}
		Let $f$ be a slice regular function on a domain $U \subseteq \mathbb{H}$ and  $q \in \mathbb{H}$. There exists a sequence $\{a_n\}_{n \in \mathbb{Z}}$ in $\mathbb{H}$ such that, for all $0 \le R_1 < R_2 \le +\infty$ with $\Sigma(p, R_1, R_2) \subseteq U$,
		\begin{equation}
			f(p) = \sum_{n \in \mathbb{Z}} (p - q)^{*n} a_n
		\end{equation}
		in $\Sigma(q, R_1, R_2)$.
	\end{thm}
	It is possible to compute the coefficients of Laurent expansions by means of an integral formula. 
	
		If $p \in \mathbb C_I \subset \mathbb{H}$, $0 \le R_1 < R_2 \le +\infty$ and $f(p) = \sum_{n \in \mathbb{Z}} (p - q)^{*n} a_n$ in $\Sigma(q, R_1, R_2)$, then for each $n \in \mathbb{Z}$
		\begin{equation}
			a_n = \frac{1}{2\pi I} \int_{\gamma_I} ds (s - q)^{-n-1} f(s)
		\end{equation}
		where $\gamma_I(s) = q + Re^{2\pi I s}$ for any $R$ such that $R_1 < R < R_2$.

	 In \cite{laurent} authors, discussed another Laurent expansion, spherical Laurent series expansion which  are in terms of a second degree polynomial. 
	
		Let $q \in \mathbb{H}$. For any $\{a_n\}_{n \in \mathbb{Z}}$ in $\mathbb{H}$ the series
		\begin{equation}
			\sum_{n \in \mathbb{Z}} Q_p^n(q) [a_{2n} + (p - q)a_{2n+1}],
		\end{equation}
		is called the spherical Laurent series centered at $q$.
		
The convergence of this series is established in Cassini shell $$U(q, R_1, R_2) = \{p \in \mathbb{H} : R_1^2 < |(p - q_0)^2 + q_1^2| < R_2^2\}$$ where $q= q_0 + I q_1 \in \mathbb{H}$, with $q_0 \in \mathbb{R}$, $q_1 \in \mathbb{R}$, $I \in \mathbb{S}$ and $0 < R_1 < R_2$.  For $\{a_n\}_{n \in \mathbb{Z}} \subseteq \mathbb{H}$ we set
		\[
		R_1 = \limsup_{n \rightarrow \infty} |a_{-n}|^{\frac{1}{n}}, \quad \text{and} \quad \frac{1}{R_2} = \limsup_{n \rightarrow \infty} |a_n|^{\frac{1}{n}}.
		\]
		Then the spherical Laurent expansion  converges absolutely and uniformly on the compact subsets of the Cassini shell $U(q, R_1, R_2)$ \cite{laurent}.
\subsection{Singularities of Slice regular Function}
In this section we discuss the  classification of the singularities of regular functions. Readers can refer \cite{poles}.

	Let $f$ be a slice regular function on a axially symmetric slice domain $U$. A point $q \in \mathbb{H}$ is a singularity for $f$ if there exists $R > 0$ such that $\Sigma(q, 0, R) \subseteq U$ (so that the Laurent expansion of $f$ at $q$, $f(p) = \sum_{n \in \mathbb{Z}} (p - q)^{*n} a_n$, has $0$ as its inner radius of convergence, and it has a positive outer radius of convergence) \cite{poles}.

	For all $R > 0$ and for all $q = x + yI \in \mathbb{H}$, if $\Sigma(q, R)$ denotes the $\sigma$-ball,  then
	\begin{equation}
		\Sigma(q, 0, R) = \Sigma(q, R) \setminus [(x + y\mathbb{S}) \setminus \{\bar{q}\}] = [\Sigma(q, R) \setminus (x + y\mathbb{S})] \cup \{\bar{q}\}.
	\end{equation}
	In other words, $\Sigma(p, 0, R)$ consists of the punctured disk $\Delta_I(q, R) \setminus \{q\}$ and of the symmetric completion $\tilde{D}$ of $D = [\Delta_I(p, R) \setminus \{p\}] \cap [\Delta_I(\bar{p}, R) \setminus \{\bar{p}\}]$.\cite{poles}
	
	The singularities are classified in quaternionic as similarly in Complex. Let $q$ be a singularity for $f$. Then $q$ is a removable singularity if $f$ extends to a neighborhood of $q$ as a slice regular function. Otherwise, 
	$f(p) = \sum_{n \in \mathbb{Z}} (p - q)^{*n} a_n.$
	
	Also, $q$ is a pole for $f$ if there exists an $m \geq 0$ such that $a_{-k} = 0$ for all $k > m$; the minimum such $m$ is called the order of the pole and denoted by $ord_f(q)$. If $q$ is not a pole, then  it is an essential singularity for $f$ and  $ord_f(q) = +\infty$.

When $0 < R < |Im(p)|$, $\Sigma(q, 0, R)$ reduces to the punctured disk $\Delta_I(q, R) \setminus \{q\}$ (where $I \in \mathbb{S}$ is chosen so that $q \in \mathbb{C}_I$). Then,
\[f_I(z) = \sum_{n \in \mathbb{Z}} (z - q)^n a_n.\]
for all $z \in \mathbb{C}_I$. If $q$ is a pole for $f_I$, that is, if there exists $n \in \mathbb{N}$ such that
\[f_I(z) = (z - q)^{-n} g_I(z)\]
for some holomorphic $g_I : \Delta_I(q, R) \rightarrow \mathbb{H}$. On the other hand, $q$ is an essential singularity for $f_I$, that is, it is not a pole for $f_I$.  $q$ is a removable singularity for $f_I$, that is, $f_I$ extends as a holomorphic function to a neighborhood of $q$ in $\mathbb{C}_I$.\cite{poles} 

A point $q \in \mathbb C_I \subset \mathbb{H}$ is a singularity for a slice regular function $f$ on an axially symmetric slice domain if and only if it is an isolated singularity for the holomorphic function $f_I$. The function $f$ has a pole of order $n$ at $q$ if and only if $f_I$ does. Furthermore, $q$ is an essential singularity for $f$ if and only if it is an essential singularity for $f_I$. It does not hold for removable singularities, $q$ is a removable singularity for $f_I$ if and only if it is a pole of order $0$ for $f_I$ if and only if it is a pole of order $0$ for $f$, but such a pole need not be a removable singularity for $f$.\cite{poles}

\section{Analytic Capacity and Approximation Problems in Complex plane}
In this part the well versed approximation results in complex approximation theory has presented. The readers can refer \cite{analytic capacity,rational app} for more details. 
Runge's work in 1885 is generally regarded as the beginning of the systematic study of rational approximation and this theorem explains how the location of poles can be controlled by the connected components of the complement of the set on which approximation is required.
Weierstrass's classical theorem, published in the same year, established
uniform polynomial approximation for continuous functions on a compact
real interval. These results initiated two related but distinct
directions: approximation of continuous functions by polynomials and
approximation of holomorphic functions by rational functions.

\begin{thm}[Runge's Theorem]\cite{Runge}
	\label{thm:runge-rational}
	Let $E\subset\mathbb C$ be compact and let $f$ be holomorphic on an open neighborhood of $E$. Let $A\subset\widehat{\mathbb C}\setminus E$ contain at least one point from each connected component of
	$\widehat{\mathbb C}\setminus E$. Then, for every $\varepsilon>0$,
	there is a rational function $r$ whose poles belong to $A$ and such that
	\[	\sup_{z\in E}|f(z)-r(z)|<\varepsilon.\]
	In particular, if $\mathbb C\setminus E$ is connected, the pole may be 	placed at infinity, and $r$ may be chosen to be a polynomial.
\end{thm}
This result shows that rational approximation is available for every function holomorphic near a compact set, provided the poles are allowed in suitable components of the complement. 

The work of Walsh, Hartogs and Rosenthal, Lavrentiev, and Keldysh \cite{hartogrose,keldysh,laverentiev, walsh,walsh2}
 clarified the relationship between the topology of a
compact set and the class of functions that can be approximated on it.
Lavrentiev considered polynomial approximation on compact sets with
empty interior. His theorem identifies connectedness of
the complement as the exact condition for polynomial approximation.

\begin{thm} \cite{laverentiev}\label{lavrentiev}
	Let $E \subset \mathbb{C}$ be a compact set in the complex plane. 
	Every continuous function $f : E \to \mathbb{C}$ can be uniformly 
	approximated by holomorphic polynomials if and only if:
	\begin{enumerate}
		\item $E$ has no interior points (i.e., $E^\circ = \varnothing$).
		\item The complement $\mathbb{C} \setminus E$ is connected.
	\end{enumerate}
	That is, for any $f \in \mathcal C_{\mathbb C}(E)$ and any $\varepsilon > 0$, there exists 
	a polynomial $P(z)$ such that
	\[
	\sup_{z \in E} |f(z) - P(z)| < \varepsilon.
	\]
\end{thm}

Here $\mathcal C_{\mathbb C}(E)$ denotes the space of continuous complex-valued functions on $E$, and $\mathbb C[z]$ denotes the polynomials with complex coefficients. 

Keldysh obtained a related characterization for closures of plane
domains. In this case, the functions under consideration are continuous
on the closure and holomorphic in the domain.

\begin{thm}\cite{keldysh}\label{keldysh}
	If a closed domain $E$ does not divide the plane and $E=\bar{E^\circ}$, then every function continuous in $E$ and holomorphic inside $E$ admits uniform approximation by polynomials.
\end{thm}
 Mergelyan's theorem later unified these settings by allowing
an arbitrary compact set and requiring only continuity on the set and holomorphicity in its interior.

Mergelyan's theorem\cite{Mergelyan} is one of the central results in classical complex approximation theory. It gives a complete characterization of polynomial approximation for the algebra of continuous functions that are holomorphic in the interior of a compact plane set.

\begin{thm}\cite{Mergelyan}\label{Mergelyan}
	For a compact set $E$ in the complex plane, every function $f$ continuous on $E$ and holomorphic at its interior points can be uniformly approximated by polynomials, it is necessary and sufficient that the complement of $E$ be connected.
\end{thm}
 The theorem extends Runge's theorem, which starts with a function holomorphic on a neighbourhood, and incorporates the earlier works by Lavrentiev and Keldysh \cite{keldysh, laverentiev}.

Mergelyan  studied the rational approximation,
when the complement is not connected, polynomial approximation may fail because polynomials cannot reflect the topology of the holes. Rational functions with poles in the complement provide a more flexible class of approximants. 

For a compact set $E\subset\mathbb C$,  $\mathcal R_{\mathbb C}(E)$ is the uniform closure on $E$ of rational functions whose
poles lie outside $E$. Since every such rational function is
continuous on $E$ and holomorphic in its interior, one always has
$\mathcal R_{\mathbb C}(E)\subseteq \mathcal A_{\mathbb C}(E).$
A fundamental question was to determine when equality holds. If the complement has finitely many connected components, Mergelyan
approximation results imply that $\mathcal R_{\mathbb C}(E)=\mathcal A_{\mathbb C}(E)$. For arbitrary compact
sets, however, a complete description of rational approximation  lead  to capacity.

The distinction between $\mathcal R_{\mathbb C}(E)=\mathcal C_{\mathbb C}(E)$ and $\mathcal R_{\mathbb C}(K)=\mathcal A_{\mathbb C}(E)$ is important.
The first equality asked whether every continuous function on $K$ can be approximated rationally; the second asked whether every continuous function that is holomorphic in the interior can be approximated.
If $E^\circ=\varnothing$, then $\mathcal A_{\mathbb C}(E)=\mathcal C_{\mathbb C}(E)$, so the two questions
coincide. For sets with nonempty interior, they were different problems.

\subsection{Analytic capacity}
Let $\widehat{C}$ be the Riemann sphere and $\Omega_{\mathbb{C}}(E)$ be the unbounded component of $\widehat{C}\setminus E$ for a compact set $E \subset \mathbb{C}$.
The study of compact sets of analytic capacity zero goes back to
Painlev\'e's work on removable singularities. Ahlfors introduced the modern notion of analytic capacity in 1947, and the theory was subsequently developed by  Ahlfors and Beurling, Garabedian, Havinson, Ivanov,
Pommerenke,   and others\cite{ahlfors,  ahlforsbeurling, conway, Garabedian, havinsonmain, ivanov, pommer}. Analytic capacity measures how effectively a compact set can support bounded holomorphic
functions with a nonzero first-order term at infinity.

For a compact set $E\subset\mathbb C$, analytic capacity is defined by
$\gamma_{\mathbb C}(E) = \sup \{ |f'(\infty)| : f \in \mathcal{F}(E) \}$
where the family $\mathcal{F}(E)$ consists of all functions $f : \mathbb{C} \setminus E \to \mathbb{C}$ satisfying the following conditions:
\begin{enumerate}
	\item $f$ is analytic on  $\Omega_{\mathbb{C}}(E)$,
	\item $\| f\|_{\infty} \le 1$ ,
	\item $f(\infty) = 0$ admitting a Laurent expansion near infinity of the form
	$f(z) = \frac{a_1}{z} + \frac{a_2}{z^2} + \frac{a_3}{z^3} + \dots$ where $f'(\infty) = a_1$.
\end{enumerate}
and for an arbitrary set $S$ as $$\gamma_{\mathbb{C}}(S) = \sup_{\substack{E\subseteq S,\ E \text{ compact}}} \gamma_{\mathbb{C}}(E) ,$$ see \cite{analytic capacity}.
Equivalently, if
$f(z)=a_1/z+a_2/z^2+\cdots$ near infinity, then
$f'(\infty)$ is  identified with $a_1$ in the definition of analytic capacity. A function which satisfies the properties listed above is called admissible function for $E$ and for a compact set $E$, there exist a unique admissible function $f_a$ which satisfies $f_a'(\infty) = \gamma_{\mathbb{C}}(E)$ which is known to be the Planar Ahlfor's function for $E$ \cite{havinsonmain}. 

Analytic capacity measures the extent to which a set supports bounded holomorphic functions with a nonzero first-order coefficient at infinity.
\begin{thm}\cite{ahlfors}
	\label{ahlfors-removability}
	A compact set $E\subset\mathbb C$ is removable for bounded holomorphic functions if and only if $\gamma_{\mathbb{C}}(E)=0.$
\end{thm}

A compact set is removable here if every bounded function holomorphic on $\mathbb C\setminus E$ extends holomorphically across $E$. This result made analytic capacity a central object in complex function theory: it converts a question about extending holomorphic functions into a numerical value.

Analytic capacity also entered approximation theory through the work of Vitushkin. Vitushkin developed methods for relating rational approximation to local capacity estimates and for decomposing approximation problems into local pieces \cite{Vitushkin1, Vitushkin2, Vitushkin3}.

Continuous analytic capacity is a related quantity particularly suited to rational approximation. It is defined by the same extremal
principle as analytic capacity, but the admissible functions are
required to be continuous on the whole plane as well as holomorphic
off the set. Denote $\alpha_{\mathbb{C}}(E)$ for continuous analytic capacity defined as:
\[
\alpha_{\mathbb{C}}(E)=\sup\{|f'(\infty)|:
f\in \mathcal C_{\mathbb C}(\widehat{\mathbb{C}}),\
f\text{ is holomorphic on }  \widehat{\mathbb{C}}\setminus E,\
\|f\|_\infty\leq1,\ f(\infty)=0\}.
\]
The continuity requirement distinguishes $\alpha_{\mathbb C}$ from $\gamma_{\mathbb C}$ and is important when the approximation problem concerns arbitrary continuous functions.

Vitushkin established a capacity-theoretic criterion for uniform
rational approximation. In its usual form, it relates the equality $\mathcal R_{\mathbb C}(E)=\mathcal A_{\mathbb C}(E)$ to the behaviour of continuous analytic capacity on discs after removing $K$ and its interior.
The criterion expresses a global approximation property in terms of local capacity comparisons. 
\section{Scheme of Approximation}
Vitushkin's scheme is a fundamental method in the theory of rational approximation \cite{analytic capacity, rational app}. Its main idea is to localize a continuous function by means of smooth cutoff functions and then estimate the resulting localized functions using analytic capacity. This method converts a global approximation problem into a collection of local problems.
Let $E\subset\mathbb C$ be compact. The aim is to characterize when
$\mathcal R_{\mathbb C}(E)=\mathcal C_{\mathbb C}(E)$ or, more generally, when $\mathcal R_{\mathbb C}(E)=\mathcal A_{\mathbb C}(E)$ \cite{analytic capacity, rational app}.

	Let $\phi\in C_c^1(\mathbb C)$, and define
	\[
	T_\phi f(z)
	=
	\frac{1}{\pi}\int_{\mathbb C}
	\frac{f(\zeta)\,\overline{\partial}\phi(\zeta)}
	{z-\zeta}\,dA(\zeta),
	\]
	whenever the integral is defined, where
	\[
	\overline{\partial}
	=\frac12\left(\frac{\partial}{\partial x}
	+i\frac{\partial}{\partial y}\right)
	\]
	and $dA$ denotes planar area measure.
	Then $T_\phi f$ is the localization operator associated with
	the cutoff function $\phi$. Under the usual regularity assumptions,
	it satisfies the corresponding Cauchy--Green identity.
The localization principle is the basic analytic ingredient of
Vitushkin's method. 
\begin{lem}\cite{analytic capacity}
	Let $f$ be holomorphic on an open set $U\subset\mathbb C$,
	and let $\phi\in C_c^1(\mathbb C)$ satisfy $\phi=1$ on $U$.
	Then the localization operator can be chosen so that
	$f-T_\phi f$ is holomorphic on $U$.
\end{lem}
The Vitushkin scheme of Approximation is explained as follows:

	Let $\{Q_j\}$ be a locally finite covering of an open set
	$V\subset\mathbb C$ by squares. Under the standard hypotheses
	on the covering, there exists a smooth partition of unity
	$\{\phi_j\}$ subordinate to it such that
	\[
	0\leq\phi_j\leq1,\qquad
	\operatorname{supp}\phi_j\subset Q_j,\qquad
	\sum_j\phi_j=1\quad\text{on }V.
	\]
	If the side length of $Q_j$ is $\delta_j$, the functions may be
	chosen with derivative estimates of the form
	\[
	\|\nabla\phi_j\|_\infty\leq \frac{C}{\delta_j},
	\]
	where $C$ is independent of $j$.
	
This partition allows a function to be decomposed into local
pieces. In Vitushkin's method, the pieces are constructed using
localization operators rather than simply multiplying the original
function by the cutoff functions.

Let $\{\phi_j\}$ be a suitable smooth partition of unity. For a
continuous function $f$, define its localized components by
\[
f_j=T_{\phi_j}f.
\]
The Cauchy--Green representation and the identity
$\sum_j\phi_j=1$ yield the corresponding decomposition of $f$,
up to the normalization at infinity:
\[
f=f(\infty)+\sum_j f_j,
\]
whenever the representation and convergence hypotheses are
satisfied.

Each $f_j$ is associated with a localized region. 
Vitushkin's method relates local capacity estimates to uniform
rational approximation. The localization decomposition reduces
the global problem to estimating the contributions associated
with small regions. The central feature of the method is the passage from a global approximation problem to local estimates governed by analytic capacity. 

Davie's variable-radius construction refines the
covering step by allowing the sizes of the covering sets to
vary with the local geometry. Davie developed a variable-radius covering scheme that is useful in
the study of rational approximation \cite{davie}. Unlike the classical Vitushkin
method, which often uses localization on discs of a fixed scale,
Davie's construction allows the sizes of the covering sets to vary according to the local geometry. It is useful when
different parts of the domain require different scales of
localization. 
The main idea is to cover an open set $V\subset\mathbb{C}$ by a
suitable family of squares whose side lengths may vary. A smooth
partition of unity subordinate to this covering is then constructed.
The associated cutoff functions allow a function to be localized
into smaller pieces, which can be treated separately in the
approximation procedure.

\begin{prop}\cite{davie}
	\label{davie covering}
	Let $V$ be an open subset of $\mathbb{C}$, and 
	$p$ be a positive continuous function on $V$.
	Then there exists a countable family of squares
	$Q_r=Q(z_r,\delta_r),\qquad r\geq 1,$
	with centres $z_r\in V$ and side lengths $\delta_r>0$,
	satisfying the following properties:
	\begin{enumerate}
		\item The squares cover $V$:
	
		\item The covering has bounded overlap: every point of
		$V$ belongs to at most four squares.
		
		\item For each $r$, the square centered at $z_r$ with
		side length $\delta_r/3$ does not intersect any
		$Q_s$ with $s\neq r$.
		
		\item The side lengths are locally controlled by $p$:
		
		\item There exists a family of continuously differentiable
		functions $\{\phi_r\}_{r\geq1}$ subordinate to the covering
		such that
		\[
		0\leq\phi_r\leq1,\qquad
		\operatorname{supp}\phi_r\subset Q_r,
		\]
		and
		\[
		\sum_{r\geq1}\phi_r(z)=1,
		\qquad z\in V,
		\]
		with the derivative estimate
		\[
		\|\nabla\phi_r\|_\infty
		\leq \frac{100}{\delta_r}.
		\]
		
		\item There exists an absolute constant $A>0$ such that
		\[
		\sum_{r\geq1}
		\min\left\{1,
		\frac{\delta_r^3}{|z-z_r|^3}\right\}
		\leq A,
		\qquad z\in\mathbb{C},
		\]
		where a summand is interpreted as $1$ when $z=z_r$.
	\end{enumerate}
\end{prop}
The variable-radius construction provides a useful model for
developing localization methods in quaternionic approximation.
\part*{\MakeUppercase{Part I: Slice Regular Capacity}}
\vspace{0.5em}
\vspace{0.5em}
The set of slice regular functions contains polynomials and convergent power series of the quaternionic variable, provided that the coefficients are written  either on the left or the right. We are considering slice regular sets of functions in quaternionic setting which provides a natural extension of the class of holomorphic functions in the complex plane. 

In this paper we consider the class of left slice regular functions ($\mathcal{SR}_L(U)$). Slice regular functions possess interesting properties on specific open sets such as  axially symmetric sets and  slice domains. Slice regular functions satisfy the representation formula on these domains, which allows to reconstruct the values of the function once that we know its values on a complex plane $\mathbb{C}_I$. We denote $\widehat{\mathbb H} =\mathbb{H} \cup \{\infty\}$, one-point compactification  of $\mathbb{H}$.

A function $f$ is left slice regular at infinity if it is left slice regular on the set $\{q \in \mathbb{H} : r < |q|\}$ for some $r > 0$ and  $\lim_{q \to \infty} f(q)$ exists in $\mathbb{H}$. We define
\[f(\infty) = \lim_{q \to \infty} f(q).\] 

In this part we define slice regular capacity and study some of its geometric properties. Also we studied quaternionic analogue of $\gamma_{+}$-capacity and Cauchy capacity. The relationship of these capacities with removable singularity of slice regular functions are also discussed.
\section{Slice regular capacity}
Let $U \subseteq \mathbb{H}$ be axially symmetric open set and $f$ be a slice regular function on $U\subseteq \mathbb{\widehat{H}}$ satisfying $f(\infty)=0.$  Then $f$ admits a quaternionic Laurent expansion of the form
\[f(q) = {q^{-1}}{a_1} + {q^{-2}}{a_2} + q^{-3}{a_3} +\cdots.\]
and hence $\partial_s f(\infty)= a_1$. Also, $\partial_s f(\infty) = lim_{q\to\infty} q(f(q)-f(\infty))$ \cite{laurent}. For a compact axially symmetric set $K\subset\mathbb H$, we denote by $\Omega_{\mathbb{H}}(K)$ the  component containing $\infty$ of $\widehat{\mathbb H}\setminus K.$  Since $K$ is axially symmetric it is clear that,  the component $\Omega_{\mathbb{H}}(K)$ containing $\infty$ is also axially symmetric. This motivates the following definitions of slice regular capacity and continuous slice regular capacity as quaternionic analogous of analytic capacity and continuous analytic capacity.

\begin{defn}
	Let $K\subset\mathbb H$ be a compact axially symmetric set .
	The slice regular capacity of $K$ is defined by
	\[
		\gamma_{\mathbb H}(K) :=\sup_f |\partial_s f(\infty)|,
	\]
	where the supremum is taken over all left slice-regular functions
	$f:\Omega_{\mathbb H}(K)\longrightarrow\mathbb H$ satisfying $f(\infty)=0,$ and $\|f\|_{\infty,\Omega(K)}\leq 1.$
\end{defn}
For an arbitrary set $S\subseteq\mathbb H$, we define,
$$\gamma_{\mathbb H}(S):= \sup_{\substack{K\subseteq S,\ K \text{ compact and axially symmetric}}} \gamma_{\mathbb H}(K).$$
We consider the functions $f:\Omega_{\mathbb H}(K)\longrightarrow\mathbb H$ which is left slice regular satisfying $f(\infty)=0,$ and $\|f\|_{\infty,\Omega(K)}\leq 1.$ We call it as admissible functions satisfying  definition of slice regular capacity.
\begin{defn}
	Let $S\subset\mathbb H$. Then the quaternionic  continuous slice regular capacity of $S$ is defined by
	\[
	\alpha_{\mathbb H}(S): = \sup_{f}|\partial_s f(\infty)|.\] 
	where supremum is taken over all  functions $f$ satisfying:
	
	\begin{enumerate}
		\item $f\in \mathcal C_L(\mathbb{\widehat{H}})$ and $f(\infty)=0$;
		\item $\|f\|_{\infty,\Omega_{\mathbb H}(K)}\leq 1$;
		\item there exists a compact axially symmetric set $K\subset S$ such that
		$f\in \mathcal{SR}_L(\widehat{\mathbb H} \setminus K).$
	\end{enumerate}
	\end{defn}
	From the above definitions it is  clear that $\alpha_{\mathbb H}(S)\leq\gamma_{\mathbb H}(S)$. 
	
	\begin{exa}
		Let $\mathbb B(0,1)$ be the unit ball in $\mathbb{H}$ and  $S =K=\overline{\mathbb B}(0,1)$. Then $K$ is compact and axially symmetric. Then $\gamma_\mathbb{H}(S) = 1$. 
		For, let  $f(q)=q^{-1}$,  $f(0) = 0$. Then $f$ is slice regular on $\{|q|>1\}$ and at $\infty$. Hence $f$ is left slice regular on $\Omega(K) = \{|q|>1\}\cup\{\infty\}$.  Also, $f(\infty) = 0$. Since $|f(q)|=1/|q|<1$ for $|q|>1$, we have $\|f\|_{\infty, \Omega(K)}\leq 1$. Moreover, $\partial_s f(\infty)=\lim_{q\to\infty}q\,(q^{-1})=1 .$
		Hence $\gamma_\mathbb H(K)\ge 1$.
		Let $f$ be admissible for $K$. Then $g(q)=f(1/q)$,
		$g(0)=0$, is left slice regular on $B(0,1)$ with $|g|\le1$. Then by Quaternionic Schwarz lemma, $$|\partial_s f(\infty)|=|\partial_s g(0)|\le 1.$$
		Taking the supremum, $\gamma_\mathbb H(K)\le 1$.
		
		Also, we can see that $\alpha_\mathbb H(S)=1$. 
		For, since $\alpha_{\mathbb H}(S)\leq\gamma_{\mathbb H}(S) = 1$ we have $\alpha_\mathbb H(S)\leq1$. Let,
		\[
		f(q)=\begin{cases} q^{-1}, & |q|> 1,\\ \bar q, & |q|\leq 1,\end{cases}
		\qquad f(\infty)=0 .
		\]
		Then, on $|q|=1$ we have $q\bar q=1$, so $\bar q=q^{-1}$ and the two pieces agree. Each piece is continuous and $f(q)\to0$ as $q\to\infty$. Hence $f\in \mathcal C_L(\widehat{\mathbb H})$ and $f(\infty)=0$.
		Since, $|f(q)|=1/|q|\leq 1$ for $|q|> 1$ and $|f(q)|=|q|\le1$ for
		$|q|\le 1$, we have $\|f\|_{\infty,\Omega(K)}\le1$.
		Also, for $\{|q|>1\}$, $f(q)=q^{-1}$ is left slice regular, and near $\infty$, $f(1/q)=q$ is left slice regular at $0$ with value $0$. Hence $f\in\mathcal{SR}_L(\widehat{\mathbb H}\setminus K)$. 	Thus. $f$ is admissible for $\alpha_\mathbb{H}(S)$, and $\partial_s f(\infty)=\lim_{q\to\infty}q\,(q^{-1})=1$. Thus, $\alpha_\mathbb{H}(S)\ge 1$. Hence $\alpha_\mathbb{H}(S) = 1$
	\end{exa}
In contrast to classical complex analysis, the algebraic properties governing slice regular functions strictly depend upon restricting the parameters $a$ and $r$ to the real field $\mathbb{R}$. Since real translations and real homogeneity preserve  the underlying axial symmetry, the class of admissible slice regular functions remains invariant under the transformation $q \mapsto rq - a$, where $r, a \in \mathbb{R}$. Conversely, a general quaternionic transformation $q \mapsto tq - b$ with non-real parameters  doesn't preserve the individual complex slices and also the standard invariance arguments within the slice-regular framework.

		\begin{prop}\label{subsetmonogamma}
		\begin{enumerate}
			\item If $S_1\subseteq S_2,$ then $\gamma_{\mathbb H}(S_1)\leq \gamma_{\mathbb H}(S_2).$
			\item Let $a\in\mathbb R$ and $r\in\mathbb R\setminus\{0\}$. Then, ${\gamma_{\mathbb H}(rS+a)=	|r|\gamma_{\mathbb H}(S).}$
		\end{enumerate}
	\end{prop}
	\begin{proof}
		(1) If $K\subseteq S_1$, then, since $S_1\subseteq S_2$,
		$K\subseteq S_2.$	Thus, $\gamma_{\mathbb H}(S_1)\leq\gamma_{\mathbb H}(S_2).$
		
		(2)  Since $a\in\mathbb R$, real translation preserves compactness and axial symmetry.  
	\end{proof}
	
	\begin{prop}\label{monotoniccap}
		\begin{enumerate}
			\item If $S_1\subseteq S_2,$ then $\alpha_{\mathbb H}(S_1)\leq \alpha_{\mathbb H}(S_2).$
			\item Let $a\in\mathbb R$ and $r\in\mathbb R\setminus\{0\}$. Then, ${\alpha_{\mathbb H}(rS+a)=	|r|\alpha_{\mathbb H}(S).}$
		\end{enumerate}
	\end{prop}
	\begin{proof}
		The required result follows in  a similar manner as in Proposition \ref{subsetmonogamma}.
	\end{proof}
	The axially symmetric domain  allows us to reconstruct the  the function after its values on a complex plane $\mathbb{C}_I$. Hence we have the following result which connects slice regular capacity to analytic capacity.
	
\begin{thm}\label{complexslice quaternionic}
	Let $S\subset\mathbb H$ be a axially symmetric set and for 
	$I\in\mathbb S$ such that $	S_I = S\cap\mathbb C_I.$ Then, 
	\begin{enumerate}
		\item  $\gamma_{\mathbb H}(S)=\gamma_{\mathbb C}(S_I).$
		\item $	\alpha_{\mathbb H}(S)= 	\alpha_{\mathbb C}(S_I).$		
	\end{enumerate}
\end{thm}
\begin{proof}
		(1) Let $K\subset\mathbb H$ be a compact axially symmetric subset of $S$ and for $I\in\mathbb S$ such that $	E =K\cap\mathbb C_I.$  Let $ z \in \Omega_{\mathbb {C}_I}(E)$. That is, $z$ belongs to the component containing infinity in the complex slice, there is a  path joining $z$ to $\infty$ that stays in $\widehat{\mathbb C}_I \setminus E = (\widehat{H}\setminus K) \cap \widehat{\mathbb C}_I$. Thus, it is contained in the connected component of
		$\widehat{\mathbb H}\setminus K$ containing $\infty$. Hence,  $\Omega_{\mathbb {C}_I}(E)\subset\widehat{\mathbb H}\setminus K.$
	  		
		Let $z\in\Omega_{\mathbb H}(K)\cap\mathbb C_I.$
		Since $\Omega_{\mathbb H}(K)$ is a connected open subset of the locally path-connected space $\widehat{\mathbb H}$, it is path connected. Thus there exists a
		continuous path $\gamma:[0,1]\longrightarrow		\widehat{\mathbb H}\setminus K$	such that $	\gamma(0)=z,\ \gamma(1)=\infty.$	For $t\in[0,1)$, write
		$$\gamma(t)=x(t)+y(t)J(t),$$ where $x(t) = \operatorname{Re}(\gamma(t)),\ 	y(t) =|\operatorname{Im}(\gamma(t))| \geq0,$ and $ J(t)\in\mathbb S$ whenever $y(t) >0$. Define the projection of $\gamma$ onto the fixed slice $\mathbb C_I$ by 
		$$\widetilde\gamma(t)	= 	x(t)+y(t)I$$ and set $\widetilde\gamma(1) = \infty.$
		
		 Since, $x(t)$ and $y(t)$  are continuous, 	the map $\widetilde\gamma$ is continuous on $[0,1)$. Moreover, $|\widetilde{\gamma}(t)|^2 = x(t) ^2 +y(t) ^2 = |\gamma(t)|^2$. Thus, $	\widetilde\gamma(t) \to \infty$ as $t\to1$ and hence $\widetilde{\gamma}(t)$ is continuous at $t=1$ in the one-point compactification.	
		
		We have to show that 	$\widetilde\gamma(t)\notin E$ for every $t\in[0,1].$ Suppose, to the contrary, that for some $t_0\in[0,1)$,
		$$\widetilde\gamma(t_0)	=	x(t_0)+y(t_0)I	\in E.$$	
		Since $E\subset K$ and $K$ is axially symmetric, the whole
		2-sphere $\left\{x(t_0)+y(t_0)J:J\in\mathbb S \right\}$
		is contained in $K$. In particular, $\gamma(t_0) = 		x(t_0)+y(t_0)J(t_0)	\in K,$	which contradicts $	\gamma([0,1])\subset\widehat{\mathbb H}\setminus K.$
		 Thus, $$\widetilde\gamma([0,1])	\subset	\widehat{\mathbb C}_I\setminus E.$$
		 Write $z=x_0 +y_0 I$, with $y_0 \in \mathbb{R}$. Then, 
		 \[\widetilde{\gamma}(0) =  x_0 +|y_0|I = \begin{cases}
		 	z, & y_o \geq 0\\
		 	\overline z, & y_0 <0.
		 \end{cases}\]
Thus either $z$ or $\overline{z}$ can be joined to $\infty$ 
		by a path in $\widehat{\mathbb C}_I\setminus E$.
		Since $K$ is axially symmetric, $E$ is invariant under complex conjugation. Thus, complex conjugation preserves the connected component $\Omega_{\mathbb{C}_I}(E)$ containing $\infty$. 
		It follows that $	z\in\Omega_{\mathbb C_I}(E).$
		Hence, $\Omega_{\mathbb H}(K)\cap\mathbb C_I\subseteq\Omega_{\mathbb{C}_I}(E).$ Therefore, $$\Omega_{\mathbb H}(K)\cap\mathbb C_I=\Omega_{\mathbb{C}_I}(E).$$
		
		Now, we show that $\gamma_{\mathbb H}(K)\leq\gamma_{\mathbb C}(E).$ For, let $f$ be an arbitrary admissible function in the definition of
		$\gamma_{\mathbb H}(K)$; that is, $f$ is left slice regular on
		$\Omega_{\mathbb{H}}(K)$, satisfies $f(\infty)=0$ and $\|f\|_{\infty}\leq 1.$ Consider the restriction $f_I = f|\Omega_{\mathbb{C}_I}(E).$ Choose $J\in\mathbb S$ such that $J\perp I$. Then by the Lemma \ref{splitting lemma}, there exist holomorphic functions $F,G$ on
		$\widehat{\mathbb C}_I\setminus E$ such that $$f_I(z)=F(z)+G(z)J.$$ Since $f(\infty)=0$, we have $f_I(\infty)=0,$	and hence $	F(\infty)=G(\infty)=0.$
		
		Since $F$ and $G$ are holomorphic near $\infty$, their Laurent expansions have the form $$F(z)=\frac{a}{z}+O\!\left(\frac1{z^2}\right),\  
		G(z)=\frac{b}{z}+O\!\left(\frac1{z^2}\right)$$ as $z\to\infty$. Hence, $$f_I(z)=\frac{a+bJ}{z}+O\!\left(\frac1{z^2}\right),$$ and $f_I'(\infty)=a+bJ.$   Thus $|f_I'(\infty)|=\sqrt{|a|^2+|b|^2}.$ If $f_I'(\infty)=0$, then the result holds.
		Suppose $f_I'(\infty)=a+bJ=c\neq0.$ Define for $\omega,\nu\in\mathbb C_I$, $g:\mathbb H\longrightarrow\mathbb C_I$ by $$g(\omega+\nu J) =	\frac{\overline a}{|c|}\,\omega +	\frac{\overline b}{|c|}\,\nu,$$ which is  $\mathbb C_I$-linear. Now, by the Cauchy-Schwarz inequality, we have
		$$	|g(\omega+\nu J)|\leq\frac{\sqrt{|a|^2+|b|^2}}{|c|}	\sqrt{|\omega|^2+|\nu|^2} = |\omega+\nu J|.$$
		Hence, $\|g\|\leq1$ and $$g(c) = \frac{\overline a\,a+\overline b\,b}{|c|} = |c|.$$ Now, for $z\in\Omega_{\mathbb C_I}(E)$ define $h(z)=g(f_I(z))$.
		Since $g$ is $\mathbb C_I$-linear and $f_I$ is holomorphic,
		$h$ is a $\mathbb C_I$-valued holomorphic function. Also, $$|h(z)|\leq|f_I(z)|\leq1.$$ Thus, $\|h\|_\infty\leq1.$  Also, $h(\infty)=0$
		and $$h'(\infty) = g(f_I'(\infty)) = g(c) = 	|c|	= |f_I'(\infty)|.$$
		Thus $h$ is an admissible function for  $\gamma_{\mathbb C}(E)$,
		and hence $|f_I'(\infty)|=|h'(\infty)|	\leq\gamma_{\mathbb C}(E).$ Also since $f_I'(\infty)=a+bJ$ from the Laurent series of $f$, we have $\partial_s f(\infty)=f_I'(\infty)=a+bJ$, Thus, $$|\partial_s f(\infty)| = |f_I'(\infty)|=|h'(\infty)|	\leq\gamma_{\mathbb C}(E).$$
		Now, by taking the supremum over all admissible $f$ we have $$\gamma_{\mathbb H}(K)\leq\gamma_{\mathbb C}(E).$$
		
		To, prove the reverse inequality, let $h$ be an arbitrary admissible holomorphic function in the definition of $\gamma_{\mathbb C}(E)$. Thus $g:\Omega_{\mathbb C_I}(E)\longrightarrow\mathbb C_I,$ such that $h(\infty)=0$ and $\|h\|_\infty\leq1.$
		Let its Laurent expansion at infinity be $$h(z)=\frac{a}{z}+O\!\left(\frac1{z^2}\right)$$ as $z\to\infty.$ Hence, $h'(\infty)=a.$ If $a=0$, there is nothing to prove. Suppose $a\neq0$ and define $\lambda=\frac{\overline a}{|a|}\in\mathbb C_I.$ Then, $|\lambda|=1,$ and $\lambda a=|a|$.	Since $E$ is conjugation symmetric define, $h^*(z) = \overline{g(\overline z)}.$ Then $h^*$ is holomorphic on
		$\Omega(E)$, $\|h^*\|_\infty\leq1,$ and $h^*(\infty)=0.$
		Also, $$h^*(z) = \frac{\overline a}{z} +	O\!\left(\frac1{z^2}\right),$$ as $z\to\infty.$	Define $$g(z) =	\frac12\left(\lambda h(z)+\overline{\lambda}\,h^*(z)\right).$$ Thus, $g(z) =	\frac12\left(\lambda h(z)+\overline{\lambda}\,\overline{h(\overline z)}\right).$
		Hence $g$ is holomorphic on $\Omega_{\mathbb{C}_I}(E)$ and satisfies
		$g(\overline z)=\overline{g(z)}.$ Thus, $g$ is intrinsic and $$|g(z)|	\leq	\frac12|\lambda|\,|h(z)| +	\frac12|\overline\lambda|\,|h^*(z)|\leq1$$ and hence we have, $\|g\|_\infty\leq1.$ Also, $g(\infty)=0$ and $g'(\infty)=|a|=|h'(\infty)|.$
		
		Since $g$ is intrinsic, for $z=x+Iy\in \Omega_{\mathbb C_I}(E)$, write
		$g(z)=u(x,y)+Iv(x,y),$ for $u,v$ real-valued function.	For any $J\in\mathbb S$, define $$f(x+Jy) =u(x,y)+Jv(x,y).$$ Since $h$ is holomorphic, we have $	u_x=v_y$ and $u_y=-v_x.$ Since 	$h(\overline z)=\overline{h(z)}$, we have $u(x, -y) = u(x,y)$ and $v(x, -y) = -v(x,y)$, it is clear that $f$ is well defined and independent of the choice of $J$.  Now, by the Cauchy-Riemann equations,  the function $f(x+Jy)=u(x,y)+Jv(x,y)$ is  slice regular on $\Omega_{\mathbb H} (K).$  Since, $$	|f(x+Jy)|	=\sqrt{u(x,y)^2+v(x,y)^2}= |g(x+Iy)|,$$ we have $	\|f\|_\infty=\|g\|_\infty\leq1.$ Also, $f(\infty)=0$ and from the Laurent expansion of $h$ we have $	\partial_s f(\infty)=g'(\infty).$
		Thus $f$ is admissible for $\gamma_{\mathbb H}(K)$ and $|\partial_s f(\infty)| = |g'(\infty)| = 	|h'(\infty)|.$
		Since $f$ is admissible for $\gamma_{\mathbb H}(K)$, $$|h'(\infty)| = |\partial_s f(\infty)|\leq	\gamma_{\mathbb H}(K).$$  Thus, by taking the supremum over all admissible $g$ gives the reversed inequality. Hence, $\gamma_{\mathbb H}(S)=\gamma_{\mathbb C}(S_I).$
		
		\noindent (2)  The required result  follows in a similar manner as (1).
\end{proof}

By the previous result,  for axially symmetric sets slice regular capacity is completely determined by the analytical capacity of any complex slice. Thus we have from the proof, for each $I\in\mathbb S$ and $S_I=S\cap\mathbb{C}_I,$
$$\gamma_{\mathbb H}(S)	=	\sup_{I\in\mathbb S}\gamma_{\mathbb C_I}(S_I)\ \text{and}\ \alpha_{\mathbb H}(S) =	\sup_{I\in\mathbb S}\gamma_{\mathbb C_I}(S_I)$$	
where $\gamma_{\mathbb H}, \alpha_{\mathbb H}$ denotes the slice regular capacity, continuous slice regular capacity  and
$\gamma_{\mathbb C_I}, \alpha_{\mathbb C_I}$ denotes the  analytic capacity,  continuous analytic capacity in the slice $\mathbb C_I$ respectively. 

In general, $\gamma_\mathbb{H}$ and $\alpha_\mathbb{H}$ need not be same.
\begin{exa}
	Let $S=\{q\in\mathbb H:|q|=1\}$ which  is a compact,  axially symmetric subset of $\mathbb H$.  We first show that
	$\gamma_{\mathbb H}(S)>0.$	Define
	\[f(q)=
	\begin{cases}
		0, & |q|\leq 1,\\[2mm]
		q^{-1}, & |q|>1
	\end{cases}
	\]
	The function $f$ is left slice regular on each component of
	$ \widehat{\mathbb H}\setminus S$ . Also,
	\[	\|f\|_{\infty, \Omega_{\mathbb H}(S)}\leq 1,\qquad f(\infty)=0.
	\]
	For $|q|>1$, the quaternionic Laurent expansion of $f$
	at infinity is $f(q)=\sum_{n=1}^{\infty}q^{-n}a_n,$
	where $a_1=1,\ a_n=0\quad(n\geq2).$
	Hence, the first coefficient  is nonzero, and
	thus $\gamma_{\mathbb H}(S)\geq1.$	

	We now prove that $\alpha_{\mathbb H}(S)=0.$
	Let $K\subset S$ be any compact axially symmetric set, and
	let $g$ be an admissible function in the definition of
	$\alpha_{\mathbb H}(S)$. Thus $g$ is continuous on
	$\widehat{\mathbb H}$ and left slice regular on
	$\mathbb H\setminus K$, with $g(\infty)=0.$	
	Fix $I\in\mathbb S$ and consider the restriction
$g_I=g|_{\mathbb C_I}.$ Then $g_I$ is continuous on the Riemann sphere
	$\widehat{\mathbb C}_I$ and holomorphic on
	$\mathbb C_I\setminus E$, where 	$$E=K\cap\mathbb C_I
	\subset\{z\in\mathbb C_I:|z|=1\}.$$
	
		In particular, $E$ has planar area zero. Hence by the 
	removability theorem for continuous functions holomorphic
	outside a set of planar area zero, $g_I$ extends
	holomorphically to the whole Riemann sphere $\widehat{\mathbb C}_I$. Thus $g_I$ is constant. Since $g_I(\infty)=0$, we have $g_I\equiv0.$	As $I\in\mathbb S$ was arbitrary, it follows that $	g\equiv0\quad\text{on }\widehat{\mathbb H}.$
	Thus every admissible function in the definition of
	$\alpha_{\mathbb H}(S)$ has zero coefficient at infinity.
	Hence, $\alpha_{\mathbb H}(S)=0.$ Thus,
		\[\alpha_{\mathbb H}(S)=0<1\leq\gamma_{\mathbb H}(S).\]
\end{exa}
 Equality holds if $S$ is open:
\begin{thm}
	Let $S\subset\mathbb H$ be open. Then $\gamma_{\mathbb H}(S)=\alpha_{\mathbb H}(S).$
\end{thm}

\begin{proof}
 First we prove this for  an open axially symmetric set.  Let $U$ be an open axially symmetric set and let $U_I=U\cap\mathbb C_I.$
	Since $U$ is axially symmetric, by Theorem \ref{complexslice quaternionic}
	$$\alpha_{\mathbb H}(U)=\alpha_\mathbb{C}(U_I)\ \text{and} \gamma_{\mathbb H}(U)=\gamma_\mathbb{C}(U_I).$$ Since $U$ is open we have $U_I$ is open. Thus by [\cite{analytic capacity}prop 10.4] we have  $\alpha_\mathbb{C}(U_I) = \gamma_\mathbb{C}(U_I)$. Hence $\alpha_{\mathbb H}(U) = \gamma_{\mathbb H}(U)$ for any  open axially symmetric set.
	
	 Let $K\subset S$ be an arbitrary compact axially symmetric set. Since $K$ is compact and $S$ is open, there exists $r>0$ such that
	$$K\subset U_r=\{q\in\mathbb H:\operatorname{dist}(q,K)<r\}\subset S.$$ We claim that $U_r$ is an open axially symmetric set. For, let $q=x+Iy\in U_r,$ for $I\in\mathbb S.$ Since $K$ is axially symmetric, the distance from $q$ to $K$ depends only on $x$ and $|y|$. Consequently, for every $J\in\mathbb S$,
	$	\operatorname{dist}(x+Jy,K) = 	\operatorname{dist}(x+Iy,K)<r.$ Thus
	$x+Jy\in U_r$ for every $J\in\mathbb S,$ and hence $U_r$ is axially symmetric.
	
	We already proved that for open axially symmetric sets,  $\gamma_{\mathbb H}(U_r)=\alpha_{\mathbb H}(U_r).$	Now, by using the monotonicity of $\gamma_{\mathbb H}$ and $\alpha_{\mathbb H}$, we have
	\begin{align*}
	\gamma_{\mathbb H}(K)	\leq 	\gamma_{\mathbb H}(U_r)
		=\alpha_{\mathbb H}(U_r)
	\leq	\alpha_{\mathbb H}(S).
	\end{align*}
	Since $K\subset S$ was an arbitrary compact axially symmetric set,
	taking the supremum over all such $K$ gives $$\gamma_{\mathbb H}(S) =
		\sup_{\substack{K\subset S\\
				K\text{ compact and axially symmetric}}}
		\gamma_{\mathbb H}(K) \leq \alpha_{\mathbb H}(S).$$
	Then by the definition, we have for any open set $\gamma_{\mathbb H}(S)=\alpha_{\mathbb H}(S).$

\end{proof}
 In connection with Vitushkin's problem on rational approximation, there was a conjecture whether analytic capacity  is semi-additive, that is  whether there exists, a universal constant $C$ such that $\gamma_{\mathbb C}(F\cup G) \leq C\bigl(\gamma_{\mathbb C}(F)+\gamma_{\mathbb C}(G)\bigr)$ for all compact sets $F,G\subset\mathbb C$. Tolsa eventually proved the semi-additivity of analytic capacity and continuous analytic capacity, in fact in the stronger countable form for borel sets $G_j$, $$\gamma_{\mathbb C}\left(\bigcup_j G_j\right)\leq C\sum_j\gamma_{\mathbb C}(G_j),\ \text{and}\
 \alpha_{\mathbb C}\left(\bigcup_j G_j\right)\leq C\sum_j\alpha_{\mathbb C}(G_j),$$ with a universal constant $C$ \cite{tolsa, tolsaalpha}. This result was an important development in the theory of analytic capacity and had significant consequences for rational approximation. 
Here, we establish the corresponding  countable semi-additivity property for the slice regular capacity and the continuous slice regular capacity of compact axially symmetric subsets of $\mathbb H$ using Theorem \ref{complexslice quaternionic}.
	\begin{thm}
		\label{thm:semiadditivity-quaternionic}
		Let $K\subset\mathbb H$ be a compact axially symmetric set, and let
		$\{K_j\}_{j=1}^{\infty}$ be a sequence of axially symmetric Borel sets such that $K=\bigcup_{j=1}^{\infty}K_j.$
		Then there exists an absolute constant $C>0$ such that
		\[\gamma_{\mathbb H}(K)\le C\sum_{j=1}^{\infty}\gamma_{\mathbb H}(K_j)	\]
		and
		\[\alpha_{\mathbb H}(K)\le C\sum_{j=1}^{\infty}\alpha_{\mathbb H}(K_j).\]
	\end{thm}
	
	\begin{proof}
		Fix $I\in\mathbb S$ and $E=K\cap\mathbb C_I,\ 
		E_i=K_i\cap\mathbb C_I,\ i\ge1.$
		Since $K=\bigcup_{j=1}^{\infty}K_j$, we have
		$E=\bigcup_{j=1}^{\infty}E_j.$
		Thus, by the semi-additivity theorems for analytic capacity and continuous analytic capacity , there exist absolute
		constants $C>0$ such that
		\[
		\gamma_{\mathbb C}(E)
		\le C\sum_{i=1}^{\infty}
		\gamma_{\mathbb C}(E_j)
		\]
		and
		\[
		\alpha_{\mathbb C}(E)
		\le C\sum_{i=1}^{\infty}
		\alpha_{\mathbb C}(E_j).
		\]
	Hence by Theorem \ref{complexslice quaternionic}, we have,
		\[
		\begin{aligned}
			\gamma_{\mathbb H}(K)
			&=\gamma_{\mathbb C_I}(E)\\
			&\le C_\gamma\sum_{j=1}^{\infty}
			\gamma_{\mathbb C_I}(E_j)=C_\gamma\sum_{j=1}^{\infty}
			\gamma_{\mathbb H}(K_j),
		\end{aligned}
		\]
		and 
		\[
		\begin{aligned}
			\alpha_{\mathbb H}(K)
			&=\alpha_{\mathbb C_I}(E)\\
			&\le C_\alpha\sum_{j=1}^{\infty}
			\alpha_{\mathbb C_I}(E_j)=C_\alpha\sum_{j=1}^{\infty}
			\alpha_{\mathbb H}(K_j).
		\end{aligned}
		\]
		
	\end{proof}

In complex plane case, there is a qualitative relation connecting analytic capacity and diameter of a compact connected set which provides a measure of geometric size of the set. 
\begin{thm}\cite{rational app}\label{compdiam} 
	If $E\subset \mathbb C$ is compact and connected , then
	$\gamma(E) \leq \operatorname{diam}(E) \leq 4\gamma(E).$
\end{thm}
Since for a axially symmetric set, Theorem \ref{complexslice quaternionic} allows the complex estimate to be transferred to a complex slice it remains to relate diameter  of a compact subset of quaternionic set to the diameter of its slice. 
\begin{thm}\label{diameterac}
	Let $K\subset\mathbb H$ be a compact axially symmetric slice connected
	set. Then
	\[
	\gamma_{\mathbb H}(K)
	\leq \operatorname{diam}_{\mathbb H}(K)
	\leq 4\gamma_{\mathbb H}(K),
	\]
	where, $	\operatorname{diam}_{\mathbb H}(K) =\sup_{q,p\in K}|q-p|.$
\end{thm}

\begin{proof}
	Let $E=K\cap\mathbb C_I.$	Since $K$ is slice connected, $E$ is connected. Hence $E$ is a compact connected subset of $\mathbb C_I$.  Hence by Theorem \ref{complexslice quaternionic} and \ref{compdiam} we have $$\gamma_{\mathbb H}(K)\leq\operatorname{diam}_{\mathbb C_I}(E)	\leq 4\gamma_{\mathbb H}(K).$$
	
	Now, it is enough to prove that $\operatorname{diam}_{\mathbb H}(K) =
	\operatorname{diam}_{\mathbb C_I}(E).$ Since $E\subseteq K$, we immediately have $\operatorname{diam}_{\mathbb C_I}(E) \leq \operatorname{diam}_{\mathbb H}(K).$ For the converse, let $q=x+yJ,\ p=u+vL$ be arbitrary points of $K$, where $J,L\in\mathbb S$ and $y,v\geq0$.
	 Now, $$|q-p|^2	\leq (x-u)^2+y^2+v^2+2yv = (x-u)^2+(y+v)^2.$$
	Since $K$ is axially symmetric, and $x+yJ,\ u+vL \in K$ we have $x+yI, u-vI\in E$. Also, $\left|(x+yI)-(u-vI)\right|^2	=	(x-u)^2+(y+v)^2.$
	Hence,	\[
	|q-p|
	\leq
	\left|(x+yI)-(u-vI)\right|
	\leq
	\operatorname{diam}_{\mathbb C_I}(E).
	\]
	Thus, by taking the supremum over all $q,p\in K$ gives $\operatorname{diam}_{\mathbb H}(K)	\leq	\operatorname{diam}_{\mathbb C_I}(E).$ Thus, we have $	\operatorname{diam}_{\mathbb H}(K)	= \operatorname{diam}_{\mathbb C_I}(E).$ Hence,  $$	\gamma_{\mathbb H}(K)\leq
	\operatorname{diam}_{\mathbb H}(K)	\leq		4\gamma_{\mathbb H}(K).$$
\end{proof}

In the complex plane,  the  area-capacity estimate 
$\operatorname{Area}(E)\leq 4\pi\,\gamma_{\mathbb C}(E)^2,$  shows that analytic capacity controls the planar area of bounded measurable set $E\subset\mathbb C$ \cite{rational app}.

It is evident that since, $\mathbb{H} \simeq \mathbb{R}^4$ the appropriate geometric quantity associated with a measurable quaternionic set is its four dimensional Lebesgue volume. This naturally leads us to seek a corresponding geometric estimate in the quaternionic setting where planar area is replaced by four-dimensional volume (hyper volume) for a bounded measurable axially symmetric set. 

Let $S$ be an bounded measurable axially symmetric set. Fix $I\in\mathbb S$ and let $S_I=S\cap\mathbb C_I$ and $S^+_I=S_I\cap\{x+Iy:y>0\}.$

For an axially symmetric set, the whole set has been determined by the slice. That is, every point $x+Iy\in S^+_I,$ generates the two-sphere
$[x+Iy]=\{x+Jy:J\in\mathbb S\}.$ Thus, geometrically, $S$ is obtained by revolving $S^+_I$ around the real axis. This spherical completion is behind  the difference between the complex and quaternionic estimates. Since,  $[x+Iy]$ is a two-sphere of radius $y$, whose surface area is $4\pi y^2$.
Thus, since $\mathbb H\simeq\mathbb R^4$, the volume element is
$dV_4=4\pi y^2\,dx\,dy,$ and hence 
$$\operatorname{Vol}_4(S)=4\pi\iint_{E_+}y^2\,dx\,dy,$$ where $\operatorname{Vol}_4(S)$ denotes the 4-dimensional volume(hyper-volume).
Thus, unlike the planar area, the quaternionic volume is represented by a weighted area integral, with the additional weight $y^2$. A natural way to control this factor is through the diameter of $S$. Thus, by using area-analytic capacity estimate and Theorem \ref{complexslice quaternionic} we have this result.
\begin{thm}\label{volumeac}
	Let $S\subset\mathbb H$ be a bounded measurable axially symmetric set. Then $\operatorname{Vol}_4(S)	\leq 2\pi^2\operatorname{diam}(S)^2	\gamma_{\mathbb H}(S)^2.$
\end{thm}
\begin{proof}
	We have,$\operatorname{Vol}_4(S) =	4\pi\iint_{{S^+}_I} y^2\,dx\,dy.$ Since $S_I$ is symmetric with respect to the real axis, $\iint_{S_I} y^2\,dA	= 	2\iint_{S^+_I}y^2\,dx\,dy.$ Thus, 
	\begin{equation}\label{volume}
		\operatorname{Vol}_4(S)	= 2\pi\iint_{S_I} y^2\,dA.
	\end{equation}
	 Let $D=\operatorname{diam}_\mathbb H(S).$  For every $x+Iy\in E$, axial symmetry gives $x+Iy, x-Iy\in K.$ Hence, $2|y|	= 	\left|(x+Iy)-(x-Iy)\right|	\leq D.$ Thus, $y^2\leq\frac{D^2}{4}.$	Now by equation \ref{volume}, we have
	\[	\operatorname{Vol}_4(S)\leq2\pi\frac{D^2}{4}\iint_{S_I} dA=\frac{\pi D^2}{2}\operatorname{Area}(S_I).	\] Since, by the  area estimate for analytic capacity [Theorem 2.3\cite{rational app}], $$\operatorname{Area}(S_I)	\leq4\pi\,\gamma_{\mathbb C}(E)^2,$$ we have $\operatorname{Vol}_4(S)\leq2\pi^2D^2\gamma_{\mathbb C}(S_I)^2.$	
	Now, by Theorem \ref{complexslice quaternionic} we have
	\[\operatorname{Vol}_4(S)\leq	2\pi^2\operatorname{diam}_\mathbb H(S)^2
	\gamma_{\mathbb H}(S)^2.\]
	\end{proof}

The volume inequality in terms  of diameter is  a consequence of the four-dimensional geometry rather than  the choice of
complex slice. In fact, for an axially symmetric set, the diameter is independent of the chosen slice since,
$\operatorname{diam}_{\mathbb H}(S) = \operatorname{diam}_{\mathbb C_I}(S\cap\mathbb C_I),$ for $I\in\mathbb S.$
This is not the same case for a general non-axially symmetric
quaternionic set, for which the diameters of different slices need not coincide.  Here is an example illustrating this Theorem. 

\begin{exa}
	For, fix a complex slice $\mathbb C_I$ and consider $E_R = \overline{\Delta(iR,1)}\cup\overline{\Delta(-iR,1)},$ for $R>1$ where $\Delta(a,r)$ be a disc of center $a$ and radius $r$. Let $K_R$ be its axially symmetric completion, that is, $K_R = \cup_{q\in E_R}[q]$. Then $K_R$ is axially symmetric and $K_R \cap \mathbb C_I = E_R.$ Now, by the semi-additivity of analytic capacity and the fact that the analytic capacity of a  disc of radius $r$ is $r$, 
	\[
	\begin{aligned}
		\gamma_{\mathbb C}(E_R) &\leq	C(\gamma_{\mathbb C}(\overline{\Delta(iR,1)})
		+ \gamma_{\mathbb C}(\overline{\Delta(-iR,1)}))&=2C
	\end{aligned}
	\]
	for some arbitrary constant $C$. Hence by Theorem \ref{complexslice quaternionic} we have $\gamma_{\mathbb H}(K_R) = \gamma_{\mathbb C}(E_R) \leq2C,$ uniformly in $R$.
	Consider the upper disc  of ${\Delta(iR,1)}$. For each point $x+Iy$ its axial completion  generates a two-dimensional sphere of radius $y$. Since for every point in the disc, $R-1 <y<R+1$  and the disc has fixed positive area, the volume ${Vol}_4(K_R) =	4\pi\iint_{{_R}_+}y^2\,dx\,dy$ gives $$4\pi^2 (R-1)^2 \leq\operatorname{Vol}_4(K_R)\leq 4\pi^2 (R+1)^2 .
	$$ Thus, $\frac{\operatorname{Vol}_4(K_R)}{\gamma_{\mathbb H}(K_R)^4}
	\longrightarrow\infty$ as $R\longrightarrow\infty.$ 
\end{exa}
If $S$ is compact, axially symmetric, and slice connected,
then the previously established diameter-capacity estimate $$\operatorname{diam}_{\mathbb H}(K)\leq4\gamma_{\mathbb H}(S)$$
allows the diameter to be absorbed into the slice regular capacity.
\begin{cor}
	Let $K\subset\mathbb H$ be a compact axially symmetric
	slice-connected set. Then $\operatorname{Vol}(K) \leq 32\pi^2\gamma_{\mathbb H}(K)^4.$
\end{cor}
\begin{proof}
	By the theorem \ref{volumeac} we have $	\operatorname{Vol}(K)\leq
	2\pi^2\operatorname{diam}(K)^2 \gamma_{\mathbb H}(K)^2.$
	Since $K$ is compact, axially symmetric and slice-connected, by Theorem \ref{diameterac} we have $\operatorname{diam}(K)\leq4\gamma_{\mathbb H}(K).$
	Hence
	\[
	\begin{aligned}
		\operatorname{Vol}(K)
		&\leq
		2\pi^2
		\left(4\gamma_{\mathbb H}(K)\right)^2
		\gamma_{\mathbb H}(K)^2
		&=
		32\pi^2\gamma_{\mathbb H}(K)^4.
	\end{aligned}
	\]
	Hence $	\operatorname{Vol}(S)\leq	32\pi^2\gamma_{\mathbb H}(S)^4.$
\end{proof}
Analytic capacity provides estimates for admissible function  defined on its complement of  compact set $E\subset\mathbb C$. It shows that the value of an admissible function and its derivative are controlled by the analytic capacity and the distance from the point to the complement set.
\begin{thm}\label{esticomp}\cite{rational app}
	Let $f$ be admissible for $K$ and $z \in \Omega_{\mathbb C}(E)$, then $|f(z)| \leq \frac{\gamma_{\mathbb{C}}(E)}{d(z,E)},$ and
		$|f'(z)| \leq \frac{4\gamma_{\mathbb{C}}(E)}{d(z,E)^2},$
	where $d(z, E)$ is the distance from $z$ to $E$.
\end{thm}
 We now establish corresponding estimates in the quaternionic setting for compact axially symmetric sets. 
\begin{thm} \label{distancequatern}
	Let $K\subset\mathbb H$ be a compact axially symmetric set, and let 
	$f$ be admissible for $K$. Then, for every $q\in\Omega_{\mathbb{H}}(K)$,
	\[	|f(q)|	\leq	\frac{\gamma_{\mathbb H}(K)}	{d_{\mathbb H}([q],K)}\ \text{and}\]
	\[	|\partial_s f(q)|	\leq \frac{4\gamma_{\mathbb H}(K)}	{d_{\mathbb H}([q],K)^2},\]
	where $\partial_s f$ denotes the slice derivative of $f$ and $d_{\mathbb H}([q],K)= \operatorname{inf}_{p \in K} |q-p|$ .
\end{thm}

\begin{proof}
	Let $q_0\in\Omega_{\mathbb H}(K)$. Choose $I\in\mathbb S$ such that	$q_0\in\mathbb C_I,$ and set $E=K\cap\mathbb C_I.$ 	
	Next,  to prove  $d_{\mathbb H}(q_0,K) =	d_{\mathbb C_I}(q_0,E).$ Let, $q_0=x_0+y_0I,\ y_0\geq0$ and $	p=x+yJ\in K,\  y\geq0.$  Since $K$ is axial symmetric, $p_I=x+yI\in E.$ Now, 
	\[
	\begin{aligned}
		|q_0-p|^2 &=	|x_0-x+y_0I-yJ|^2	&=	(x_0-x)^2+y_0^2+y^2
		-2y_0y\,\langle I,J\rangle,
	\end{aligned}
	\]
	\[\text{and}\	|q_0-p_I|^2	= 	(x_0-x)^2+(y_0-y)^2.	\] where $\langle I,J\rangle$ denotes the Euclidean inner product in
	the imaginary subspace of $\mathbb H$.  Thus,
	\[	\begin{aligned}
		|q_0-p|^2-|q_0-p_I|^2	&=	2y_0y\bigl(1-\langle I,J\rangle\bigr) &\geq0,
	\end{aligned}\]
	since, $\langle I,J\rangle\leq1.$
	Hence $	|q_0-p_I|\leq|q_0-p|.$
	
	By taking the infimum over all $p\in K$, we have $	d_{\mathbb C_I}(q_0,E) \leq d_{\mathbb H}(q_0,K).$
	Conversely, since  $E\subset K,$	we immediately have
	$d_{\mathbb H}(q_0,K) \leq 	d_{\mathbb C_I}(q_0,E).$
	
	Since, $\Omega_{\mathbb{H}}(K)\cap\mathbb C_I	= \Omega_{\mathbb{C}_I}(E)$, the restriction $f_I=f|_{\Omega_{\mathbb C_I}(E)}$ is well defined. Since $f$ is slice regular, the  Lemma \ref{splitting lemma} gives an imaginary unit $J\in\mathbb S$, $J\perp I$, and holomorphic functions $F,G:\Omega_{\mathbb C_I}(E)\longrightarrow\mathbb C_I$ such that $$f_I(z)=F(z)+G(z)J.$$ Since $f$ is admissible for $K$, $$	\|f\|_{\infty,\Omega(K)}\leq1,\ f(\infty)=0.$$ Thus, $	\|f_I\|_{\infty,\Omega_{\mathbb C_I}(E)}\leq1,\	f_I(\infty)=0.$ We first prove the estimate for $f(q_0)$.
	
	It is sufficient to assume that  $f(q_0)\neq0.$ Let, $f(q_0)=a+bJ,\ a,b\in\mathbb C_I.$ Define $N_0:\mathbb H\longrightarrow\mathbb C_I$	by
	$$N_0(w_1+w_2J) 	= 	\frac{\overline a\,w_1+\overline b\,w_2}
	{|f(q_0)|},\ w_1,w_2\in\mathbb C_I$$ which is $\mathbb C_I$-linear functional. Now, by the Cauchy--Schwarz inequality,
	\[
	\begin{aligned}
		|N_0(w_1+w_2J)| &\leq	\frac{|a||w_1|+|b||w_2|}	{|f(q_0)|}\\
		&\leq\frac{\sqrt{|a|^2+|b|^2}	\sqrt{|w_1|^2+|w_2|^2}}
		{|f(q_0)|}	&=	|w_1+w_2J|
	\end{aligned}
	\]
	and thus, $\|N_0\|\leq1.$ Also, 
	\[
	\begin{aligned}
		N_0(f(q_0)) &=	\frac{\overline a\,a+\overline b\,b}
		{|f(q_0)|}\\
		&=	\frac{|a|^2+|b|^2}{|f(q_0)|} &=
		|f(q_0)|.
	\end{aligned}
	\]
	Define $h_0(z)=N_0(f_I(z)).$ 	Since $N_0$ is $\mathbb C_I$-linear and $f_I$ is holomorphic, $h_0$ is a $\mathbb C_I$-valued holomorphic function on	$\Omega_{\mathbb C_I}(E)$. Moreover, $|h_0(z)|\leq	|f_I(z)|\leq1$ 	and $h_0(\infty)=0.$ Thus $h_0$ is admissible for $E$ in the definition of analytic capacity. Also, since $	N_0(f(q_0)) = 	|f(q_0)|$, we have $|h_0(q_0)|=|f(q_0)|.$ Hence by Theorem \ref{esticomp}, $|h_0(q_0)|\leq	\frac{\gamma_{\mathbb C}(E_I)}	{d_{\mathbb C_I}(q_0,E_I)}.$ Thus, $|f(q_0)|\leq	\frac{\gamma_{\mathbb C}(E_I)}	{d_{\mathbb H}(q_0,K)}.$ Hence by Theorem \ref{complexslice quaternionic},
	\[	|f(q_0)|	\leq 	\frac{\gamma_{\mathbb H}(K)}
	{d_{\mathbb H}(q_0,K)}.	\]
	
	We now prove the estimate for the slice derivative.	It is enough to prove for the case that $\partial_s f(q_0)\neq0.$ Let, $\partial_s f(q_0)=a+bJ,\	a,b\in\mathbb C_I.$	Define the $\mathbb C_I$-linear functional $N_1:\mathbb H\longrightarrow\mathbb C_I$ by $N_1(w_1+w_2J) 	= 	\frac{\overline a\,w_1+\overline b\,w_2}{|\partial_s f(q_0)|}.$ 	Again, by the Cauchy--Schwarz inequality, $\|N_1\|\leq1$ and $$N_1(\partial_s f(q_0))=
	\frac{|a|^2+|b|^2}{|\partial_s f(q_0)|} = |\partial_s f(q_0)|.$$
	Define, $h_1(z):=N_1(f_I(z)).$	Then $h_1$ is holomorphic on $\Omega_{\mathbb{C}_I}(E)$, and $|h_1(z)| \leq1,\ h_1(\infty)=0.$
	Hence $h_1$ is admissible for $E$ in the definition of analytic capacity.
	Since $N_1$ is $\mathbb C_I$-linear, we have $h_1'(q_0) = 	N_1(f_I'(q_0)).$ Since the slice derivative is given by $$\partial_s f(q_0) 	= 	\frac{\partial f_I}{\partial z}(q_0)	= 	f_I'(q_0),$$ we have $h_1'(q_0)	= 	N_1(\partial_s f(q_0)).$  Hence, $|h_1'(q_0)|
	= |\partial_s f(q_0)|.$	Thus by the  derivative estimate for analytic capacity $|h_1'(q_0)|	\leq 	\frac{4\gamma_{\mathbb C}(E)}
	{d_{\mathbb C_I}(q_0,E)^2}.$ Hence, 
	\[
	\begin{aligned}
		|\partial_s f(q_0)|
		&=
		|h_1'(q_0)|	&\leq	\frac{4\gamma_{\mathbb C}(E)}
		{d_{\mathbb C_I}(q_0,E)^2} &=	\frac{4\gamma_{\mathbb H}(K)}
		{d_{\mathbb H}(q_0,K)^2}.
	\end{aligned}
	\]	
	Since $q_0\in\Omega_{\mathbb H}(K)$ was arbitrary, the two estimates hold for every $q\in\Omega_{\mathbb H}(K)$. 
\end{proof}
The point-wise estimate described above leads to quantitative bounds for the Laurent coefficients of admissible functions. In the  complex plane controls not only the values and derivatives of
admissible functions away from $E$, but also the coefficients
appearing in their expansions at infinity. Such coefficient estimates are particularly useful in  Vitushkin's approximation theory, where control of the tail of an expansion plays an important role in quantitative approximation arguments.

  We establish a corresponding estimate for its quaternionic Laurent coefficients using Laurent expansion at infinity.  The resulting estimate provides a quantitative link between the size of the coefficients and the quaternionic slice
regular capacity. Instead of  disc in complex setting we have the alternative  ball $\mathbb B_{\mathbb H}(0,R)$ centered at 0 and of radius $R$. 
\begin{thm}
	Let $K\subset B_{\mathbb H}(0,R)$ be a compact axially symmetric set.
	Let $f$ be admissible for $K$, and suppose that
	\[	f(q)=q^{-1}\partial_s f(\infty)+\sum_{n=2}^{\infty}q^{-n}a_n,
	\qquad z\in\mathbb H\setminus K.\]
	Then, for every $n\geq2$,
	\[	|a_n|\leq enR^{\,n-1}\gamma_{\mathbb H}(K).	\]
\end{thm}

\begin{proof}
	For $n\geq2$, choose $R_0>R$. Let $I\in\mathbb S$ and consider the
	circle $C_I(R_0)=\{z\in\mathbb C_I:|z|=R_0\}.$	Since $K\subset B_{\mathbb H}(0,R)$, 
	$$|z-q|	\geq ||z|-|q|| 	\geq R_0-R$$ for $z\in C_I(R_0)$ and $q\in K$. Hence, $d_{\mathbb H}(z,K)\geq R_0-R.$	Thus, by Theorem \ref{diameterac}, $$|f(z)| \leq	\frac{\gamma_{\mathbb H}(K)}{d_{\mathbb H}(z,K)}\leq
	\frac{\gamma_{\mathbb H}(K)}{R_0-R}.$$
	Now restrict $f$ to the slice $\mathbb C_I$. Since $f$ is slice
	regular and admissible, its restriction is holomorphic and has the
	Laurent expansion
	\[
	f(z)=\frac{\partial_s f(\infty)}{z} +\sum_{n=2}^{\infty}\frac{a_n}{z^n}.\]
	Now, by Cauchy's coefficient estimate,
	\[|a_n|	\leq 	R_0^{\,n-1}	\max_{|z|=R_0}|f(z)|.
	\]
	Hence, $|a_n| \leq 	\frac{R_0^{\,n-1}}{R_0-R}\,	\gamma_{\mathbb H}(K).$
	Choose, $R_0=\frac{nR}{n-1}.$ Then $R_0-R=\frac{R}{n-1},$ and hence
	$$|a_n| \leq 	nR^{\,n-1} 	\left(\frac{n}{n-1}\right)^{n-1}
	\gamma_{\mathbb H}(K).$$	Since, $\left(1+\frac1{n-1}\right)^{n-1}<e,$
	we have $|a_n| \leq enR^{\,n-1}\gamma_{\mathbb H}(K).$
\end{proof}
 In the complex setting, covering annuli by
discs of a prescribed radius is  important in the
localization and approximation arguments. The annuli $\{z\in\mathbb C:nr\leq |z-z_0|\leq(n+1)r\}$; $z_0\in \mathbb C$ is covered by $Cn$ discs of radius $r$.

In the quaternionic setting, we identify $\mathbb H$ with $\mathbb R^4$, and the corresponding annuli  covered by four-dimensional balls. Consider the annulus $A_n(q_0,r) = \{q\in\mathbb H:nr\leq |q-q_0|\leq(n+1)r\}$; $q_0\in \mathbb H$ and $r>0$. Then, its four-dimensional volume is  $\frac{\pi ^2}{2}r^4(4n^3+O(n^2)),$ whereas a ball of radius $r$ has volume  $\frac{\pi ^2}{2}r^4$. Thus, volume considerations naturally lead to a covering number of order $n^3$. The following Lemma makes this estimate precise and provides the geometric covering.
\begin{thm}\label{covering}
	There exists a universal constant $C_0>0$ such that for all
	$q_0\in\mathbb H$, $r>0$, and every integer $n\geq 1$, the annulus
	\[A_n(q_0,r) = 	\left\{	q\in\mathbb H:	nr\leq |q-q_0|\leq (n+1)r
	\right\}\] can be covered by at most $C_0n^3$ balls of radius $r$.
\end{thm}
\begin{proof}
	We identify $\mathbb H$ with $\mathbb R^4$ equipped with the Euclidean norm.	Fix $q_0\in\mathbb H$, $r>0$ and $n\geq1$. Since $A_n(q_0,r)$ is closed and bounded subset of $\mathbb R^4$, it is compact. 	
	Choose $q_1\in A_n(q_0,r)$. If there is another point in the annulus at distance at least $r$ from $q_1$, choose it as $q_2$ and continue this way. Since the annulus is compact, it will stop after some point. Let it be $q_N$. We have chosen a maximal $r$-separated subset $\{q_1,\ldots,q_N\}\subset A_n(q_0,r).$	We will show that
	$A_n(q_0,r)\subset \bigcup_{i=1}^{N}B(q_i,r).$
	
	Take an arbitrary point $q\in A_n(q_0,r)$. Suppose for contradiction that $q$ does not belong to any of the balls. Then $|q-q_i|\geq r, i=1,\ldots,N.$	Hence we can add $q$ to our family, contradicting the maximality of $\{q_1,\ldots,q_N\}.$ Hence there must exist some $i$ such that $q\in B(q_i,r).$ Thus, $A_n(q_0,r)\subset \bigcup_{i=1}^{N}B(q_i,r).$  
	
	Also since  $|q_i-q_j|\geq r, i\neq j,$ we get disjoint balls $	B\left(q_i,\frac r2\right), 	B\left(q_j,\frac r2\right).$ For, if some $x$ in both balls, then
	$|q_i-x|<\frac r2$ and $|x-q_j|<\frac r2$ and hence $q_i -q_j < r$ which is a	contradiction. Thus the balls $	B\left(q_i,\frac r2\right),\ i=1,\ldots,N$	are pairwise disjoint. Since $nr \leq |q_i-q_0| \leq (n+1)r$, every point $q\in B\left(q_i,\frac r2\right),$ satisfies 
	\[|q-q_0| \leq |q-q_i|+|q_i-q_0| <\frac r2+(n+1)r = 	\left(n+\frac32\right)r.\]
	Also, $	|q-q_0| \geq |q_i-q_0|-|q-q_i| \geq nr-\frac r2	= 	\left(n-\frac12\right)r.$ Thus, for all $i$,
	\[B\left(q_i,\frac r2\right)\subset\left\{q\in\mathbb H:	\left(n-\frac12\right)r	\leq |q-q_0|\leq\left(n+\frac32\right)r	\right\}\ \text{and thus}.\] 
	\[B\left(q_i,\frac r2\right)\subset B\left(q_0,(n+\frac 32)r\right)\]
	Thus, since the balls are pairwise disjoint,
	\[ \operatorname{Vol}\left(\bigcup_{i=1}^{N}B(q_i,\frac r2).\right) = \sum_{i=1}^{N}\left(B(q_i,\frac r2)\right) = N\operatorname{Vol}\left(B(q_i,\frac r2)\right).\]
	Hence, 	\[	N\operatorname{Vol}\left(B\left(0,\frac r2\right)
	\right)\leq \operatorname{Vol} \left(B\left(q_0,\left(n+\frac32\right)r\right)	\right).\]
	
	Since, the volume of a ball $B(0,r)$ in $\mathbb R^4$ is
	$\operatorname{Vol}(B(0,r))	= 	\frac{\pi^2}{2}r^4.$
	Hence, $N\frac{\pi^2}{2}\left(\frac r2\right)^4	\leq \frac{\pi^2}{2}
	\left(n+\frac32\right)^4r^4,$ and thus, $N\leq 16\left(n+\frac32\right)^4.$
	
	We know that the disjoint balls $B\left(q_i,\frac r2\right)$ are contained in the annulus \[\left\{q\in\mathbb H:\left(n-\frac12\right)r<
	|q-q_0|	< \left(n+\frac32\right)r\right\}.\]	Thus, 
	\[	N\frac{\pi^2}{2}\left(\frac r2\right)^4	\leq 	\frac{\pi^2}{2}
	\left[\left(n+\frac32\right)^4	- 	\left(n-\frac12\right)^4
	\right]r^4.\] Hence, 
	\[N\leq16\left(8n^3+12n^2+14n+5\right)\] and thus \[N\leq16n^3\left(8+\frac 12n+ \frac 14n^2+ \frac 5n^3\right) .\] 
	For $n\geq1$, $	N\leq C_0n^3$ with a universal constant $C_0$ after enlarging the constant to cover the small values of $n$. Thus,the annulus $A_n(q_0,r)$ can be covered by at most $C_0n^3$ balls of radius $r$. 
	Also, since $n\geq1$, we have $N\leq 16(32n^3) .$ Thus, $N\leq 624n^3$	and hence the annulus $A_n(q_0,r)$ can be covered by at most $624n^3$ balls of radius $r$. Thus we may take $C_0=624.$
\end{proof}
The number of balls required to cover the annulus in quaternions is of order $n^3$, which reflects the four-dimensional geometry of the annulus.  This estimate is particularly important
in Vitushkin's theory, where a set is decomposed into localized pieces and one have to control the total contribution of those pieces. 
\begin{thm}\cite{rational app}
	Let $p > 0$ be an integer, and let $\{E_j\}$ be a family of subsets of the bounded set $E \subset \mathbb{C}$, such that every disc $\Delta(z; \gamma(E))$ meets at most $p$
	of the sets $\{E_j\}$. Then $\sum \gamma(E_j) \leq 100 p \gamma(E).$ If $g_j$ is an analytic capacity admissible function for $E_j$, then $	\sum |g_j(z)| \leq 100 p.$
\end{thm}
We establish the corresponding statement for slice regular capacity in $\mathbb H$. The role of discs in the complex plane is now played by four-dimensional quaternionic balls. Combining the above covering estimate with Theorem \ref{estimate pointwise} allows  the members of the family to be summed over the annular decomposition.
\begin{thm}\label{coveringgamma}
	Let $K\subset\mathbb H$ be bounded and let 	$\{K_j\}$
	be a family of axially symmetric subsets of $K$.  Assume that every ball $\mathbb B(q,r)$ which meets atmost $p$ of the sets $K_j$ where $p>0$ is an integer. Then there exists a
	 constant $C_0>0$  such that \[\sum_{j\in J}\gamma_{\mathbb H}(K_j)	\leq 	C_0\,\gamma_{\mathbb H}(K).\]
	If $f_j$ is an admissible function for $K_j$, then
	\[\sum_{j\in J}|f_j(q)|\leq	C_1p\]	
	where  $C_1>0.$
\end{thm}
\begin{proof}
	Let $q_0\in\mathbb H$ and put $r=\gamma_{\mathbb H}(K).$
	If $r=0$, then $\gamma_{\mathbb H}(K_j) = 0$ for every $j$ by the monotonicity and the result holds. 
	
	Suppose that $\gamma>0$. 	Consider the ball  $B(q_0,r)$ and for $n\geq 1$ the annuli, 
	$$A_n(q_0,r)	= 	\{q\in\mathbb H:nr\leq|q-q_0|<(n+1)r\}.$$
	Thus, by the Theorem \ref{covering}, each annulus can be covered by at most $624n^3$	balls of radius $r$. Since, by hypothesis, every ball of radius $r$ meets at most $p$ of the sets
	$K_j$, it follows that at most $624pn^3$ of the sets $K_j$ meet the annulus $A_n(q_0,r)$. Denote these sets by $K_{n,1}, K_{n,2}, \cdots K_{n,m_n}$ where $m_n \leq 624pn^3$. For each set we have $d(q_0, K_{n,l})\geq nr$. 	Let $f_{n,l}$ be admissible for the corresponding $K_{n,l}$. Thus by theorem \ref{distancequatern} and since $	d(q_0,E_{n,k})\geq nr,$
	\[	|f_{n,l}(q_0)|	\leq 	\frac{\gamma_{\mathbb H}(K_{n,l})}	{d(q_0,K_{n,l})}	\leq 	\frac{\gamma_{\mathbb H}(K_{n,l})}{nr}.
	\]
	Choose $c_n=\sum_l\frac{\gamma_{\mathbb H}(K_{n,l})}{\gamma_{\mathbb H}(K)}$	with the summation over all sets meeting the $n$-th annulus. Since $\gamma_{\mathbb H}(K_{n,l})	\leq 	\gamma_{\mathbb H}(K),$
	we have $0\leq c_n\leq624pn^3.$
	Also, if $a = \sum_j\frac{\gamma_{\mathbb H}(K_j)}{\gamma_{\mathbb H} (K)}$ we have $\sum_{n=1}^{\infty} c_n \leq a$. Since $c_n\leq 624pn^3$, we have $\sum_{n=1}^{N} \frac{c_n}{n} \leq 624pN^3$ and $\sum_{n>N} \frac{c_n}{n} \leq \frac aN$. Thus for any integer $N\geq1$, 
	\[\sum_{n=1}^{\infty} \frac{c_n}{n} = \sum_{n=1}^{N} \frac{c_n}{n} +\sum_{n>N}\frac{c_n}{n} \leq 624pN^3 + \frac aN.\]
	If $a=0$, the desired estimate is immediate. Assume $a>0$ and let $x=\left(\frac{a}{624p}\right)^\frac 14.$ If $x\geq1$, choose $N=[x]$. Then, $x\leq N\leq 2x$  and thus, $624pN^3 \leq 8(624px^3)$ and $\frac aN \leq 624px^3$. Hence $\sum_{n=1}^{\infty} \frac{c_n}{n} \leq 9(624px^3) = (624p)^{\frac 14}a^{\frac 34} < 45p^{\frac 14}a^{\frac 34}$
	
	If $0<x<1$, we take $N=1$ then, $\sum_{n=1}^{\infty} \frac{c_n}{n} \leq a$ and $a\leq 624px^3 <45 p^{\frac 14}a^{\frac 34}.$ Thus in both cases, $\sum_{n=1}^{\infty} \frac{c_n}{n}<45 p^{\frac 14}a^{\frac 34}.$
	Since $|f_{n,l}(q_0)|\leq 	\frac{\gamma_{\mathbb H}(K_{n,l})}{nr}$, 
	\[\sum_{n,l} |f_{n,l}(q_0)|\leq  \sum_{n=1}^{\infty}\sum_l \frac{\gamma_{\mathbb H}(K_{n,l})}{nr} = \sum_{n=1}^{\infty}\frac{c_n}{n}\leq 45 p^{\frac 14}a^{\frac 34}.\]
	
	Since, $B(q_0,r)$  meets atmost $p$ of the sets $K_j$, the set $K_j$ which meet $B(q_0,r)$ contribute atmost $p$. For these sets we use $|f_{j}(q_0)| \leq 1.$ Hence, $$\sum_j|f_j(q_0)|\leq p+45p^{\frac 14}a^{\frac 34}.$$ For all $j$, choose  $h_j(\infty)=\frac{\gamma_{\mathbb H}(K_j)}{2}$, Let, $M = p+45p^{\frac 14}a^{\frac 34} $. Then $h(q) = \frac{\sum_{j}h_j(q)}{M}$ is admissible for $K$ and $|h(q)| \leq 1$. Also, $h(\infty) = \frac{1}{2M}\sum_j\gamma_{\mathbb H}(K_j) = a\gamma_\mathbb H(K) = \frac{ar}{2M}$. By the definition of  slice-regular capacity, $\frac{a}{2M}\leq\gamma_{\mathbb H}(E)=r$.
	Thus $a\leq2M = 2p + 90p^{\frac 14}a^{\frac 34}$. 	Hence $a \leq 91^{\frac 14}p.$ since, $a= \sum_j\frac{\gamma_{\mathbb H}(K_j)}{\gamma_{\mathbb H}(K)}$ we have \[\sum_{j\in J}\gamma_{\mathbb H}(K_j)	\leq 	C_0\,\gamma_{\mathbb H}(K).\]
	for a universal constant $C_0$. 
	
	Now, finally using $a<91^4p$, in $\sum_j|f_j(q_0)|\leq p+45p^{\frac 14}a^{\frac 34}$ we have $\sum_j|f_j(q_0)|\leq (1+45.91^3)p.$ So we write for some constant $C_1$ such that 
	\[\sum_{j\in J}|f_j(q)|\leq	C_1p.\]
\end{proof}
\section{Quaternionic analogue of $\gamma_+$-capacity and Cauchy capacity}
In this section we introduce a measure-theoretic counterparts of slice regular capacity, the positive quaternionic capacity ($\gamma^{+}_{\mathbb H}$) and quaternionic Cauchy capacity $(\gamma^{c}_{\mathbb H})$.

Analytic capacity plays a fundamental role in the theory of rational approximation and removable singularities. A related quantity is positive analytic capacity, denoted by
$\gamma_+(E)$. It is defined by restricting the extremal problem
to positive finite measures supported on $E$ whose Cauchy
transforms satisfy a uniform bound outside $E$. More precisely,
\[\gamma_+(E) = \sup\left\{
\mu(E):
\begin{array}{l}
	\mu\geq 0,\quad \operatorname{supp}\mu\subset E,\\[2pt]
	\|\mathcal C\mu\|_{L^\infty(\mathbb C\setminus E)}\leq 1
\end{array}
\right\},
\]
where
\[
\mathcal C\mu(z)
=
\int \frac{1}{z-w}\,d\mu(w)
\]
is the Cauchy transform, with the normalization fixed as above.

The importance of positive analytic capacity comes from the fact
that it permits analytic capacity to be studied through positive
measures and singular integral estimates. In particular, Tolsa
proved that analytic capacity and positive analytic capacity are
comparable \cite{tolsa}. Before Tolsa's work, the geometric characterization of sets of positive analytic capacity was a major open problem.  
\begin{thm}\label{tolsa1.1}\cite{tolsa}
	There exists an absolute constant $A$ such that $\gamma(E) \leq A\gamma_{+}(E)$ for any compact set $E$.
\end{thm}
This theorem is a fundamental result in the modern theory of
analytic capacity. It allows estimates for positive measures and
their Cauchy transforms to be used in the study of analytic
capacity itself.
Tolsa  proved the semi-additivity of analytic capacity, resolving a longstanding problem posed by Vitushkin.

Also, Tolsa characterized continuous analytic
capacity using curvature estimates for measures with zero linear
density and proved its countable semi-additivity. This result also
established the inner boundary conjecture.
Tolsa's results on analytic capacity and continuous analytic
capacity resolved two important open problems in complex
analysis. The first concerned semi-additivity of analytic capacity
and the characterization of removable singularities. The second
concerned semi-additivity of continuous analytic capacity and
the inner boundary conjecture, with consequences for uniform
rational approximation.

$\gamma^{+}_{\mathbb H}$ is defined through positive finite Radon measures supported on a compact set whose quaternionic Cauchy transforms satisfy a suitable uniform boundedness condition. In the  complex theory, positive measures helped to study analytic capacity more geometrically. In particular, Theorem \ref{tolsa1.1} by Tolsa \cite{tolsa}, establishes the comparability between analytic capacity and its positive-measure counterpart. In this section we establish the corresponding result in the quaternionic setting. 

Let $K\subset\widehat{\mathbb H}$ be a compact axially symmetric set. For a positive finite Radon measure $\mu$ supported on $K$, define its quaternionic Cauchy transform by
\begin{equation*}
	\mathcal C_{\mathbb H} ^\mu(q)=\int_K S_L^{-1}(s,q)\,d\mu(s),\qquad q\in\mathbb H\setminus K,
\end{equation*}
where $S_L^{-1}(s,q)$ is the left slice Cauchy kernel.
Now we can define Quaternionic $\gamma_+$-capacity, an analogue of $\gamma_{+}$ capacity.
\begin{defn}
	Let $K\subset\widehat{\mathbb H}$ be a compact axially symmetric set. The
	quaternionic $\gamma_+$-capacity of $K$ denoted by $\gamma^{+}_{\mathbb H}(K)$ is defined by
	\begin{equation*}
		\gamma^{+}_{\mathbb H}(K):=\sup_{\mu}\mu(K),
	\end{equation*}
	where the supremum is taken over all positive finite Radon measures
	$\mu$ satisfying $\operatorname{supp}\mu\subset K$
	and
	$\left\|\mathcal C_{\mathbb H}^\mu\right\|_{L^\infty(\widehat{\mathbb H}\setminus K)}\leq 1.$
\end{defn}
\begin{lem}
	\label{cauchy-kernel-infinity}
	Let $K$ be a compact set. Then, uniformly for $s\in K$,
	\begin{equation*}
		S_L^{-1}(s,q) = -q^{-1}+O(|q|^{-2}),\ |q|\longrightarrow\infty.
	\end{equation*}
\end{lem}
\begin{proof}
	Let $s_0 =\operatorname{Re}(s).$ Since $s_0,|s|^2\in\mathbb R$, they commute with $q$. Hence,
	\begin{align*}
		q^2-2\operatorname{Re}(s)q+|s|^2 &=
		q^2-2s_0q+|s|^2 &= q^2\left(1-2s_0q^{-1}+|s|^2q^{-2}\right).
	\end{align*}
	Thus, $$\left(q^2-2s_0q+|s|^2\right)^{-1} = \left(1-2s_0q^{-1}+|s|^2q^{-2}\right)^{-1}q^{-2}.$$
	Hence,
	\begin{align*}
		S_L^{-1}(s,q)
		&=
		-\left(
		1-2s_0q^{-1}+|s|^2q^{-2}
		\right)^{-1}
		q^{-2}(q-\overline{s})\\ 
		&=
		-\left(
		1-2s_0q^{-1}+|s|^2q^{-2}
		\right)^{-1}
		\left(q^{-1}-q^{-2}\overline{s}\right).
	\end{align*}
	Since $K$ is compact, there exists $M>0$ such that $	|s|\leq M,$ for $s\in K.$
	Thus, for sufficiently large $|q|$,
	\begin{align*}
		\left|2s_0q^{-1}-|s|^2q^{-2}\right|&\leq		2M|q|^{-1}+M^2|q|^{-2}<\frac12
	\end{align*}
	uniformly for $s\in K$.
	Thus, the quaternionic Neumann series gives
	\begin{align*}
		\left(
		1-2s_0q^{-1}+|s|^2q^{-2}
		\right)^{-1}
		&=
		1+O(|q|^{-1}),
	\end{align*}
	uniformly for $s\in K$. Hence,
	\begin{align*}
		S_L^{-1}(s,q)
		&=
		-\left(1+O(|q|^{-1})\right)
		\left(q^{-1}-q^{-2}\overline{s}\right) \notag\\
		&=
		-q^{-1}+O(|q|^{-2}),
	\end{align*}
	uniformly for $s\in K$. 
\end{proof}
We next obtain the corresponding expansion for the quaternionic
Cauchy transform of a finite measure.
\begin{lem}\label{cauchy transform infinity}
	Let $K$ be compact and let $\mu$ be a finite positive
	Radon measure supported on $K$. Then
	\[\mathcal C_{\mathbb H} ^\mu(q)
	=-\mu(K)q^{-1} + O(|q|^{-2}), \qquad |q|\longrightarrow\infty.\]
	In particular, the coefficient of $q^{-1}$ in the expansion of 	$\mathcal C_{\mathbb H}^\mu$ at infinity is $-\mu(K)$.
\end{lem}
\begin{proof}
	By Lemma~\ref{cauchy-kernel-infinity}, uniformly for $s\in K$,
	\[S_L^{-1}(s,q)= -q^{-1}+O(|q|^{-2}).\]
	Hence,
	\begin{align}
		\mathcal C_{\mathbb H} ^\mu(q)
		&=\int_K S_L^{-1}(s,q)\,d\mu(s) \notag\\
		&= 	\int_K\left(-q^{-1}+O(|q|^{-2})\right)d\mu(s)\\
		&=-q^{-1}\mu(K)+O(|q|^{-2})\mu(K).
	\end{align}
	Since $\mu(K)<\infty$, the last term is	$O(|q|^{-2})$. Thus,
	\[\mathcal C_{\mathbb H} ^\mu(q)=-\mu(K)q^{-1}+O(|q|^{-2}).\]
\end{proof}
From this lemma, it is clear that $\mathcal C_{\mathbb H}^\mu (\infty) = 0$ and $|\partial_s(\mathcal C_{\mathbb H}^\mu)(\infty)| = \mu(K).$ Next lemma gives the  comparison between the quaternionic $\gamma^{+}_{\mathbb H}$-capacity and slice regular capacity.
\begin{lem}	\label{gamma plus less gamma}
	Let $K$ be a compact axially symmetric set. Then,
	\[\gamma^{+}_{\mathbb H}(K)\leq\gamma_{\mathbb H}(K).\]
\end{lem}

\begin{proof}
	Let $\mu$ be any measure admissible in the definition of
	$\gamma^{+}_{\mathbb H}(K)$. Thus, $\mu\geq0,\ 	\operatorname{supp}\mu\subset K,$ and $	\left\|
	\mathcal C_{\mathbb H}^\mu\right\|_{L^\infty(\mathbb H\setminus K)}\leq1.$  Also, $\mathcal C_{\mathbb H}^\mu$ is slice regular on $\widehat{\mathbb H}\setminus K$, bounded  by $1$, and vanishes at infinity. Thus, it is an admissible function in the definition of the slice regular capacity $\gamma_{\mathbb H}(K)$. As, $|\partial_s(\mathcal C_{\mathbb H}^\mu)(\infty)| = \mu(K)$ we have, 	
	\[\mu(K)\leq\gamma_{\mathbb H}(K).\]
	Since the above inequality holds for every measure $\mu$ admissible 	for $\gamma^{+}_{\mathbb H}(K)$, taking the supremum over all such 	measures gives
	\[\gamma^{+}_{\mathbb H}(K)\leq\gamma_{\mathbb H}(K).\]
\end{proof}
We next establish the relation between the complex positive analytic capacity of the planar compact set $E$ and the quaternionic positive capacity of compact axially symmetric set $K$. 	
\begin{lem}\label{slice gammaplus}
	Let $K$ be a compact axially symmetric set and let $I\in\mathbb S$. Let $E=K\cap\mathbb C_I.$ Then
	\begin{equation*}
		\gamma_+(E)\leq	2\gamma^{+}_{\mathbb H}(K).
	\end{equation*}
\end{lem}
\begin{proof}
	Let $\varepsilon>0$ and $\nu$ be an admissible  measure in the definition of $C\nu$, the planar Cauchy transform. Then by the definition of the  positive analytic capacity, there exists a positive finite Radon measure $\nu$ supported on $E$ such that
	\[
	\|C\nu\|_{L^\infty(\widehat{\mathbb C}_I\setminus E)}
	\leq 1\ \text{and}\ \nu(E)>\gamma_+(E)-\varepsilon.
	\]
	We consider $\nu$ as a positive finite Radon measure on $\mathbb H$ supported on $K$ and consider its quaternionic Cauchy transform
	\[F(q) =\mathcal C_{\mathbb H}^\nu(q)
	=
	\int_E S_L^{-1}(s,q)\,d\nu(s),
	\qquad q\in\mathbb H\setminus K.
	\]
	Then for, $s,q\in\mathbb C_I$, the quaternionic Cauchy kernel reduces to 	the  complex Cauchy kernel. That is,
	\[S_L^{-1}(s,q)=\frac{1}{s-q},\]
	and hence 
	\[F(q)	= 	\int_E\frac{1}{s-q}\,d\nu(s) =	C\nu(q),
	\qquad q\in \widehat{\mathbb C}_I\setminus E.\]
	Thus, $\sup_{q\in \widehat{\mathbb C}_I\setminus E}|F(q)|
	\leq1.$
	
	Since $K$ is axially symmetric, if $q=x+Jy\in\mathbb H\setminus K,$ for $J\in\mathbb S,$	then $x+Iy\notin K
	\ \text{and}\ x-Iy\notin K.$ Thus by  the slice representation formula 
	\[F(x+Jy) = \frac12 \left[(1-JI)F(x+Iy)
	+(1+JI)F(x-Iy)\right].\]
	Hence,
	\begin{align*}
		|F(x+Jy)|
		&\leq\frac12 \left(|1-JI|\,|F(x+Iy)|+|1+JI|\,|F(x-Iy)|\right).\\
		&\leq 2\sup_{\zeta\in\widehat{\mathbb C}_I\setminus E}|F(\zeta)|.
	\end{align*}
	Thus $\|\mathcal C_{\mathbb H}^\nu\|_{L^\infty(\widehat{\mathbb H}\setminus K)}\leq2.$		
	Hence the measure $\mu=\frac12\nu$
	is admissible in the definition of 	$\gamma^{+}_{\mathbb H}(K)$, since, $\left\|\mathcal C_{\mathbb H}^{\mu}
	\right\| _{L^\infty(\mathbb H\setminus K)}\leq1.$
	Also, $\mu(K)	= 	\frac12\nu(E).$
	Hence, by the definition of $\gamma_{+} ^\mathbb H$,
	\[\gamma^{+}_{\mathbb H}(K)\geq\frac12\nu(E).\]
	Since $\nu(E)>\gamma_+(E)-\varepsilon,$ we have
	$\gamma^{+}_{\mathbb H}(K)\geq\frac12		\bigl(\gamma_+(E)-\varepsilon\bigr)$ and for $\varepsilon\to0$ gives 		\[
	\gamma^{+}_{\mathbb H}(K)
	\geq
	\frac12\gamma_+(E).
	\]
\end{proof}
From Lemma \ref{gamma plus less gamma} we have
\begin{equation}\label{gammaplus}
	\gamma^{+}_{\mathbb H}(K)\leq\gamma_{\mathbb H}(K).
\end{equation}
Next we prove  Quaternionic analogue of Tolsa's Theorem \ref{tolsa1.1} which is the converse part of \eqref{gammaplus}.
\begin{thm}
	\label{quaternionic tolsa}
	Let $K\subset\mathbb H$ be a compact axially symmetric set. Then there exists an absolute constant $A_{\mathbb H}>0$ such that
	\begin{equation*}
		\gamma_{\mathbb H}(K)
		\leq
		A_{\mathbb H}\gamma^{+}_{\mathbb H}(K).
	\end{equation*}
\end{thm}
\begin{proof}
	Fix $I\in\mathbb S$ and set $E=K\cap\mathbb C_I.$
	Since $K$ is axially symmetric, then by Theorem \ref{complexslice quaternionic}
	\[\gamma_{\mathbb H}(K)=\gamma_{\mathbb C}(E).\]
	Now by Theorem \ref{tolsa1.1}, there exists an absolute constant $A>0$ such that 
	\[\gamma_{\mathbb{C}}(E)\leq A\gamma_+(E).\]
	On the other hand, by Lemma \ref{slice gammaplus},
	\[\gamma_+(E)\leq2\gamma^{+}_{\mathbb H}(K).\]
	Thus,
	\begin{align*}
		\gamma_{\mathbb H}(K)
		&=\gamma_{\mathbb C}(E)\\
		&\leq A\gamma_+(E)\\
		&\leq 2A\gamma^{+}_{\mathbb H}(K).
	\end{align*}
	Thus, by taking $A_{\mathbb H} = 2A,$ we have the result.
\end{proof}
From these all results we have $$\gamma^{+}_{\mathbb H}(K)\leq\gamma_{\mathbb H}(K)\leq A\gamma^{+}_{\mathbb H}(K)$$
and $$\gamma_+(E)\leq\gamma^{+}_{\mathbb H}(K)\leq\gamma_{\mathbb H}(K)=\gamma_{\mathbb{C}}(E) \leq B\gamma_+(E)$$
for some constants $A$ and $B$.

Now we  introduce quaternionic analogue of  Cauchy capacity  and discuss the Painlev\'e problem \cite{murai12, younis} regarding the equality of Cauchy capacity and analytic capacity and discuss the  quaternionic version of  Painlev\'e problem.

Following Younsi \cite{younis}, the Cauchy capacity of $E$ is defined by
\[
\gamma_c(E)
=
\sup\left\{
|\mu(E)|:
\operatorname{supp}\mu\subset E,\
|\mathcal C\mu(z)|\leq 1
\text{ for all }z\in\mathbb C\setminus E
\right\},
\]
where the supremum is taken over all finite complex Borel measures
supported on $E$ \cite{younis}.

The Cauchy capacity is related to analytic capacity and positive
analytic capacity $\gamma_+$ by
\[
\gamma_+(E)\leq\gamma_c(E)\leq\gamma(E)
\leq C\gamma_+(E),
\]
where $C$ is an absolute constant. The comparison between
$\gamma$ and $\gamma_+$ follows from Tolsa's theorem \cite{tolsa}.

\begin{thm}\cite{havinsonmain}
	\label{havinson}
	Let $E\subset\mathbb C$ be compact and have finite Painlev\'e
	length $k$. If $f$ is bounded and holomorphic on
	$\widehat{\mathbb C}\setminus E$, with $f(\infty)=0$, then
	there exists a finite complex Borel measure $\mu$ supported on 	$E$ such that $	f(z)=\mathcal C\mu(z),\  z\in\widehat{\mathbb C}\setminus E,$
	and $\|\mu\|\leq \frac{k}{2\pi}\|f\|_\infty.$
\end{thm}

Applying this representation theorem to the Ahlfors function
of $E$ yields the following consequence.
\begin{thm}
	\label{gamma-gammac-finite-length}
	If $E\subset\mathbb C$ is compact and has finite Painlev\'e
	length, then $\gamma(E)=\gamma_c(E).$
\end{thm}

The finite Painlev\'e-length hypothesis is important. In general,
analytic capacity and Cauchy capacity need not be equal. Younsi
studied sufficient conditions for their equality and showed that
equality of the two capacities does not necessarily imply that
the Ahlfors function is the Cauchy transform of a measure
supported on the compact set.

In quaternionic we define the quaternionic analogue of Cauchy capacity is defined as follows:	
\begin{defn}
	Let $K\subset\mathbb H$ be a compact axially symmetric set. The quaternionic Cauchy capacity of $K$ is defined by
	\[\gamma_{c}^ \mathbb H(K) = \sup_{\mu} |\mu(K)|\] where $\mu$ is a finite quaternion-valued Borel measure, $\operatorname{supp}\mu\subset K,$ $\mathcal C_\mu^{\mathbb H}(q)|\leq1 $ for all $q\in\widehat{\mathbb H}\setminus K$.	
\end{defn} 
In the following result we establish monotonicity of the quaternionic Cauchy capacity.
\begin{lem}
	\label{monotonicitycauchy}
	Let $K_1,K_2\subset\mathbb H$ be compact axially symmetric sets with $K_1\subset K_2.$	Then
	\[\gamma^{c}_{\mathbb H}(K_1)\leq \gamma^{c}_{\mathbb H}(K_2).\]
\end{lem}

\begin{proof}
	Let $\mu$ be any measure admissible for $\gamma^{c }_{\mathbb H}(K_1)$. Then 	\[\operatorname{supp}\mu\subset K_1\ \text{and}\ |\mathcal C_\mu^{\mathbb H}(q)|\leq1, \	q\in \widehat{\mathbb H}\setminus K_1.\]
	Since $K_1\subset K_2$, we have $\operatorname{supp}\mu\subset K_2$ and  $\widehat{\mathbb H}\setminus K_2 \subset\widehat{\mathbb H}\setminus K_1.$
	Thus, $|\mathcal C_\mu^{\mathbb H}(q)|\leq1,$ for $
	q\in \widehat{\mathbb H}\setminus K_2.$
	Thus $\mu$ is also admissible for $\gamma^{\mathbb H}_{c}(K_2)$.
	
	Also, since $\mu$ is supported on $K_1$, $
	\mu(K_2)=\mu(K_1).$ Hence
	\[|\mu(K_1)| = 	|\mu(K_2)| 	\leq 	\gamma^{c}_{\mathbb H}(K_2).
	\]
	Taking the supremum over all measures $\mu$ admissible for
	$\gamma^{c}_{\mathbb H}(K_1)$, we have 
	\[\gamma^{c}_{\mathbb H}(K_1)\leq \gamma^{c}_{\mathbb H}(K_2).\]
\end{proof}
We first establish the elementary inequality between the
quaternionic Cauchy capacity and the slice regular capacity.
\begin{lem}\label{cauchy less analytic}
	For every compact axially symmetric set $K\subset\mathbb H$,
	\[\gamma^{c }_{\mathbb H}(K)\leq\gamma_{\mathbb H}(K).\]
\end{lem}
\begin{proof}
	Let $\mu$ be admissible for $\gamma^{c }_{\mathbb H}(K)$. Thus
	$\operatorname{supp}\mu\subset K$ and $|\mathcal C_\mu^{\mathbb H}(q)|\leq1,$ for $q\in \widehat{\mathbb H}\setminus K.$	 Then, $\mathcal C_\mu^{\mathbb H}(q)$ 
	is slice regular on $\widehat{\mathbb H}\setminus K$ and satisfies
	\[\|\mathcal C_\mu^{\mathbb H}\|_{L^\infty(\widehat{\mathbb H}\setminus K)}\leq1.\]
	Also, $\mathcal C_\mu^{\mathbb H}(\infty)=0$ and $|\partial_s(\mathcal C_\mu^{\mathbb H}(\infty))|=|\mu(K)|.$	Hence  $\mathcal C_\mu^{\mathbb H}$ is admissible for the slice regular capacity, and thus
	\[	|\mu(K)|\leq\gamma_{\mathbb H}(K).\]
	Taking the supremum over all admissible measures $\mu$ gives
	\[	\gamma^{c }_{\mathbb H}(K)\leq\gamma_{\mathbb H}(K).\]
\end{proof}
Let $f_K$ denote the function of $K$ which is admissible in the definition of $\gamma_{\mathbb H}(K)$ and extremal function which attains the supremum. That is, $f_K\in \mathcal C_L(\mathbb{\widehat{H}}), f_K(\infty)=0, \|f_K\|_{\infty,\Omega(K)}\leq 1$, $f_K\in \mathcal{SR}_L(\mathbb H\setminus K) $ and $\|f_K\|_\infty = \gamma_{\mathbb H}(K).$ We say it as Quaternionic Ahlfors function of $K$. 

A compact set $E\subset\mathbb C$ has finite Painlevé length if there is a number $k$ such that every open neighborhood $U$ of $E$ contains a cycle $\Lambda$ surrounding $E$ that consists of finitely many disjoint analytic Jordan curves and has length less than $l$. The infimum of such numbers $l$ is called the Painlevé length of $E$ \cite{younis}.

We are defining Quaternionic Painlev\'e length for proving quaternionic Painlev\'e problem [Corollary 1.3, Theorem 1.4
\cite{younis}] related to Cauchy capacity.

Let $K\subset\mathbb H$ be a compact axially symmetric set. For $I\in\mathbb S$, put $E=K\cap\mathbb C_I.$ The quaternionic Painlev\'e length of $K$ is defined by
\[\mathcal L_{\mathrm P,\mathbb H}(K) =\frac{1}{4\pi}
\int_{\mathbb S} \ell_{\mathrm P}(E)\,dS(I),\]
where $\ell_{\mathrm P}(E)$ denotes the complex Painlev\'e length of the compact set $E\subset\mathbb C_I$ and $dS(I)$ denote the surface area measure on $\mathbb S$.
We say that $K$ has finite quaternionic Painlev\'e length if
$\mathcal L_{\mathrm P,\mathbb H}(K)<\infty.$

It is well defined since it is independent of choice $I\in \mathbb{S}$. For, let $I,J\in\mathbb S$. Since $K$ is axially symmetric, the map 
\[\Phi_{I,J}:\mathbb C_I\longrightarrow\mathbb C_J,\qquad
\Phi_{I,J}(x+Iy)=x+Jy,\]
is an isometry satisfying $\Phi_{I,J}(K\cap\mathbb C_I) 	= 	K\cap\mathbb C_J.$
Hence,
\[\ell_{\mathrm P}(K\cap\mathbb C_I) =
\ell_{\mathrm P}(K\cap\mathbb C_J).	\]
Since $\frac{1}{4\pi}dS(I) =1$, for $I\in \mathbb{S}$, it follows that
\[\mathcal L_{\mathrm P,\mathbb H}(K) =
\ell_{\mathrm P}(K\cap\mathbb C_I),
\qquad I\in\mathbb S.
\]
Thus the definition is independent of the choice of the
distinguished slice.
Hence if $K$ has finite quaternionic Painlev\'e length, 
\[\ell_{\mathrm P}(K\cap  C_I) 	= 	\mathcal L_{\mathrm P,\mathbb H}(K)<\infty.\]
The following lemma establishes the quaternionic analogue of the Havinson representation theorem.

It provides an estimate for the total variation of the representing measure. This estimate will be used subsequently to relate the quaternionic Painlev\'e length to the quaternionic Cauchy capacity.
\begin{lem}
	\label{quaternionic havinson}
	Let $K\subset\mathbb H$ be a compact axially symmetric set with finite quaternionic Painlev\'e length.
	Suppose that the quaternionic Ahlfors function $f_K$ is unique. Then there 	exists a finite quaternion valued  Borel measure $\mu$ supported on $K$ such that
	\[f_K(q) = 	\int_E S_L^{-1}(s,q)\,d\mu(s),\qquad
	q\in\widehat{\mathbb H}\setminus K.\]
	Moreover, 
	\[\|\mu\|\leq\frac{\mathcal L_{\mathrm P, \mathbb H}(K)}{2\pi}.\]
\end{lem}

\begin{proof}
	Fix $I\in\mathbb S$ and put $E=K\cap\mathbb C_I.$
	Since $K$ has finite quaternionic Painlev\'e length, $E$ has finite Painlev\'e length.
	Let $f_E$ be the planar Ahlfors function associated with $E$.  Then,
	\[\|f_E\|_\infty\leq1, 	\qquad 	f_E(\infty)=0, \text{and}\ f'_E(\infty)= \gamma_{\mathbb C}(E)>0. \]
	Since $E$ is axially symmetric, the function $g(z)=\overline{h_E(z)}$ which is also a planar Ahlfors function. Since, planar Ahlfors function is unique, $h_E(\bar{z}) = \overline{h_E(z)}$. Hence, $h_E$ is intrinsic. Thus, $h_E(x+Iy)=u(x,y)+Iv(x,y),$ and define $H$ by
	\[	H(x+Jy)=u(x,y)+Jv(x,y),\qquad J\in\mathbb S.\]
	Thus, $H$ is slice regular on $\Omega_{\mathbb H}(K)$. Also,
	\[
	\begin{aligned}
		|H(x+Jy)|^2
		&=u(x,y)^2+v(x,y)^2\\
		&=|h_E(x+Iy)|^2.
	\end{aligned}
	\]
	Thus, $H(\infty)=0,\ \text{and }\|H\|_{\infty,\Omega_{\mathbb H}(K)}\leq 1.$
	Then, by Theorem \ref{complexslice quaternionic}
	\[
	\begin{aligned}
		\partial_s H(\infty)=h_E'(\infty)=\gamma_{\mathbb C}(E)=\gamma_{\mathbb H}(K).
	\end{aligned}
	\]
	Hence $H$ is a normalized quaternionic Ahlfors function.
	Then by assumption of  uniqueness of $f_K$, we have
	$H=f_K.$
	In particular, $f_K(z)=h_E(z),$ for $ z\in\Omega_{\mathbb H}(K)\cap\mathbb C_I.$
	
	Since $\ell_p(E)<\infty$, then by  Theorem \ref{havinson}  yields a finite complex Borel 	measure $\nu$, supported on $E$, such that
	\[h_E(z)=\int_E\frac{1}{s-z}\,d\nu(s)\]
	and
	\begin{equation*}\label{eq:havinson-estimate}
		\|\nu\|\leq\frac{\ell_p(E)}{2\pi}.
	\end{equation*}
	Identifying $\mathbb C_I$ with a subfield of $\mathbb H$,
	we regard $\nu$ as a quaternion-valued Borel measure
	$\mu$ supported on $E$, by defining $\mu(A)=\nu(A)$
	for every Borel subset $A$ of $E$. In particular,
	$|\mu\|=\|\nu\|.$
	
	Define
	\[	F(q)=\int_E S_L^{-1}(s,q)\,d\mu(s),\qquad q\in\Omega_{\mathbb H}(K).\]
	Then,	$F$ is left slice regular on $\Omega_{\mathbb H}(K)$.
	For $z\in\Omega_{\mathbb H}(K)\cap\mathbb C_I$, we have
	\[	S_L^{-1}(s,z)=\frac{1}{s-z}.	\]
	Thus,
	\[
	\begin{aligned}
		F(z)&=\int_E S_L^{-1}(s,z)\,d\mu(s)\\
		=\int_E\frac{1}{s-z}\,d\nu(s)\\
		&=h_E(z) =f_K(z).
	\end{aligned}
	\]
	Since both $F$ and $f_K$ are left slice regular on the
	axially symmetric slice domain $\Omega_{\mathbb H}(K)$, by Theorem \ref{identity principle}, 
	$F(q)=f_K(q),$ for $q\in\Omega_{\mathbb H}(K).$	Hence,
	\[f_K(q)=\int_E S_L^{-1}(s,q)\,d\mu(s),\qquad q\in\Omega_{\mathbb H}(K).\]
	Also,
	\[\begin{aligned}
		\|\mu\|
		&=\|\nu\|\leq\frac{\ell_p(E)}{2\pi}=\frac{\mathcal L_{P,\mathbb H}(K)}{2\pi}.
	\end{aligned}
	\]
	\end{proof}
The question regarding next theorem in classical setting was raised  first time by Murai \cite{murai12} was later studied by Havinson, Younsi \cite{havinson8,havinson9, younis} etc. M Younsi proved that Cauchy capacity and slice regular capacity of $E\subset \mathbb C$ are equal if $E$ has finite Painlev\'e length. We  hereby state its quaternionic analogue regarding the equality of quaternionic Cauchy capacity and slice regular capacity.
\begin{cor}
	\label{quaternionic corollary}
	Let $K\subset\mathbb H$ be a compact axially symmetric set with finite quaternionic Painlev\'e length. Then
	\[\gamma^{c }_{\mathbb H}(K)=\gamma_{\mathbb H}(K).\]
\end{cor}

\begin{proof}
	Fix $I\in\mathbb S$ and put $E=K\cap\mathbb C_I.$
	Let $K\subset\mathbb H$ be a compact axially symmetric set with finite quaternionic Painlev\'e length.Then by
	Lemma \ref{quaternionic havinson}, there exists a finite quaternion	Borel measure $\mu$ supported on $K$ such that
	\[f_K(q) = 	\int_E S_L^{-1}(s,q)\,d\mu(s),
	\qquad q\in\widehat{\mathbb H}\setminus K.
	\]
	Since $f_K$ is the quaternionic Ahlfors function,
	\[|f_K(q)|\leq1, \qquad q\in\widehat{\mathbb H}\setminus K.
	\]
	Hence $\mu$ is admissible in the definition of
	$\gamma^{c }_{\mathbb H}(K)$. Thus,
	\[\gamma^{c }_{\mathbb H}(K) \geq |\mu (K)|.\]
	On the other hand by Lemma \ref{cauchy transform infinity},
	\[f_K(q) = 	-q^{-1}\mu(K)+O(|q|^{-2}),\qquad |q|\longrightarrow\infty.\]
	Thus, the coefficient at infinity 
	\[|\partial_s(f_K)(\infty)| = 	|\mu(K)|.\]
	Since $f_K$ is the quaternionic Ahlfors function,
	\[|\partial_s(f_K)(\infty)| = \gamma_{\mathbb H}(K).\]
	Thus
	\[|\mu(K)| = \gamma_{\mathbb H}(K)\] and hence 
	\[\gamma^{c }_{\mathbb H}(K) \geq 	|\mu(K)| = 		\gamma_{\mathbb H}(K).\]
	On the other hand, by Lemma \ref{cauchy less analytic},
	\[\gamma^{c }_{\mathbb H}(K)\leq\gamma_{\mathbb H}(K).\]
	Hence $\gamma^{c }_{\mathbb H}(K)=\gamma_{\mathbb H}(K).$
\end{proof}
\begin{cor}
	Let $K\subset\mathbb H$ be compact and axially symmetric and let $E=K\cap\mathbb C_I.$	If $K$ has finite Painlev\'e length, then
	\[\gamma^{c }_{\mathbb H}(K) = \gamma_{\mathbb H}(K) =
	\gamma_{\mathbb{C}_I}(E).\]
\end{cor}

\begin{proof}
	From Theorem \ref{quaternionic corollary}, we have
	\[\gamma^{c }_{\mathbb H}(K) = \gamma_{\mathbb H}(K).\]
	Also, by Theorem \ref{complexslice quaternionic},		
	\[	\gamma_{\mathbb H}(K)=\gamma(E).\]
	Thus,
	\[ \gamma^{c }_{\mathbb H}(K) 	= 	\gamma_{\mathbb H}(K)
	=	\gamma(E).	\]
\end{proof}
Younsi\cite{younis} proved the generalization of this result to sets of finite Painlev\'e length and Carath\'{e}odory convergence.

Let $(\Omega_k)_{k\in\mathbb N}$ be a sequence of axially symmetric slice domains in $\widehat{\mathbb H}$ containing $\infty$, and let $\Omega\subset\widehat{\mathbb H}$ be an axially symmetric slice domain containing $\infty$.
For $I\in\mathbb S$, set $\Omega'_{k,I}=\Omega_{k, I}\cap\widehat{\mathbb C}_I,$ and $\Omega'_I=\Omega\cap\widehat{\mathbb C}_I.$
We say that $(\Omega_k)_{k\in\mathbb N}$ converges to $\Omega$ in the quaternionic Carath\'eodory sense if, for some $I\in\mathbb S$, the sequence of  planar domains $(\Omega'_{k,I})$	converges to $(\Omega'_I)$ in the complex Carath\'eodory kernel sense.

Equivalently, $\Omega'_I$ is the Carath\'eodory kernel of
$(\Omega'_{k,I})_{k\in\mathbb N}$ with respect to $\infty$; that is, for every compact set $E\subset\Omega'_I,$ where $\Omega'_I$ is the largest domain containing $\infty$,	there exists $k_0\in\mathbb N$ such that $E\subset\Omega'_{k,I},$ for $ k\geq k_0,$ if such a domain exists. 
The above definition is independent of the choice of the slice $\mathbb C_I$.
For, let $I,J\in\mathbb S$. Since the domains $\Omega_k$ and $\Omega$ are axially symmetric, the map
\[\Phi_{I,J}:\mathbb C_I\longrightarrow\mathbb C_J,\qquad
\Phi_{I,J}(x+Iy)=x+Jy,\]
is an isometry satisfying $\Phi_{I,J}(\Omega'_{k,I})=\Omega'_{k,J}$	and $\Phi_{I,J}(\Omega'_I)=\Omega'_I.$

\noindent Hence the complex Carath\'eodory kernel convergence of
$(\Omega'_{k,I})_{k\in\mathbb N}$ to $\Omega'_I$ is equivalent to the Carath\'eodory kernel convergence of $(\Omega'_{k,J})_{k\in\mathbb N}$ to $\Omega'_J$. Thus, for axially symmetric domains, we have $\Omega_k\xrightarrow{\mathrm{Car}}\Omega$ without specifying a  slice.

Also, here we mean $\Omega$ is non-degenerate m-connected slice-domains means for $I\in \mathbb S$, the planar domain $\widehat{C}_I\setminus E$ has $m$ connected components and none of its connected component of the complement is a singleton. 
In the next theorem we are generalizing Theorem \ref{quaternionic corollary} to sets of finite Quaternionic Painlev\'e length.
\begin{thm}
	Let $K\subset\mathbb H$ be a compact axially symmetric set. Suppose that there exists a sequence $(K_k)_{k\in\mathbb N}$ of compact 	axially symmetric subsets of $K$ satisfying the following properties:
	\begin{enumerate}
		\item Each $K_k$ has finite quaternionic Painlev\'e length;
		
		\item There exists an integer $m\geq1$ such that  the unbounded components of $\widehat{\mathbb H} \setminus K_k $ and $\widehat{\mathbb H}\setminus K$ denoted by $\Omega_k$ and $\Omega$ are non-degenerate $m$-connected slice domains;
		
		\item The sequence $(\Omega_k)_{k\in\mathbb N}$ converges 	to $\Omega$ in the quaternionic Carath\'eodory sense.
	\end{enumerate}
	Then
	\[\gamma^{c }_{\mathbb H}(K)=\gamma_{\mathbb H}(K).\]
\end{thm}

\begin{proof}
	Fix $I\in\mathbb S$ and put $E=K\cap\mathbb C_I,$ $		E_k=K_k\cap\mathbb C_I.$ Since $K$ and $K_k$ are axially symmetric, the corresponding slice domains are
	$\Omega_I = \Omega\cap\mathbb C_I$	and $		\Omega_{k,I} =	\Omega_k\cap\mathbb C_I.$
	
	Now, by the definition of the quaternionic Carath\'eodory convergence, $\Omega_{k,I}\xrightarrow{\mathrm{Car}}
	\Omega_I$ in the complex Carath\'eodory sense in
	$\mathbb C_{\infty,I}$. Also, by assumption, each $K_k$ has finite quaternionic Painlev\'e length. Thus, by Corollary \ref{quaternionic corollary},
	\[\gamma^{c }_{\mathbb H}(K_k)	=\gamma_{\mathbb H}(K_k),
	\qquad k\in\mathbb N.\]
	
	Let $f_k$ and $f$ denote the complex Ahlfors functions associated 	with the domains $\Omega_{k,I}$ and $\Omega_I$, respectively. Since	the domains $\Omega_{k,I}$ are nondegenerate $m$-connected domains
	and converge to $\Omega_I$ in the sense of Carath\'eodory, the complex convergence theorem for Ahlfors functions gives $f_k\longrightarrow f$
	locally uniformly on $\Omega_I$.
	In particular, the derivatives at infinity converge,
	\[\partial_s(f_K)(\infty)\longrightarrow  \partial_s f(\infty).\]
	Since $\partial_s(f_K)(\infty)=\gamma(E_k)$ and $\partial_s f(\infty)=\gamma(E),$ we have $	\gamma(E_k)\longrightarrow\gamma(E).$
	
	Since by Theorem \ref{complexslice quaternionic},
	\[\gamma_{\mathbb H}(K_k)=\gamma(E_k)\ \text{and}\ 		\gamma_{\mathbb H}(K)=\gamma(E).
	\]
	Hence, $\gamma_{\mathbb H}(K_k) \longrightarrow
	\gamma_{\mathbb H}(K).$ Since, $K_k\subseteq K,$
	by the monotonicity of the quaternionic Cauchy capacity we have $\gamma^{c }_{\mathbb H}(K_k)\leq
	\gamma^{c }_{\mathbb H}(K).$
	Now by  Corollary \ref{quaternionic corollary}, we have
	\[\gamma_{\mathbb H}(K_k) = \gamma^{c }_{\mathbb H}(K_k).
	\]
	Hence
	\[\gamma_{\mathbb H}(K_k) \leq 	\gamma^{c }_{\mathbb H}(K).
	\]
	Taking the limit as $k\to\infty$ and using
	$\gamma_{\mathbb H}(K_k) \longrightarrow\gamma_{\mathbb H}(K),$ we obtain
	\[\gamma_{\mathbb H}(K) \leq \gamma^{c }_{\mathbb H}(K).
	\]
	Thus by Lemma \ref{cauchy less analytic}  
	\[\gamma^{c }_{\mathbb H}(K)=\gamma_{\mathbb H}(K).\]
\end{proof}
\section{ Slice-Regular Capacity and Removable Singularity of Slice Regular Functions}
A fundamental problem in complex analysis is to characterize those compact sets across which bounded holomorphic functions defined on the complement of the set can be extended holomorphically.  In the complex setting, this question is closely related to analytic capacity. It is proved that a compact set has zero analytic capacity  when it is removable for bounded holomorphic functions \cite{analytic capacity}.
We now formulate the corresponding quaternionic removability result.
\begin{thm}
	\label{Removability and quaternionic analytic capacity}
	Let $K\subset\mathbb H$ be a compact axially symmetric set, and let $U\subset\mathbb H$ be an axially symmetric slice domain such that $K\subset U$. Then $\gamma_{\mathbb H}(K)=0$
	if and only if every bounded slice-regular function $ 	f:U\setminus K\longrightarrow\mathbb H$	extends to a slice-regular function on $U$.
\end{thm}

\begin{proof}
	Assume that $\gamma_{\mathbb H}(K)=0.$ 	Fix $I\in\mathbb S$, and put $E=K\cap\mathbb C_I,\ U_I=U\cap\mathbb C_I.$ Since  by Theorem \ref{complexslice quaternionic},
	$\gamma_{\mathbb H}(K) 	= 	\gamma_{\mathbb C_I}(E),$ we have $	\gamma_{\mathbb C_I}(E)=0.$
	Thus, $E$ is removable for bounded holomorphic functions by Theorem \ref{ahlfors-removability}.
	
	Now let, $f:U\setminus K\longrightarrow\mathbb H$ be bounded and slice regular. Choose $J\in\mathbb S$ such that $IJ=-JI.$
	Then by  Lemma \ref{splitting lemma}, the restriction of $f$ to $\mathbb C_I$ has the form $f_I(z)=F(z)+G(z)J,\ z\in U_I\setminus E,$ where
	$F,G:U_I\setminus E\longrightarrow\mathbb C_I$	are holomorphic.
	Since $f$ is bounded, the functions $F$ and $G$ are bounded.
	Indeed, the Lemma \ref{splitting lemma} gives $	|F(z)|\leq |f_I(z)|,
	\ |G(z)|\leq |f_I(z)|.$	Hence, since $E$ is removable for bounded holomorphic functions, $F$ and $G$ extend holomorphically to $U_I$. Let their extensions be denoted  by $\widetilde F, \widetilde G.$ 	
	
	Define $\widetilde f_I(z)	= 	\widetilde F(z)+\widetilde G(z)J,\ z\in U_I.$ Then $\widetilde f_I=f_I\ \text{on }U_I\setminus E.$ We now extend $\widetilde f_I$ from the slice $\mathbb C_I$ to the whole domain $U$. For $L\in\mathbb S$, define $$\widetilde f(x+Ly)	=  	\frac12(1-LI)\widetilde f_I(x+Iy)  	+  	\frac12(1+LI)\widetilde f_I(x-Iy).$$	Since $U$ is axially symmetric, whenever $x+Ly\in U$, the points $x+Iy, x-Iy$	belong to $U_I$. Thus the above equation is well-defined. Now, by the representation formula for slice-regular functions, $\widetilde f$ is slice regular on $U$.
	It remains to verify that $\widetilde f$ agrees with $f$ on
	$U\setminus K$.	Let $q=x+Ly\in U\setminus K.$	Since $K$ is axially symmetric, $x+Ly\notin K$ implies $x+Iy\notin E,  	x-Iy\notin E.$ Thus, $\widetilde f_I(x+Iy)=f_I(x+Iy)$ and $
	\widetilde f_I(x-Iy)=f_I(x-Iy).$	Hence,
	\[
	\begin{aligned}
		\widetilde f(x+Ly)
		&=
		\frac12(1-LI)f_I(x+Iy)
		+
		\frac12(1+LI)f_I(x-Iy)\\
		&=f(x+Ly),
	\end{aligned}
	\]
	where the last equality follows from the representation formula applied to the slice-regular function $f$. Thus $\widetilde f|_{U\setminus K}=f.$ Hence, $f$ extends slice regularly to $U$.
	Thus, if $\gamma_{\mathbb H}(K)=0$ then	$K$ is removable for bounded slice-regular functions.
	
	\noindent Conversely, assume that every bounded slice-regular function $f:U\setminus K\longrightarrow\mathbb H$ extends slice regularly to $U$. We shall prove that $\gamma_{\mathbb H}(K)=0$. Fix $I\in\mathbb S$, and put $E=K\cap\mathbb C_I$, $U_I=U\cap\mathbb C_I$. We first show that $E$ is removable for bounded holomorphic functions on $U_I$. Suppose, to the contrary, that $E$ is not removable. Then there exists a bounded holomorphic function $g:U_I\setminus E\longrightarrow\mathbb C_I$ which does not admit a holomorphic extension to $U_I$. 
	Since $\overline{E}=E$, the function $g^*(z) := \overline{g(\overline{z})}$ is also bounded and holomorphic on $U_I\setminus E$. Define $$h_1(z) = \frac{g(z)+g^*(z)}{2}\ \text{and}\ h_2(z) = \frac{g(z)-g^*(z)}{2I}.$$ Then $h_1$ and $h_2$ are holomorphic on $U_I\setminus E$, and $h_k(\overline{z})=\overline{h_k(z)}$ for $k=1,2$. Also $g = h_1 + Ih_2$. At least one of $h_1$ and $h_2$ cannot extend holomorphically to $U_I$. (If both extended, it would give a holomorphic extension of $g$, yielding a contradiction). Suppose that $h_1$ does not extend. Let, for $z=x+Iy$, $$h_1(x+Iy) = u(x,y) + Iv(x,y),$$ where $u(x,y)$ and $v(x,y)$ are real-valued. Since, $ h_1$ is holomorphic and $h_1(\overline{z})=\overline{h_1(z)}$, we have 
	\begin{equation}\label{intrinsic}
		u(x,-y) = u(x,y),\ v(x,-y) = -v(x,y).
	\end{equation}
	For $q = x+Ly$ with $L\in\mathbb{S}$, define $$H(x+Ly) = u(x,y) + Lv(x,y).$$ Then by \eqref{intrinsic}, it is well-defined.   Hence $H:U\setminus K\to\mathbb{H}$ is slice regular and $$|H(q)|^2 = u(x,y)^2 + v(x,y)^2 = |h_1(x+Iy)|^2.$$ Hence $$\|H\|_\infty = \|h_1\|_\infty \le \|g\|_\infty < \infty.$$ Thus, $H$ is a bounded slice-regular function on $U\setminus K$.  Then, $H$ extends to a slice-regular function  $\tilde{H}:U\to\mathbb{H}$ by assumption. Since the restriction of a slice-regular function to a slice is holomorphic, $\tilde{H}_I = \tilde{H}|_{U_I}$ is holomorphic on $U_I$. But on $U_I\setminus E$, we have $\tilde{H}_I = H|_{\mathbb{C}_I}$. Then, for $L=I$, $$H(x+Iy) = u(x,y) + Iv(x,y) = h_1(x+Iy).$$  Thus, $ \tilde{H}_I$ is a holomorphic extension of $h_1$ from $U_I\setminus E$ to $U_I$, which is a contradiction, since $h_1$ was not extendable. Hence our assumption is wrong. Thus, $E$ is removable for bounded holomorphic function. Hence, $\gamma_{\mathbb{C}}(E) = 0$ by Theorem \ref{ahlfors-removability}. Now,  since $ \gamma_{\mathbb{H}}(K) = \gamma_{\mathbb{C}}(E)$, we have $\gamma_{\mathbb{H}}(K) = 0$. 
\end{proof}
A central ingredient in Tolsa's work is the relationship between
the Cauchy transform and the Menger curvature of measures. For
three distinct points $z_1,z_2,z_3\in\mathbb C$, let
$c(z_1,z_2,z_3)$ denote the reciprocal of the radius of the circle passing through these points. If the points are collinear, set $c(z_1,z_2,z_3)=0.$
The  total Menger curvature of  a positive finite measure $\mu$, is
\[c^2(\mu)=\iiint c(z_1,z_2,z_3)^2\, d\mu(z_1)\,d\mu(z_2)\,d\mu(z_3).\]
See, \cite{tolsa, melnikov}.
The curvature integral measures the geometric distribution of
the measure. Melnikov's curvature identity connects this
quantity with the $L^2(\mu)$ behaviour of the Cauchy transform.
This connection made it possible to characterize analytic
capacity in terms of measures satisfying linear growth and
finite curvature. For a positive measure $\mu$, linear growth means that there exists a constant $C_0>0$ such that
$\mu(B(z,r))\leq C_0r$ for every $z\in\mathbb C$ and $r>0$.
A main result in the theory of analytic capacity is Tolsa's characterization of non-removable sets in terms of positive measures with finite linear growth and finite curvature.
\begin{thm}\cite{finitetolsa}
	\label{tolsa-curvature}
	Let $E\subset\mathbb C$ be compact. Then $\gamma_{\mathbb C}(E)>0$ if and only if there exists a nonzero Borel measure $\mu$ supported on $E$ such that $\mu(B(z,r))\leq C_0r$ for all $z\in\mathbb C$ and $r>0$, and
	$	c^2(\mu)<\infty.$
\end{thm}

The Painlev\'e problem asks for a geometric characterization of
compact sets that are removable for bounded holomorphic functions. A compact set $E\subset\mathbb C$ is removable if every bounded function holomorphic on $\mathbb C\setminus E$ extends holomorphically across $E$.

The connection between removability and analytic capacity is
expressed by $E$  is removable if and only if $\gamma(E)=0.$
Tolsa proved these results on removability and its correspondence with $\gamma_{+}$ capacity.
\begin{thm}\label{tolsa1.2}\cite{tolsa}
	A compact set $E \subset \mathbb{C}$ is non removable for bounded analytic functions if and only if it supports a positive Radon measure with linear growth and finite curvature.
\end{thm}
\begin{thm}\label{tolsa1.3}\cite{tolsa}
	A compact set $E \subset \mathbb{C}$ is removable for bounded analytic functions if and only if there exists a finite positive Radon measure $\mu$ on $\mathbb{C}$ such that for all $x \in E$ either $\Theta_{\mu}^*(x) = \infty$ or $c_{\mu}^2(x) = \infty$.
\end{thm}
This result gives a geometric characterization of positive analytic capacity. It shows that the existence of a measure with linear growth and finite Menger curvature detects whether the compact set has positive analytic capacity. This result provides a  connection between removability, analytic capacity, and the geometric properties of measures supported on the set.

The purpose of this section is to establish the corresponding result in the quaternionic setting. For a compact axially symmetric set $K$, we reduce the question of removability for bounded slice regular functions to the corresponding planar removability problem on a complex slice. It allows us to apply Tolsa's result \cite{tolsa} and then lift the resulting measure from the slice to the quaternionic setting.

The main issue is to formulate the geometric conditions in a quaternionic way. In particular, the classical linear growth condition is replaced by spherical linear growth, while the Menger curvature is replaced by the corresponding spherical curvature associated with triples of quaternions. 

Let $K\subset\mathbb H$ be a compact axially symmetric set.
For $q,s\in\mathbb H$, define the spherical distance between the spherical classes $[q]$ and $[s]$ by
\[
d_{\mathbb S}([q],[s])
=
\inf\{|u-v|:u\in[q],\ v\in[s]\}.
\]
Equivalently,
\[
d_{\mathbb S}([q],[s])
=
\sqrt{
	\big(\operatorname{Re}q-\operatorname{Re}s\big)^2
	+
	\big(|\operatorname{Im}q|-|\operatorname{Im}s|\big)^2
}.
\]
In particular, if $q,s\in\mathbb C_I$, then
\[
d_{\mathbb S}([q],[s])
=
\min\{|q-s|,|q-\overline{s}|\}.
\]
For $q,s\in K$, we say that $q$ and $s$ are spherically
	equivalent if $[q]=[s]$. Let $\tilde{K}$ denotes the quotient space of $K$ with respect to this equivalence relation; that is, its elements are the spherical equivalence classes
\[\tilde{K} = \{[q]:q\in K\}\]
Thus, each point of $\tilde{K}$ represents one spherical class in $K$. For $[q]\in \tilde{K}$ and $r>0$, define
\[B_{\mathbb S}([q],r) = \{s\in \tilde{K}:d_{\mathbb S}([q],[s])<r\}.\]
and \[B_{\mathbb S}^K([q],r) = \{s\in K:d_{\mathbb S}([q],[s])<r\}.\]
Equivalently, for	$I\in\mathbb S$ and $E=K\cap\mathbb C_I.$
\[	B_{\mathbb S}^{E}([q],r) =	\bigl(B_{\mathbb C_I}(q,r)
\cup	B_{\mathbb C_I}(\overline q,r)\bigr)\cap E,	\]
where $B_{\mathbb C_I}$ is the ordinary Euclidean ball in the slice $\mathbb C_I$.

Let $\pi:K\longrightarrow \tilde{K},  \pi(q)=[q],$
be the  projection onto the spherical classes. If $\mu$ is a	positive finite Radon measure supported on $K$, we denote by $\nu(A)=\mu\bigl(\pi^{-1}(A)\bigr)$ be the measure on the space of spherical classes
for every Borel set $A\subset \tilde{K}$.
Equivalently, for every bounded Borel function $\varphi$ on
$\tilde{K}$,
\[
\int_{\tilde{K}}\varphi([q])\,d(\nu([q]))
=
\int_K\varphi(\pi(q))\,d\mu(q).
\] 
Now, the upper spherical density of $\mu$ is defined by
\[	\Theta_{\mathbb H,\mu}^{*}([q])
=
\limsup_{r\to0}
\frac{\nu\bigl(B_{\mathbb S}([q],r)\bigr)}{r}.
\]
Equivalently,
\[
\Theta_{\mathbb H,\mu}^{*}([q])
=
\limsup_{r\to0}
\frac{\mu\bigl(B_{\mathbb S}^{K}([q],r)\bigr)}{r}.
\]

We say that $\mu$ has spherical linear growth if there exists a
constant $C>0$ such that 
\[	\nu\bigl(B_{\mathbb S}([q],r)\bigr)\leq Cr,
\qquad [q]\in \tilde{K},\quad r>0,\]
Equivalently,
\[	\mu\bigl(B_{\mathbb S}^{K}([q],r)\bigr)\leq Cr,
\qquad [q]\in \tilde{K},\quad r>0.\]
For $z_1,z_2,z_3\in\mathbb C_I$, let
$c(z_1,z_2,z_3)$ denote the complex Menger curvature. We define the spherical Menger curvature of three spherical classes by
\begin{equation}
	\label{spherical menger final}
	c_{\mathbb H}([q_1],[q_2],[q_3])^2
	=
	\frac18
	\sum_{\varepsilon_1,\varepsilon_2,\varepsilon_3\in\{0,1\}}
	c\left(
	q_1^{\varepsilon_1},
	q_2^{\varepsilon_2},
	q_3^{\varepsilon_3}
	\right)^2,
\end{equation}
where $q_i^0=q,\	q_i^1=\overline q.$
The factor $1/8$ is chosen so that the spherical curvature agrees 	with the complex Menger curvature after passing to a
conjugation-symmetric measure on a fixed complex slice.
For $\nu$, the spherical measure we define the  spherical curvature by
\begin{equation}
	\label{local spherical curvature final}
	c_{\mathbb H,\mu}^{\,2}([q])
	=
	\iint_{\tilde{K}}
	c_{\mathbb H}([q],[s],[t])^2
	\,d\nu([s])\,d\nu([t]),
\end{equation}
and the total spherical curvature by
\begin{equation}
	\label{total spherical curvature final}
	c_{\mathbb H}^{\,2}(\mu)
	=
	\int
	c_{\mathbb H,\mu}^{\,2}([q])
	\,d\nu([q]).
\end{equation}	
Next lemma establishes a correspondence between measures on the space of spherical equivalence classes and conjugation invariant measures on a complex slice. We construct a planar measure which recovers the original spherical measure. This result allows us to relate quaternionic problems to their planar one.
\begin{lem}
	\label{symmetric slice lift final}
	Let $K\subset\mathbb H$ be compact and axially symmetric, and let $E=K\cap\mathbb C_I.$
	Let $\nu$ be a positive finite Radon measure on the spherical classes
	of $K$. Put	$E^+ = \{x+Iy\in E:\ y\geq0\}.	$
	Then there exists a conjugation-invariant positive finite Radon measure $\lambda$ supported on $E$ such that
	\begin{equation}
		\label{slice lift final}
		\int_E\varphi(z)\,d\lambda(z)
		=
		\int_{E^+}
		\frac{\varphi(z)+\varphi(\overline z)}{2}
		\,d\nu([z])
	\end{equation}
	for every bounded Borel function $\varphi$ on $E$.
	Moreover, $	\lambda \circ\pi^{-1}=\nu, \pi(z)=\pi(\overline z)= [z]$ where $\pi$ be the projection from $E$ on to $\tilde{K}$.
\end{lem}

\begin{proof}
	For a Borel set $A\subset E$,  define,
	\[\lambda(A) = 	\int_{E^+}	\frac{\chi_A(z)+\chi_A(\overline z)}{2}\,d\nu([z]),\]
	where $\chi_{A_j}$ is the characteristic function. Then $\lambda$ is a positive finite Radon measure on $E$.
	Since  for a simple function $\phi$,
	\[\phi(z) = \sum_{j=1}^{n} a_j\chi_{A_j}(z),\] for $A_j \subset E$,
	we have $$\int_E\varphi(z)\,d\lambda(z)
	=
	\int_{E^+}
	\frac{\varphi(z)+\varphi(\overline z)}{2}
	\,d\nu([z])$$ 
	
	By replacing $A$ by $\overline A$ yields $		\lambda(\overline A)=\lambda(A).$
	Hence, $\lambda \circ\pi^{-1}=\nu.$
	Thus the measure $\lambda$ is  conjugation invariant and represents the measure $\nu$ on the spherical classes.
\end{proof}
Using this lemma we proved the equality of linear growth and spherical linear growth.
\begin{prop}\label{finitelength}
	Let $K\subset\mathbb H$ be a compact axially symmetric set. Fix $I\in\mathbb S$ and $E=K\cap\mathbb C_I.$
	Let $\lambda$ be a positive finite Radon measure on $E$ and let $\mu$ be positive finite Radon measure  on $K.$
	
	Then $\lambda$ has  linear growth if and only if
	$\mu$ has  spherical linear growth. More precisely,
	\[\lambda(B_{\mathbb C_I}(z,r)\cap E)\le Cr,\qquad z\in E,\ r>0,\]
	if and only if
	\[
	\nu(B_{\mathbb S}([q],r))\le C'r,
	\qquad [q]\in\widetilde K,\ r>0.
	\]
	where $\nu$ be the positive finite Radon measure on $\tilde{K}$.
\end{prop}

\begin{proof}
	Let $\lambda$ be a positive finite Borel measure on $E$ and let $\pi:E\longrightarrow\widetilde K,\ \pi(z)=[z]$	be the quotient map. Define the spherical measure $\nu$ on $\widetilde K$ by
	$\nu(A)=\lambda(\pi^{-1}(A))$
	for every Borel set $A\subset\widetilde K$.
	For $q,s\in\mathbb C_I$, the spherical distance
	is given by
	\[d_{\mathbb S}([q],[s])=\min\{|q-s|,|q-\overline{s}|\}.\]
	Hence, $d_{\mathbb S}([q],[s])<r$
	if and only if $|q-s|<r$ or $|q-\overline{s}|<r.$
	Since $	|q-\overline{s}|=|\overline q-s|,$
	we have
	\begin{equation*}\label{eq:spherical-ball-preimage}
		\pi^{-1}(B_{\mathbb S}([q],r))
		= 	\big(B_{\mathbb C_I}(q,r)
		\cup B_{\mathbb C_I}(\overline q,r)\big)\cap E.
	\end{equation*}
	First, suppose that $\lambda$ has finite linear growth.
	Then there exists a constant $C>0$ such that
	\[
	\lambda(B_{\mathbb C_I}(z,r)\cap E)
	\le Cr,
	\qquad z\in E,\quad r>0.
	\]
	Then,
	\[
	\begin{aligned}
		\nu(B_{\mathbb S}([q],r))
		&=\lambda\big(\pi^{-1}(B_{\mathbb S}([q],r))\big)\\
		&=\lambda\left(
		\big(B_{\mathbb C_I}(q,r)
		\cup B_{\mathbb C_I}(\overline q,r)\big)\cap E
		\right)\\
		&\le
		\lambda(B_{\mathbb C_I}(q,r)\cap E)
		+
		\lambda(B_{\mathbb C_I}(\overline q,r)\cap E)\\
		&\le Cr+Cr =2Cr.
	\end{aligned}
	\]
	Thus, $\nu(B_{\mathbb S}([q],r))\le 2Cr,$
	which shows that $\mu$ has finite spherical linear growth.
	
	Conversely, suppose that  $\mu$ supported on $K$ such that  $\mu$ has spherical linear growth. Then there exist spherical measure $\nu$ supported on $\tilde{K}$ such that  $\nu(B_{\mathbb S}([q],r))\le 2Cr.$  By Lemma \ref{symmetric slice lift final}, let $\lambda$ be the
	corresponding conjugation-symmetric positive finite Radon measure	supported on $E=K\cap\mathbb C_I.$		
	For $q\in E$, we have 
	\[B(q,r)\cap E 	\subset B_{\mathbb S}^{K}([q],r).\]
	Thus,
	\[\lambda(B(q,r)) \le 	\mu(B_{\mathbb S}^{K}([q],r))
	\le Cr.\]
	Thus $\lambda$ has finite linear growth.
	Hence, $\lambda$ has finite linear growth if and only if
	$\nu$ has finite spherical linear growth.
\end{proof}
Next lemma establish a relationship between the spherical upper density of a quaternionic measure and the upper density of its associated slice measure. 
\begin{lem}
	\label{density correspondence final}
	Let $\mu$ be a positive finite Radon measure supported on $K$,  $\nu$ be the spherical class measure
	and let $\lambda$ be the conjugation-symmetric slice measure from $\nu$ given by Lemma \ref{symmetric slice lift final}. Then, for every $z\in E$,
	\[\Theta_{\mathbb H,\mu}^{*}([z])=\infty\ \text{if and only if}\ \Theta_{\lambda}^{*}(z)=\infty,\]
	where $\Theta_{\lambda}^{*}(z) 	= \limsup_{r\to0}
	\frac{\lambda(B(z,r))}{r}.$
\end{lem}

\begin{proof}
	Suppose  that $z\in E\setminus\mathbb R$. We have
	\[B_{\mathbb S}^{E}([z],r) 	= 	\bigl(B(z,r)\cup B(\overline z,r)\bigr)\cap E.\]
	Since $Im(z)\neq 0$, we can always choose sufficiently small $r$. But, for sufficiently small $r>0$, the two balls are disjoint. Since $\lambda$ is invariant under conjugation,
	\[\lambda(B(z,r)) = \lambda(B(\overline z,r)).\]
	Thus, \[\nu(B_{\mathbb S}([z],r)) 	= 	2\lambda(B(z,r)).
	\]
	Hence
	\[\Theta_{\mathbb H,\mu}^{*}([z]) =  	2\Theta_{\lambda}^{*}(z),
	\]
	and the result follows.
	
	If $z\in\mathbb R$, then $[z]=\{z\}$ and
	$\nu(B_{\mathbb S}([z],r)) = \lambda(B(z,r)).$
	Thus, $\Theta_{\mathbb H,\mu}^{*}([z]) = 	\Theta_{\lambda}^{*}(z).$
\end{proof}
Next we establish a connection between spherical and classical Menger curvature.  For every $z\in E$, we show that the spherical curvature associated with the conjugacy class $[z]$ coincides with the planar Menger curvature associated with $z$.  Since the measure $\lambda$ is conjugation invariant we can express the planar curvature integral as an average over the conjugate representatives of each point.
\begin{lem}
	\label{curvature correspondence final}
	For every $z\in E$,
	$$c_{\mathbb H,\mu}^{\,2}([z]) = c_{\lambda}^{\,2}(z),$$
	where, $c_{\lambda}^{\,2}(z) = \iint c(z,w,t)^2\,d\lambda(w)d\lambda(t).$
	Then,
	$$c_{\mathbb H}^{\,2}(\mu)	= c^2(\lambda).$$
\end{lem}

\begin{proof}
	Since $\lambda$ is conjugation invariant, for every positive Borel function $\psi$ on $E\times E$,
	\[\iint\psi(w,t)\,d\lambda(w)d\lambda(t)=
	\frac14 \sum_{\varepsilon_2,\varepsilon_3\in\{0,1\}}
	\iint \psi(w^{\varepsilon_2},t^{\varepsilon_3})
	\,d\nu([w])d\nu([t])\]
	by defining $\lambda$ as in Lemma \ref{slice lift final} where $w^0 =w,\ w^1 =\bar{w}$.
	By taking, $\psi(s,t)=c(q,w,t)^2$ for fixed $z$ we get
	\begin{equation}
		\label{planar curvature expansion final}
		c_{\lambda}^{\,2}(z) =
		\frac14 	\sum_{\varepsilon_2,\varepsilon_3\in\{0,1\}}
		\iint			c(z,w^{\varepsilon_2},t^{\varepsilon_3})^2
		\,d\nu([w])d\nu([t]).
	\end{equation}
	On the other hand, from
	\eqref{spherical menger final} and
	\eqref{local spherical curvature final},
	\begin{align*}
		c_{\mathbb H,\mu}^{\,2}([z])
		&=
		\frac18
		\sum_{\varepsilon_1,\varepsilon_2,\varepsilon_3\in\{0,1\}}
		\iint
		c(z^{\varepsilon_1},
		w^{\varepsilon_2},
		t^{\varepsilon_3})^2
		\,d\nu([w])d\nu([t]).
	\end{align*}
	Since the Menger curvature is invariant under complex conjugation, the contributions corresponding to $\varepsilon_1=0$ and
	$\varepsilon_1=1$ are equal. Thus,
	\begin{align*}
		c_{\mathbb H,\mu}^{\,2}([q])
		&=
		\frac14
		\sum_{\varepsilon_2,\varepsilon_3\in\{0,1\}}
		\iint
		c(q,s^{\varepsilon_2},t^{\varepsilon_3})^2
		\,d\nu([s])d\nu([t]).
	\end{align*}
	Comparing this expression with
	\eqref{planar curvature expansion final}, we obtain
	\[
	c_{\mathbb H,\mu}^{\,2}([z])
	=
	c_{\lambda}^{\,2}(z).
	\]
	Integrating with respect to the corresponding measure gives
	\[
	c_{\mathbb H}^{\,2}(\mu)
	=
	c^2(\lambda).
	\]
\end{proof}	
From all these results we have this corollary for a positive finite Radon measure supported on $K$.
\begin{cor}
	\label{spherical planar transfer final}
	Let $K\subset\mathbb H$ be compact and axially symmetric, and let $E=K\cap\mathbb C_I.$ Let $\mu$ be a positive finite Radon measure supported on $K$.
	Then its spherical-class measure admits  a planar measure $\lambda$ supported on $E$ such that
	\[
	\Theta_{\mathbb H,\mu}^{*}([q])=\infty\ \text{if and only if}\ \Theta_{\lambda}^{*}(q)=\infty
	\]
	for every $q\in E$, and $c_{\mathbb H,\mu}^{\,2}([q])	=c_{\lambda}^{\,2}(q).$
	In particular,
	\begin{equation}
		\label{eq:total-transfer-final}
		c_{\mathbb H}^{\,2}(\mu)<\infty
		\ \text{if and only if}\ 
		c^2(\lambda)<\infty.
	\end{equation}
\end{cor}
\begin{proof}
	From Lemma \ref{symmetric slice lift final},
	Lemma \ref{density correspondence final} and 
	Lemma \ref{curvature correspondence final} the required result follows.
\end{proof}	
Tolsa \cite{finitetolsa} proved that for a compact set in $\mathbb{C}$, analytic capacity is positive if and only if there exists a positive finite Radon measure supported on the compact set with linear growth and finite curvature.
The following theorem shows that non-removability of $K$ for bounded slice regular functions yields a nontrivial positive measure on $K$ satisfying these quaternionic geometric estimates which is the analogue of [Theorem 1.2 \cite{tolsa}.]
\begin{thm}\label{remotolsa}
	Let $K\subset\mathbb H$ be a compact axially symmetric set. Then $K$ is non-removable for bounded slice regular functions if and only if there exists a nonzero positive finite Radon measure $\mu$ supported on $K$ such that $\mu$ has spherical linear growth and finite total spherical curvature.
\end{thm}

\begin{proof}		
	Suppose that $K$ is non-removable. Then by Theorem \ref{Removability and quaternionic analytic capacity}, $\gamma_{\mathbb H}(K)>0$. Hence by Theorem \ref{complexslice quaternionic}, $\gamma_{\mathbb H}(E)>0$. Then there exists a positive Radon measure $\lambda$ on $E$ such that $\lambda$ has finite linear growth and finite total spherical curvature. Then by Proposition \ref{finitelength} and Corollary \ref{spherical planar transfer final} there exist a spherical measure $\nu$ such that $\nu$ has finite spherical linear growth and finite total spherical curvature.
	
	\noindent Similarly we can prove the converse part. 
\end{proof}
The next theorem is a dual of Theorem \ref{remotolsa}.
\begin{thm}
	Let $K\subset\mathbb H$ be a compact axially symmetric set and fix $I\in\mathbb S$. Put $E=K\cap\mathbb C_I.$
	Then $K$ is removable for bounded slice regular functions if and only if there exists a finite positive Radon measure $\mu$ on $\mathbb H$ such that, for every $q\in K$,
	\[\Theta_{\mathbb H,\mu}^{*}([q])=\infty
	\qquad\text{or}\qquad
	c_{\mathbb H,\mu}^{\,2}([q])=\infty.
	\]
\end{thm}
\begin{proof}
	The required result follows from Lemma \ref{density correspondence final} and Lemma \ref{curvature correspondence final}.
\end{proof}
\part*{\MakeUppercase{Part II:Quaternionic Rational Approximation}}
\section{Approximation of Slice Regular Functions in Compact Sets}
Let $K\subset \mathbb{H}$ be a  compact set. Let $\mathcal C(K, \mathbb H)$ be the set of all continuous quaternion valued functions  and $\mathcal C_L(K)$ be the set of continuous slice functions on  $K$ respectively. Since every continuous slice functions are continuous functions, we have $\mathcal C_L(K) \subseteq \mathcal C(K, \mathbb H)$. But not every continuous function need not be continuous slice function, for example consider, $f: \mathbb{H} \to \mathbb{H}$ defined by $f(q) = x_1$, where $q = x_0 + x_1 i + x_2 j + x_3 k$.

Now we define,
\[\widetilde{\mathcal A}_L(K)=\left\{f\in \mathcal C(K, \mathbb H):
f\text{ is slice regular on }K^\circ\right\},
\]
and
\[\mathcal A_L(K)=\left\{f\in \mathcal C_L(K):
f\text{ is slice regular on }K^\circ\right\}.
\]
Also, $\mathcal{R}_L(K)$ is the set of all functions on $K$ uniformly approximable by left rational functions with poles off $K$.

\noindent Then the following holds from the definition:
\begin{lem}\label{approxsubset}
	Let $K \subset \mathbb{H}$ be compact and axially symmetric. Then 
	\begin{enumerate}
		\item $\mathcal{R}_L(K) \subset \mathcal C(K, \mathbb H)$,
		\item $\mathcal{R}_L(K) \subset \mathcal C_L(K)$,
		\item $\mathcal{R}_L(K) \subseteq  \mathcal A_L(K)$,
		\item $\mathcal{R}_L(K) \subseteq \widetilde{\mathcal A}_L(K)$.
	\end{enumerate}
\end{lem}

\begin{proof}
	Since every slice rational functions are slice regular functions, it follows.
\end{proof}
Next lemma is immediate:
\begin{lem}\label{alcl}
	Let $K\subseteq\mathbb{H}$ be a compact set. 
	Then
	\[\widetilde{\mathcal A}_L(K)=\mathcal C(K, \mathbb H)\ \text{and } {\mathcal A}_L(K)=\mathcal C_L(K)\] if and only if $K^\circ=\varnothing.$
\end{lem}

\begin{proof}
	Suppose that $K^\circ=\varnothing$.
	Then by the definition of $\widetilde{\mathcal A}_L$  every continuous quaternion-valued function on $K$ belongs to $\widetilde{\mathcal A}_L(K)$. Thus,
	\[\widetilde{\mathcal A}_L(K)=\mathcal C(K, \mathbb H).\]
	
	Suppose that $K^\circ\neq\varnothing$.
	Consider the function $	f:K\longrightarrow\mathbb{H},\ f(q)=\overline{q}.$	Then, $f\in \mathcal C(K, \mathbb H)$. 
	We show that $f$ is not slice regular on
	$K^\circ$. For any $I\in\mathbb{S}$, write
	$q=x+Iy$. Then $f_I(x+Iy)=x-Iy.$ Thus,
	\[\frac12\left(\frac{\partial}{\partial x} +I\frac{\partial}{\partial y}\right)f_I(x+Iy)
	=\frac12\bigl(1+I(-I)\bigr)=1\neq 0.\]
	Thus $f$ is not slice regular on any nonempty open subset
	of $\mathbb{H}$. In particular,
	$f\notin \widetilde{\mathcal A}_L(K)$. Hence,
	\[\widetilde{\mathcal A}_L(K)\neq \mathcal C(K, \mathbb H).\]
\end{proof}
The quaternionic analogue of Runge’s Theorem for  slice regular functions on an axially symmetric compact set $K$ was proved by Fabrizio Colombo, Irene Sabadini, and Daniee C. Struppa in \cite{runge quaternionic};
\begin{thm}\cite{runge quaternionic}
	Let $A$ be an axially symmetric set such that $A\cap C\neq\varnothing$ 	for every connected component $C$ of
	$(\mathbb H\cup\{\infty\})\setminus K$. If $f$ is left slice regular on an axially symmetric open set
	$U\supset K$, then, for every $\varepsilon>0$, there exists a left rational function $r$ whose poles lie in $A$ such that
	\[
	\sup_{q\in K}|f(q)-r(q)|<\varepsilon.
	\]
\end{thm}
Thus under the hypotheses of Quaternionic Runge's theorem, if $f$ is slice regular in an axially symmetric open neighborhood of compact set $K$, then the restriction $f|_K$ belongs to $\mathcal R_L(K)$.
Also, if $\mathcal O_L(K)$ denote the set of all left slice
regular functions on an open neighborhood of $K$ and $\mathcal R_L(K)$ denotes the uniform closure of the left rational functions whose poles satisfy the hypotheses of Runge's theorem, then it is clear that $\overline{\mathcal O_L(K)}=\mathcal R_L(K).$

Runge's theorem permits rational approximation with poles in
appropriate components of the complement.  The following statements are the standard Runge-type result presented by Gal and Sabadini.
\begin{thm}\cite{quatern approx}
	\label{qrunge-polynomial}
	Suppose that $\widehat{\mathbb H}\setminus K$ and $\widehat{\mathbb C}\setminus (K \cap _I\mathbb C)$ is connected for every $I \in \mathbb{S}$. Then every slice regular on an axially symmetric open neighborhood of $K$ can be approximated by a  quaternionic polynomial.
\end{thm}

\noindent
When the complement has the required connectedness, the rational
approximants can be replaced by polynomials. The theorem is proved by reducing the approximation problem to a complex slice and using the representation formula to recover the quaternionic function.

Now we consider how the approximation works if $\widehat{\mathbb H}\setminus K$ or if $\widehat{\mathbb C}_I \setminus (K\cap \mathbb{C}_I)$ is connected.   Next is the quaternionic analogue of Lavrentiev approximation theorem \cite{laverentiev}:
\begin{thm}
	Let $K\subset\mathbb H$ be a compact, axially symmetric,
	nowhere-dense set such that $\Omega=\mathbb H\setminus K$	is an axially symmetric slice domain. Then
	\[	\overline{\mathbb H[q]}^{\, \widetilde{\mathcal C}_L (K)}=\mathcal{C}_L(K),\]
	where, 
	\[\mathbb H[q] = \left\{ \sum_{n=0}^{N}q^n a_n:
	N\geq 0,\ a_n\in\mathbb H\right\}.\]
\end{thm}

\begin{proof}
	Since every left slice-regular polynomial is a continuous left slice function on $K$ and the class of continuous left slice
	functions is closed under uniform convergence we have
	\[\overline{\mathbb H[q]}^{\,\mathcal C(K, \mathbb H)}\subseteq \mathcal C_L(K).\]
	
	Now, fix $I\in\mathbb S$ and put $E=K\cap\mathbb C_I.$
	Since $K$ is axially symmetric and nowhere dense,
	$E^\circ=\varnothing.$
	
Also, $\mathbb C_I\setminus E=(\mathbb H\setminus K)\cap\mathbb C_I	$	is connected since $\Omega$ is an axially symmetric slice
domain. Hence, by Theorem \ref{lavrentiev},
\[\overline{\mathbb C_I[z]}^{\,\mathcal C(E, \mathbb C_I)}
=\mathcal C(E, \mathbb C_I).\]

Let $f\in \mathcal{C}_L(K)$, and choose $J\in\mathbb S$ with $J\perp I$. 	Then by Lemma \ref{contslice}, there exist $g_1,g_2\in \mathcal C(E, \mathbb C_I)$ such that
\[f(z)=g_1(z)+g_2(z)J,\qquad z\in E.\]

Given $\varepsilon>0$, choose $p_1,p_2\in\mathbb C_I[z]$ such that
\[
\|g_1-p_1\|_E<\varepsilon,
\qquad
\|g_2-p_2\|_E<\varepsilon.
\]
Define $P(q)=p_1(q)+p_2(q)J.$	Then $P$ is a left slice-regular polynomial and $\|f-P\|_E<2\varepsilon.$

Set $h=f-P$. Then by Theorem \ref{representation formula}, for every $x+Ly\in K$,
\[
\begin{aligned}
	h(x+Ly)
	={}&\frac12\bigl(h(x+Iy)+h(x-Iy)\bigr)+\frac{LI}{2}
	\bigl(h(x-Iy)-h(x+Iy)\bigr).
\end{aligned}
\]
Axial symmetry ensures that $x+Iy,x-Iy\in E$. Hence,
\[
\begin{aligned}
	|h(x+Ly)|
	&\leq |h(x+Iy)|+|h(x-Iy)|<4\varepsilon.
\end{aligned}
\]
Thus, $\|f-P\|_K\leq4\varepsilon.$ Since $\varepsilon>0$ is arbitrary, we have
\[
\mathcal C_L(K)\subseteq
\overline{\mathbb H[q]}^{\,\widetilde{\mathcal C}(K)}.
\]

Hence, $\overline{\mathbb H[q]}^{\,\widetilde{\mathcal C}(K)}=\mathcal C_L(K).$
\end{proof}
Thus, if $K^\circ=\varnothing$,  every function in $\mathcal C_L(K)$ can be approximated uniformly on $K$ by polynomials. Also by Lemmas \ref{approxsubset} and \ref{alcl} we have 
\[\mathcal{C}_L(K) = \mathcal{R}_L(K) = \mathcal{A}_L(K).\]
But this doesn't holds in general in the case of all functions in $\mathcal C(K, \mathbb H)$ even when the compact set is axially symmetric and nowhere dense and its complement is an axially symmetric slice domain. For example,

\begin{exa}\label{noapprox}

Let $K=\{q\in\mathbb H:\operatorname{Re}(q)=0,\ 1\leq |q|\leq 2\}.$
That is, $K=\{Iy:I\in\mathbb S,\ 1\leq y\leq 2\}.$
Then $K$ is compact and axially symmetric. Also, $	K^\circ=\varnothing,$
so $K$ is nowhere dense in $\mathbb H$.

Also $\Omega=\mathbb H\setminus K$	is open, axially symmetric, and connected. For, every point of 	$\Omega$ with nonzero real part can be joined to the real axis without meeting $K$. Also, points on the imaginary axis that belong to 
$\Omega$ can also be joined to the real axis by paths avoiding $K$. Thus $\Omega$ is an axially symmetric slice domain.

Fix a quaternionic imaginary unit $i$, define,
\[f:K\longrightarrow\mathbb H,\qquad
f(q)=\frac{q}{|q|}\,i\,\frac{q}{|q|}.
\]
Since $|q|\geq 1$ on $K$, the function $f$ is well-defined and
continuous.

We claim that $f$ cannot be uniformly approximated on $K$ by
left slice-regular polynomials.

Suppose $P(q)=\sum_{n=0}^{N}q^n a_n,
\ a_n\in\mathbb H,$ is a left slice-regular polynomial. For every $I\in\mathbb S$ and 	every $r\in[1,2]$, we have
\[
\begin{aligned}
	P(Ir)+P(-Ir)
	&=\sum_{n=0}^{N}
	\bigl((Ir)^n+(-Ir)^n\bigr)a_n\\
	&=2\sum_{2k\leq N}(-1)^k r^{2k}a_{2k}.
\end{aligned}
\]
The expression on the right is independent of $I$. Thus,
\[
P(Ir)+P(-Ir)=P(Jr)+P(-Jr)
\]
for every $I,J\in\mathbb S$ and every $r\in[1,2]$. If a sequence of left slice-regular polynomials converges uniformly
to a function $g\in \mathcal C(K, \mathbb H)$, passing to the limit gives
\[
g(Ir)+g(-Ir)=g(Jr)+g(-Jr).
\]
Hence every function in the uniform closure of the left
slice-regular polynomials satisfies this identity.

However, for the function $f$ defined above, we have
\[
f(Ir)=IiI,\qquad f(-Ir)=IiI.
\]
Taking $I=i$ gives, $f(ir)+f(-ir)=2iii=-2i,$
whereas taking $J=j$, gives $f(jr)+f(-jr)=2jij=2i,$ which is a contradiction. So the function $f$ does not satisfy the necessary identity. Hence, $f\notin
\overline{\mathbb H[q]}^{\,\mathcal C(K, \mathbb H)},$	and thus
\[\overline{\mathbb H[q]}^{\,\mathcal C(K, \mathbb H)}
\neq \mathcal C(K, \mathbb H).\]
\end{exa}

This example shows that the classical Lavrentiev approximation
does not extend to all continuous quaternion-valued
functions. In particular, the failure occurs even $K$ is
compact, axially symmetric, and nowhere dense and that
$\mathbb H\setminus K$ is an axially symmetric slice domain.

In the above theorem we considered nowhere dense set. In the following theorem we are considering the sets which satisfy $K=\overline{K^\circ}$, which is a quaternionic analogue of  approximation theorem by Keldysh \cite{keldysh}. Keldysh proved that $\mathbb{C}\setminus E$ is connected and $E=\bar{E^\circ}$, then every function continuous in $E$ and holomorphic on $E^\circ$ admits uniform approximation by polynomials.
\begin{thm}
Let $K\subset\mathbb H$ be a compact axially symmetric set such that $K=\overline{K^\circ}$.
Then
\[\widetilde{\mathcal A}_L(K)=\overline{\mathbb H[q]}^{\,C(K,\mathbb H)}\]
if and only if $\mathbb C_I\setminus E$ is connected where, $E=K\cap \mathbb{C}_I$.
\end{thm}
\begin{proof}
Suppose that $\mathbb C_I\setminus E$ is connected. Let $f\in \widetilde{\mathcal A}_L(K)$ and $\varepsilon>0$. Then by the lemma \ref{splitting lemma}, there exist
$g_1,g_2\in C(E,\mathbb C_I)$, holomorphic on $K^\circ \cap \mathbb{C}_I$, such that
\[	f(z)=g_1(z)+g_2(z)J,\qquad z\in E.\]
Since $K^\circ$ is an axially symmetric slice domain, we have $K^\circ \cap \mathbb{C}_I=E^\circ$ and $E=\overline{K^\circ \cap \mathbb{C}_I}$.
Thus by Theorem \ref{keldysh}  for every $\delta>0$, polynomials $p_1,p_2\in\mathbb C_I[z]$ such that
\[\|g_1-p_1\|_E<\delta, \qquad \|g_2-p_2\|_E<\delta.\]
Define $P(q)=p_1(q)+p_2(q)J\in\mathbb H[q].$ Now,
\[|f(z)-P(z)|^2 = |g_1(z)-p_1(z)|^2 + |g_2(z)-p_2(z)|^2.\]
It follows that $\|f-P\|_E<\sqrt2\,\delta.$
Since $K= \overline{K^\circ}$ we have,
\[\|f-P\|_K\leq2\|f-P\|_E<2\sqrt2\,\delta.\]
Choose $0<\delta<\frac{\varepsilon}{2\sqrt2}$. Then,
\[\|f-P\|_K<\varepsilon.\]
Thus, $\widetilde{\mathcal A}_L(K)\subseteq\overline{\mathbb H[q]}^{\,C(K,\mathbb H)}.
$
Conversely, let $P_n\in\mathbb H[q]$ converge uniformly
on $K$ to $f$. Then $f\in C(K,\mathbb H)$ and $f$ is left slice regular on $K^\circ$. 
Hence $f\in \widetilde{\mathcal A}_L(K)$, and $\overline{\mathbb H[q]}^{\,C(K,\mathbb H)} \subseteq \widetilde{\mathcal A}_L(K).$ Thus,
\[\widetilde{\mathcal A}_L(K)=\overline{\mathbb H[q]}^{\,C(K,\mathbb H)}.\]

Suppose that $\widetilde{\mathcal A}_L(K)=\overline{\mathbb H[q]}^{\,C(K,\mathbb H)}$ and  $\mathbb C_I\setminus E$ is disconnected. Thus by Theorem \ref{keldysh}, there exists
$g\in C(E,\mathbb C_I), \ g|_{K^ \circ \cap \mathbb{C}_I}\ \text{ holomorphic},$ which cannot be uniformly approximated on $E$ by polynomials in $\mathbb C_I[z]$. Now since $K=\overline{K^\circ}$,  by  Theorem \ref{representation formula} there exist  $f\in \widetilde{\mathcal A}_L(K)$ such that $f|_E=g.$ Now by our assumption, there are polynomials
$P_n\in\mathbb H[q]$ such that
\[\|f-P_n\|_K\longrightarrow0.\]
Thus, $\|g-P_n\|_E\longrightarrow0.$ Now by Lemma \ref{splitting lemma} for $p_n$  on $\mathbb C_I$ such that
\[P_n(z)=p_{1,n}(z)+p_{2,n}(z)J,\qquad p_{1,n},p_{2,n}\in\mathbb C_I[z].\]
Now, 
\[\|g-p_{1,n}\|_E\leq\|g-P_n\|_E\longrightarrow0\] which is a contradiction. Hence our assumption is wrong and $\mathbb C_I\setminus E$ is connected.
\end{proof}
Here we consider all continuous quaternion valued functions on $K$ and slice regular on $K^\circ$ and we have to able to approximate all these functions by quaternionic polynomials if $\mathbb{C}\setminus E$ is connected. Also by Lemma \ref{approxsubset}, we have 
\[\overline{\mathbb H[q]}^{\,C(K,\mathbb H)} = \mathcal{R}_L(K) = \mathcal{A}_L(K).\]
In Complex setting for arbitrary compact sets these theorems are extended by Mergelyan's approximation Theorem.
For the quaternionic analogue of this we know that by Example \ref{noapprox} it doesn't holds in general in the case of all functions in $\mathcal C(K, \mathbb H)$ even when the compact set is axially symmetric and nowhere dense and its complement is an axially symmetric slice domain. 
The  complex Mergelyan theorem applies to every compact
$K\subset\mathbb C$ with connected complement. A general theorem
of this scope is not available for arbitrary quaternionic compact
sets. Gal and Sabadini in \cite{quatern approx} proved some type of Mergelyan approximation theorem in quaternionic for all functions in $\mathcal C(K, \mathbb H)$ with some geometric constraints. They proved polynomial approximation
for several classes of sets, using a quaternionic analogue of the
Riemann mapping theorem in restricted settings and explicit
approximation operators.

\begin{thm}\cite{quatern approx}
\label{qball}
Let $R>0$ and let $	K=\overline {\mathbb B}(0,R)=\{q\in\mathbb H:|q|\leq R\}.$	If $f\in \widetilde{A}_L(K)$, then for every $\varepsilon>0$ there exists a 
quaternionic polynomial $P$ such that $\|f-P\|_K<\varepsilon.$
\end{thm}

\begin{thm}\cite{quatern approx}
\label{q-cassini}
Let $q_0=x_0+Iy_0$, where $x_0,y_0\in\mathbb R$, $y_0>0$,
and let $R>0$. Consider a Cassini cell of the form
$K=\{q\in\mathbb H:|q^2-2x_0q+x_0^2+y_0^2|\leq R^2\}.$
If $f\in \widetilde{A}_L(K)$  there exists a  quaternionic
polynomial $P$ such that $\|f-P\|_K<\varepsilon.$	
\end{thm}

\begin{thm}\cite{quatern approx}
\label{q-jordan-slices}
Let $T$ be a bounded set in the class $\mathcal R(\mathbb H)$
,class of axially symmetric open sets $U$ in $\mathbb{H}$ such that $U\cap \mathbb{C}_I$ is simply connected and $T\cap\mathbb C_I$ is a Jordan region. Then,
$f \in \widetilde{A}_L(T)$  can be uniformly approximated on $T$ by  quaternionic polynomials.
\end{thm}
We hereby prove a version of this for  continuous slice functions on an arbitrary compact axially symmetric set .
\begin{thm}
Let $K\subset\mathbb H$ be a compact axially symmetric set.
Then $\mathcal A_L(K) 	= 	\overline{\mathbb H[q]}^{\,\|\cdot\|_K}$ if and only if $\mathbb C_I\setminus E$ is connected.
\end{thm}

\begin{proof}
Assume that $\mathbb C_I\setminus E$ is connected. It is enough to prove that $\mathcal A_L(K)\subseteq
\overline{\mathbb H[q]}^{\,\|\cdot\|_K}.$
Let $f\in\mathcal A_L(K)$. Then by Lemma \ref{splitting lemma} for every $z\in E$, we have
\[	f(z)=g_1(z)+g_2(z)J,\]
where $g_1,g_2:E\to\mathbb C_I$ are uniquely determined. Then  since $f$ is a continuous slice function on $K$, $g_1$ and $g_2$ are continuous on $E$.	Also, since $K$ is axially symmetric, $	E^\circ=K^\circ\cap\mathbb C_I$. Thus $g_1$ and $g_2$ is slice regular on $E^\circ$. Since $\mathbb C_I\setminus E$ is connected, by \ref{Mergelyan}, for every $\varepsilon>0$ there exist
complex polynomials $p_1$ and $p_2$ such that
\[	\sup_{z\in E}|g_1(z)-p_1(z)|<\varepsilon\ \text{and } \sup_{z\in E}|g_2(z)-p_2(z)|<\varepsilon.
\]
Now, define $P(q)=p_1(q)+p_2(q)J.$	Then $P\in\mathbb H[q]$ and for every $z\in E$,
\[	\begin{aligned}
	|f(z)-P(z)|
	&=
	|(g_1(z)-p_1(z))+(g_2(z)-p_2(z))J|\\
	&\leq |g_1(z)-p_1(z)|+|g_2(z)-p_2(z)|<2\varepsilon.
\end{aligned}
\]
Hence $	\|f-P\|_E<2\varepsilon.$	
Since, $f-P$ is a continuous left slice function
on $K$, by the Theorem \ref{representation formula}, for every $q\in K$,
\[\|f-P\|_K\leq 2\|f-P\|_E<4\varepsilon.\]
Since $\varepsilon>0$ is arbitrary, we have $	f\in\overline{\mathbb H[q]}^{\,\|\cdot\|_K}.$
Thus
\[\mathcal A_L(K)\subseteq\overline{\mathbb H[q]}^{\,\|\cdot\|_K}.\]

Suppose that $\mathbb C_I\setminus E$ is disconnected. Since, $E$ is compact, there exists a bounded connected component $U$ of $\mathbb C_I\setminus E$. Choose $a\in U$. Let, $a=\alpha+\beta I,
\ \alpha,\beta\in\mathbb R,$
and define
\[	Q(q)=q^2-2\alpha q+|a|^2.\]
which is an intrinsic polynomial, and its zero set is the
quaternionic conjugacy sphere determined by $a$; if $a$ is
real, this sphere reduces to the point $a$. Since $a\notin E$ and $K$ is axially symmetric, the zero set of $Q$ does not intersect $K$. Indeed, if $q\in K$ with $Q(q)=0$, then $q$ would be conjugate to $a$. Then by	axial symmetry  $a\in K\cap\mathbb C_I=E$, a contradiction. Hence $Q(q)\neq0,\ q\in K.$	
Since $K$ is compact, $Q$ is bounded away from zero on $K$.
Thus $f(q)=Q(q)^{-1}$ is  continuous on $K$ and slice regular on a
neighborhood of $K$ since $Q$ is intrinsic. In particular, $f\in\mathcal A_L(K).$	

Suppose that $\mathcal A_L(K) =\overline{\mathbb H[q]}^{\,\|\cdot\|_K}.$
Then there exists a sequence $(P_n)$ in $\mathbb H[q]$
such that
\[\|f-P_n\|_K\longrightarrow0.\]
Now, for $z\in\mathbb C_I$, $P_n(z)=p_n(z)+r_n(z)J,$
where $p_n$ and $r_n$ are polynomials with coefficients
in $\mathbb C_I$. 
Since $f(z)=Q(z)^{-1}\in\mathbb C_I$ for $z\in E$, 
\[	\left|\frac{1}{Q(z)}-p_n(z)\right|\leq |f(z)-P_n(z)|.\]
Hence,
\[	\sup_{z\in E}\left|\frac{1}{Q(z)}-p_n(z)\right|\longrightarrow0.
\]

Define the complex polynomial $h_n(z)=1-Q(z)p_n(z).$ Since, $Q$ is bounded on $E$, we have
\[
\begin{aligned}
	\sup_{z\in E}|h_n(z)|
	&=
	\sup_{z\in E}
	\left|Q(z)\left(\frac{1}{Q(z)}-p_n(z)\right)\right|\\
	&\leq
	\sup_{z\in E}|Q(z)|
	\sup_{z\in E}
	\left|\frac{1}{Q(z)}-p_n(z)\right|
	\longrightarrow0.
\end{aligned}
\]
On the other hand, $Q(a)=0$, and hence $h_n(a)=1$ for every $n$.	Since $U$ is a bounded connected component of
$\mathbb C_I\setminus E$, its boundary satisfies $\partial U \subseteq E.$ Thus by the  maximum modulus principle applied to $h_n$ on $U$ gives,
\[
\begin{aligned}
	1
	&=|h_n(a)|\\
	&\leq\sup_{z\in\partial U}|h_n(z)|\\
	&\leq\sup_{z\in E}|h_n(z)|\rightarrow 0
\end{aligned}
\]
as $n\to\infty$ which is a contradiction.	
Thus $\mathbb C_I\setminus E$ must be connected. 
\end{proof}
By this theorem we can say that if $\mathbb C_I\setminus E$ is connected,
$\mathcal \mathcal A_L(K) = \mathcal R_L(K) = 	\overline{\mathbb H[q]}^{\,\|\cdot\|_K}.$
Mergelyan in \cite{Mergelyan} proved that if the complement has finitely many components every continuous functions $K$ which is analytic on $K^\circ$ can be uniformly approximated by rational functions. For quaternionic version of this we have,

\begin{thm}\label{rlclal}
	Let $K\subset\mathbb H$ be compact and axially symmetric.
	Fix $I\in\mathbb S$ and put $E=K\cap\mathbb C_I.$
	Suppose that $\mathbb C_I\setminus E$ has finitely many
	connected components. Then
	\[
	\mathcal R_L(K)=\mathcal A_L(K).
	\]
\end{thm}

Before proving this we are going to prove this lemma which identifies the uniform rational approximation problem on $K$ with the corresponding complex rational approximation problem on $E = K \cap \mathbb C_I.$  Let $C_{\mathbb{C}}(E)$, $R_{\mathbb{C}}(E)$ and $A_{\mathbb{C}}(E)$ be the set of all continuous complex valued functions on $E$, rationally approximable functions on $E$ and continuous functions on $E$ which are analytic on $E^0$ respectively.
\begin{lem}\label{10}
	Let $K \subset \mathbb{H}$ be a compact axially symmetric set and let $I \in \mathbb{S}$. Let $E = K \cap \mathbb{C}_I$. Then,
	\begin{enumerate}
		\item $ \mathcal{R}_L(K)^{\|\cdot\|_K} = \mathcal{C}_L(K)$  if and only if  $\mathcal R_{\mathbb C}(E)^{\|\cdot\|_{E}} = \mathcal C_{\mathbb C}(E),$
		\item $\mathcal R_L(K) = A_L(K)$ if and only if $\mathcal R_{\mathbb C}(E) = \mathcal A_{\mathbb C}(E)$.
		
	\end{enumerate}
	
\end{lem}
\begin{proof}
	We will give the proof of (1). The proof of (2) follows as in (1).
	Suppose that $\mathcal{R}_L(K) = \mathcal{C}_L(K)$. Let $g \in \mathcal C_{\mathbb C}(E)$ and let $\epsilon > 0$. Then by the representation formula, $g$ admits a continuous left slice extension, $ \tilde{g} \in \mathcal{C}_L(K)$  such that  $\tilde{g}|_{E} = g$. Now since $\mathcal{R}_L(K) = \mathcal{C}_L(K)$, there exists $r \in \mathcal{R}_L(K)$ such that $\| \tilde{g} - r \|_K < \epsilon$. If we restrict to $E$, then $\| g - r|_{E} \|_{E} < \epsilon$. 	Then we can find a rational complex valued function $r_I$ on $\mathbb{C}_I$ such that $ \| g - r_I \|_{E} \leq \| g - r|_{E} \|_{E} < \epsilon .$	That is, $g$ has been approximated by a complex rational function. Since $g \in \mathcal C_{\mathbb C}(E)$ was arbitrary,
	\[\mathcal R_{\mathbb C}(E) = \mathcal C_{\mathbb C}(E). \]
	
	Now, conversely, suppose that $\mathcal R_{\mathbb C}(E)=\mathcal C_{\mathbb C}(E).$ 	Let $f\in \mathcal C_L(K)$ and let $\varepsilon>0$. Then by the  Lemma \ref{splitting lemma}, the restriction $f_I=f|_{E}$ can be written as
	$f_I=G_1+G_2J,\  G_1,G_2\in \mathcal C_{\mathbb C}(E),$ where $J\in\mathbb S$ is orthogonal to $I$. Since, $\mathcal R_{\mathbb C}(E)=\mathcal C_{\mathbb C}(E),$	there exist $r_1,r_2\in\mathcal R_{\mathbb C}(E)$ such that 
	\[\|G_1-r_1\|_{E}<\frac{\varepsilon}{4},\			\|G_2-r_2\|_{E}<\frac{\varepsilon}{4}.\]
	Let, $r_I=r_1+r_2J.$ Thus, by the characterization of left rational functions, there exist some  $r\in\mathcal R_L(K)$ such that have $r_I = r|E$. Also, 
	\begin{align*}
		\|f_I-r_I\|_{E}
		&=
		\|(G_1-r_1)+(G_2-r_2)J\|_{E}\\
		&\leq
		\|G_1-r_1\|_{E}
		+
		\|G_2-r_2\|_{E} <
		\frac{\varepsilon}{4}+\frac{\varepsilon}{4} = 			\frac{\varepsilon}{2}.
	\end{align*}
	Let, $h=f-r.$	Then $h$ is a left slice function on the axially symmetric set $K$. Let $q=x+yL\in K,\ L\in\mathbb S,$
	and set $z=x+yI\in E.$ Then by the representation formula,
	$h(q) = \frac12(1-LI)h(z) + \frac12(1+LI)h(\overline z).$
	Consequently,
	\begin{align*}
		|h(q)|
		&\leq
		\frac12|1-LI|\,|h(z)| + \frac12|1+LI|\,|h(\overline z)|\\
		&\leq	|h(z)|+|h(\overline z)| \leq 2\|h\|_{E}.
	\end{align*}
	Now by taking the supremum over $q\in K$, we have $		\|f-r\|_K = \|h\|_K \leq 2\|h\|_{E} = 			2\|f_I-r_I\|_{E}	<\varepsilon.$
	Thus every $f\in \mathcal C_L(K)$ can be uniformly approximated on $K$ by elements of $\mathcal R_L(K)$. Hence, $\mathcal R_L(K)=\mathcal C_L(K).$
\end{proof}  
\begin{proof}[Proof of \ref{rlclal}]
	Since $\mathbb C_I\setminus E$ has finitely many connected
	components, the Mergelyan rational approximation theorem
	we have $\mathcal R_{\mathbb C}(E)=\mathcal A_{\mathbb C}(E).$ Then by Theorem \ref{10},
	$\mathcal R_L(K)=\mathcal A_L(K).$
\end{proof}
If the complement has infinitely many components, Mergelyan's swiss cheese example shows that it is not in general the approximation holds.
Using these examples and Lemma \ref{10} we can provide counter examples in quaternionic.
\begin{exa}
	By Swiss cheese examples, there exists a compact planar Swiss-cheese set
	\[	E_0=\overline{\mathbb D}\setminus\bigcup_{n=1}^{\infty}D_n,
	\]
	where $\mathbb D$ is the open unit disk and the disks $D_n$ are
	chosen so that $E_0$ has empty interior, for which
	$R(E_0)\neq C(E_0).$ We  may assume that
	$E_0$ lies entirely in the open upper half-plane.
	
	Define $E=E_0\cup\overline{E_0},$ and $		\overline{E_0}:=\{\overline z:z\in E_0\}.$	Then $E$ is compact, has empty interior, and is invariant under
	complex conjugation. Fix $I\in\mathbb S$ and identify $\mathbb C$ with $\mathbb C_I$. Define the circularization of $E$ by
	\[K=\left\{x+yJ\in\mathbb H:x+yI\in E,\ J\in\mathbb S\right\}.
	\]
	Since, $E$ is compact and invariant under complex conjugation, $K$ is compact and axially symmetric. Also,
	$K\cap\mathbb C_I=E$ and we can prove that $K$ has same properties as $E$. Then by Lemma \ref{10} we have
	$\mathcal R_L(K)=\mathcal A_L(K)$ if and only if $\mathcal R_{\mathbb C}(K\cap\mathbb C_I)=\mathcal A_{\mathbb C}(K\cap\mathbb C_I).$ Since $K\cap\mathbb C_I=E$, and $\mathcal R_{\mathbb C}(E)\neq \mathcal A_{\mathbb C}(E)$, we have, $\mathcal R_L(K)\neq \mathcal A_L(K).$
\end{exa}
\section{Quaternionic Version of Vitushkin Approximation}
We are now in position to extend complex Vitushkin approximation to quaternionic setting. The complex Vitushkin theorem is based on a localization procedure in the complex plane. A compact set is decomposed into small pieces, and the approximation problem on the whole set is related to the corresponding problems on these localized pieces. In the complex setting, Euclidean discs provide the natural localization sets. In the quaternionic setting, however, the direct use of Euclidean balls
requires some care because of the non-commutative and axially symmetric
geometry of slice regular functions.

Since, for  axially symmetric sets  by theorem \ref{complexslice quaternionic} we have  $\gamma_{\mathbb H}(S)=\gamma_{\mathbb C}(E)$ for an axially symmetric set $S$ and $E=S\cap\mathbb C_I$,   we exploit the reduction of the quaternionic problem to a complex slice throughout this section.  Thus, in order to transfer the localization arguments from the complex Vitushkin theory to the quaternionic setting, the localization sets in a slice should be chosen in such a way that their quaternionic counterparts remain axially symmetric. 

A ball $\mathbb B_{\mathbb H}(q,r) =\{p\in\mathbb H:|p-q|<r\}$ with a non-real center $q\notin\mathbb R$ is not axially symmetric, although $q$ belongs to the ball, the whole $2$-sphere $[q]=\{x+Jy:J\in\mathbb S\},\  q=x+Iy,$
need not be contained in the ball. Consequently, such a ball is not
compatible with the axial symmetry required for the slice reduction.
Thus, one could restrict the localization to balls centered at real points, once these balls are axially symmetric. However, this would be too restrictive for the present purpose. In a fixed slice $\mathbb C_I$, the complex Vitushkin decomposition uses discs centered at arbitrary points of the slice, whereas restricting the centers to the real axis
is not acceptable for  this direct correspondence. In particular, a small
neighborhood of a non-real point $p\in\mathbb C_I$ cannot in general
be represented by a sufficiently small ball centered on the real
axis without introducing a sufficient enlargement of the localized
set. Thus we are extending this to more general compact sets, Cassini cells. For $q = x + I y \in \mathbb{H}$, with $x, y \in \mathbb{R}$, $y > 0$, $I \in \mathbb{S}$, and $r > 0$, let us denote:
\[ B(x + y \mathbb{S}; r) := \left\{ q_0 \in \mathbb{H} \mid \| (q_0 - x)^2 + y^2 \| < r^2 \right\}\ .\]
A set of the form $B(x + y \mathbb{S}; r)$ will be called a Cassini cell \cite{quatern approx}. For a fixed $J \in \mathbb{S}$, let $\mathbb{C}_J = \{ a + Jb \mid a, b \in \mathbb{R} \}$, $p_J \in \mathbb C_J$. Then for every $J \in \mathbb{S}$, the Cassini cell satisfies the following properties:
\begin{enumerate}
	\item $B(x + y \mathbb{S}; r) \cap \mathbb{C}_J = B_J(p_J, r) = \left\{ z \in \mathbb{C}_J : |z - p_J||z - \overline{p}_J| < r^2 \right\}$.
	
	\noindent For, let $z = a + Jb \in \mathbb{C}_J$; $J \in \mathbb{S}$. If $z\in B(x + y \mathbb{S}; r)$, then by definition $|(z - x)^2 + y^2 |<r^2$.   Since $z, p_J, \overline{p}_J \in \mathbb{C}_J$ and the slice $\mathbb{C}_J$ is commutative, $$(z - p_J)(z - \overline{p}_J) =  (z - x)^2 + y^2 .$$ That is, $\|(z - x)^2 + y^2 \| = \| z - p_J \| \| z - \overline{p}_J \|.$
	Hence, $ z \in B(p_J, r)$ if and only if $|z - p_J||z - \overline{p}_J| < r^2$. 
	\item $B(x + y \mathbb{S}; r)$ is axially symmetric.
	\item By the condition $\| (q - x)^2 + y^2\| < r^2$ it depends only on $x$ and $y^2$, and thus it does not depend on the particular choice of $I \in \mathbb{S}$. This proves the independence of the choice of the slice.
	\item If  $p$ is non-real, the Cassini cell in $\mathbb{C}_J$ contains 2 components, 	$\mathcal{B}_J^+(p, r),$ a component containing $p$ and $\mathcal{B}_J^-(p, r),$  a component containing  $\overline{p}.$
\end{enumerate}
For more details on Cassini cell and Cassini -pseudo metric refer, \cite{quatern approx,cassini}.

We know that complex Vitushkin theorem is formulated using disc and in quaternionic we are choosing Cassini cells. On a fixed slice $\mathbb C_I$ and for $p = x + Iy \in \mathbb{C}_I \setminus \mathbb{R}$ Cassini cells are defined by $B_I(p, r) = \left\{ z \in \mathbb{C}_I : |z - p||z - \overline{p}| < r^2 \right\}.$ The following Lemma makes a comparison between disc in $\mathbb{C}$,  and $B_I(p, r)$.
\begin{lem}\label{12}
	Let $p = x + Iy \in \mathbb{C}_I \setminus \mathbb{R}$, $y = \text{Im}(p) > 0$, and let $0 < r < y$, and
	\[ B_I(p, r) = \left\{ z \in \mathbb{C}_I : |z - p||z - \overline{p}| < r^2 \right\}. \]
	Let $B_I^+(p, r)$ denotes the component containing $p$. Put	$ \rho_p(r) = y - \sqrt{y^2 - r^2}$.
	Then, for all sufficiently small $r$,
	\[ \Delta_I\left(p, \frac{\rho_p(r)}{2}\right) \subset B_I^+(p, r) \subset \Delta_I(p, \rho_p(r)),\]
	where $\Delta_I$ is the Euclidean disc in the complex slice $\mathbb C_I$.
	Also, $$\displaystyle \rho_p(r) = \frac{r^2}{y + \sqrt{y^2 - r^2}} \sim \frac{r^2}{2y}$$ for sufficiently small $r>0$.
\end{lem}

\begin{proof}
	Let $z = p + w$ where $w\in \mathbb{C}$. Then $|z - p| = |w|$ and $|z - \overline{p}| = |2Iy + w|$. Hence, $$2y - |w| \leq |2Iy + w| \leq 2y + |w|.$$ Let $s = |z - p|$. If $z \in \mathcal{B}_I^+(p, r) \subset  B_I(p, r)$, then	$s(2y - s) < r^2$.
	The smallest positive root of $s(2y - s) = r^2$ is $s = y - \sqrt{y^2 - r^2} = \rho_p(r)$. That is $s< \rho_p(r)$ and thus, $|z-p|< \rho_p(r)$. Hence, $$\mathcal{B}_I^+(p, r) \subset \Delta_I(p, \rho_p(r)).$$ 
	
	Suppose that $|z - p| < \dfrac{\rho_p(r)}{2}$. Then, $|z-\bar{p}|\leq 2y+q$ and thus,
	$$ |z - p||z - \overline{p}| \leq s(2y + s) <\frac{\rho_p(r)}{2} \left( 2y + \frac{\rho_p(r)}{2} \right).$$
	Since,
	$$ \rho_p(r) = y - \sqrt{y^2 - r^2} = \dfrac{r^2}{y + \sqrt{y^2 - r^2}}$$
	$\rho_p(r) \sim \dfrac{r^2}{2y}$, for sufficiently small $r$. Thus,, $$ \frac{\rho_p(r)}{2} \left( 2y + \frac{\rho_p(r)}{2} \right) < r^2. $$	Hence,  $\Delta_I\left(p, \frac{\rho_p(r)}{2}\right) \subset \mathcal{B}_I^+(p, r) .$ 
\end{proof}
The preceding lemma gives a geometric comparison between the component of a Cassini cell containing $p$ and an ordinary disc centered at $p$. For the Vitushkin argument, however geometric inclusion is not enough. We need the corresponding comparison for analytic capacity. 

The next Lemma shows that the normalized analytic capacity is unchanged when ordinary discs are replaced by the appropriate Cassini components.
\begin{lem}\label{disccell}
	Let $K \subset \mathbb{H}$ be compact axially symmetric set and let $E = K \cap \mathbb{C}_I$. Let $p = x + Iy \in E \setminus \mathbb{R}$, $y = \text{Im}(p) > 0$. Then,
	\[ \limsup\limits_{\rho \rightarrow 0} \frac{\gamma_{\mathbb{C}_I}(\Delta_I(p, s) \setminus E)}{s} > 0 \text{ if and only if } \limsup\limits_{r \rightarrow 0} \frac{\gamma_{\mathbb{C}_I}(B_I(p, r) \setminus E)}{\rho_p(r)} > 0. \]
\end{lem}
\begin{proof}
	By Lemma \ref{12} we have,
	\[ \Delta_I\left(p, \dfrac{\rho_p(r)}{2}\right) \subset B_I^+(p, r) \subset \Delta_I(p, \rho_p(r)) \]
	The monotonicity of analytic capacity gives,
	\[ \gamma_{\mathbb{C}_I}\left(\Delta_I\left(p, \frac{\rho_p(r)}{2}\right) \setminus E\right) \leq \gamma_{\mathbb{C}_I}(\mathcal{B}_I^+(p, r) \setminus E) \]
	and
	\[ \gamma_{\mathbb{C}_I}(\mathcal{B}_I^+(p, r) \setminus E) \leq \gamma_{\mathbb{C}_I}(\Delta_I(p, \rho_p(r)) \setminus E). \]
	Let $s=\rho_p(r)$. Since $\rho_p(r) \rightarrow 0$ as $r \rightarrow 0$,  $r \mapsto \rho_p(r)$ is continuous and strictly increasing for $0<r<y$,  every sufficiently small positive $s$ corresponds to exactly one positive $r$. Also $s\to 0$ if and only if $r\to 0$.
	
	\noindent Thus we have, \[\limsup\limits_{s \rightarrow 0} \dfrac{\gamma_{\mathbb{C}_I}(\Delta_I(p, s) \setminus E)}{s} > 0\] is equivalent to 
	\[\limsup\limits_{r \rightarrow 0} \dfrac{\gamma_{\mathbb{C}_I}(\mathcal{B}_I^+(p, r) \setminus E)}{\rho_p(r)} > 0.\]
	Since $E$ is invariant under complex conjugation, we have
	\[ \mathcal{B}_I^-(p, r) \setminus E = \overline{\mathcal{B}_I^+(p, r) \setminus E} .\]
	Hence, $\gamma_{\mathbb{C}_I}(\mathcal{B}_I^-(p, r) \setminus E) = \gamma_{\mathbb{C}_I}(\mathcal{B}_I^+(p, r) \setminus E)$.
	
	By semi-additivity of analytical capacity, there exists a universal constant $C > 0$ such that,
	$$ \gamma_{\mathbb{C}_I}(B_I(p, r) \setminus E) \leq C \left[ \gamma_{\mathbb{C}_I}(B_I^+(p, r) \setminus E) + \gamma_{\mathbb{C}_I}(B_I^-(p, r) \setminus E) \right] .$$ 
	Thus, $ \gamma_{\mathbb{C}_I}(B_I(p, r) \setminus E) \leq 2C \gamma_{\mathbb{C}_I}(B_I^+(p, r) \setminus E). $
	
	Conversely by monotonicity, $\gamma_{\mathbb{C}_I}(B_I^+(p, r) \setminus E) \leq \gamma_{\mathbb{C}_I}(B_I(p, r) \setminus E)$.
	Thus, the positivity of the upper limit is unchanged when $B_I^+(p, r)$ is replaced by the Cassini cell $B_I(p, r)$. Hence the proof.
\end{proof}
Vitushkin \cite{Vitushkin1, Vitushkin2, Vitushkin3} reduce the problem of rational approximation of compact set $E$ with interior to a problem of analytic capacity.
\begin{thm}\cite{rational app, analytic capacity}\label{gammavit}
	The following assertions are equivalent
	\begin{enumerate}
		\item $R_{\mathbb{C}}(E) = C_{\mathbb{C}}(E)$.
		\item For every open set $D\subset C$, $\gamma(D \setminus E) = \gamma(D)$.
		\item For $dxdy$-almost all $z \in E$, 
		$		\limsup_{r \to 0} \frac{\gamma(\Delta(z; r) \setminus E)}{r} > 0.$
	\end{enumerate}
\end{thm}
We are proving the quaternionic analogue of  Vitushkin approximation theorem using $\gamma_{\mathbb{H}}$ capacity.
\begin{thm}
	Let $K \subset \mathbb{H}$ be a compact axially symmetric set. For $p = x + Iy \in K \setminus \mathbb{R}$, put $\rho_p(r) = |y| - \sqrt{y^2 - r^2}$. Then the following are equivalent for $K$:
	\begin{enumerate}[(i)]
		\item $\mathcal{R}_L(K) = \mathcal{C}_L(K)$.
		\item $\gamma_{\mathbb{H}}(U \setminus K) = \gamma_{\mathbb{H}}(U)$, for all axially symmetric open set $U$.
		\item $\limsup\limits_{r \rightarrow 0} \dfrac{ \gamma_{\mathbb{H}}(B(x + y \mathbb{S}; r)\setminus K)}{\rho_p(r)} > 0$ for $m_4$ (Lebesgue) a.e. non-real point $p \in K$.
	\end{enumerate}
	
\end{thm}
\begin{proof}
	Suppose (i) holds. Let $E = K \cap \mathbb{C}_I$. Then by Lemma \ref{10}, $\mathcal R_{\mathbb C}(E) = \mathcal C_{\mathbb C}(E)$. Now by Theorem \ref{gammavit}, $$\gamma_{\mathbb{C}_I}(D \setminus E) = \gamma_{\mathbb{C}_I}(D),$$ for every open set $D \subset \mathbb{C}_I$.	Let $U \subset \mathbb{H}$ be any axially symmetric open set and then $U_I = U \cap \mathbb{C}_I$ is an open subset of $\mathbb{C}_I$. Hence, we can take $D = U_I$. Thus,
	$ \gamma_{\mathbb{C}_I}(U_I \setminus E) = \gamma_{\mathbb{C}_I}(U_I). $ 
	Since $U$ and $K$ are axially symmetric, $(U \setminus K) \cap \mathbb{C}_I = U_I \setminus E$. Also we have $\gamma_{\mathbb{H}}(U) = \gamma_{\mathbb{C}_I}(U_I)$. Hence,
	$ \gamma_{\mathbb{H}}(U \setminus K) = \gamma_{\mathbb{C}_I}(U_I \setminus E)$ and thus, $ \gamma_{\mathbb{H}}(U \setminus K) = \gamma_{\mathbb{C}_I}(U_I) = \gamma_{\mathbb{H}}(U) .$
	Hence (ii) holds.
	
	Now, suppose that (ii) holds. Let $E = K \cap \mathbb{C}_I$ and $p = x + Iy \in E \setminus \mathbb{R}$; $y > 0$. Since $K$ is axially symmetric, $E$ is symmetric with respect to complex conjugation. Consider an open set $D \subset \mathbb{C}_I$ which is symmetric under conjugation, that is, $D = \overline{D}$. Its axially symmetric completion $U $ is an axially symmetric open subset of $\mathbb{H}$ and $U \cap \mathbb{C}_I = D$. Thus by (ii), $\gamma_{\mathbb{H}}(U \setminus K) = \gamma_{\mathbb{H}}(U)$ and by Theorem \ref{complexslice quaternionic}, $\gamma_{\mathbb{H}}(U) = \gamma_{\mathbb{C}_I}(D)$ and $\gamma_{\mathbb{H}}(U \setminus K) = \gamma_{\mathbb{C}_I}(D \setminus E)$. 
	
	Thus, $\gamma_{\mathbb{C}_I}(U \setminus E) = \gamma_{\mathbb{C}_I}(U)$. Now by  Theorem \ref{gammavit}, for almost every $p \in E \setminus \mathbb{R}$, we have
	$$\limsup_{s \rightarrow 0} \frac{\gamma_{\mathbb{C}_I}(\Delta_I(p, s) \setminus E)}{s} > 0. $$
	Now, for $p = x + Iy$, $y > 0$ we have $\rho_p(r) = y - \sqrt{y^2 - r^2}$. Thus by Lemma \ref{disccell}, we have,
	$$ \limsup_{r \rightarrow 0} \frac{\gamma_{\mathbb{H}}(B_I(p, r) \setminus K)}{\rho_p(r)} > 0 .$$ 
	since, $B(x + y \mathbb{S}; r) \cap \mathbb{C}_I = B_I(p_I, r)$, and by Theorem \ref{complexslice quaternionic},
	$\limsup\limits_{r \rightarrow 0} \dfrac{\gamma_{\mathbb{H}}(B(x + y \mathbb{S}; r) \setminus K)}{\rho_p(r)} > 0$ for $m_4$ (Lebesgue) a.e. non-real point $p \in K$. 
	
	Suppose (iii) holds. Let $E = K \cap \mathbb{C}_I$. Then by Lemma \ref{disccell},
	$$ \limsup_{r \rightarrow 0} \dfrac{\gamma_{\mathbb{C}_I}(B_I(p, r) \setminus K)}{\rho_p(r)} > 0$$ if and only if $$ \limsup_{s \rightarrow 0} \dfrac{\gamma_{\mathbb{C}_I}(\Delta_I(p, s) \setminus E)}{s} > 0 .$$ 
	Hence by Theorem \ref{gammavit}, $\mathcal R_{\mathbb C}(E) = \mathcal C_{\mathbb C}(E) $. Thus, by Lemma \ref{10},
	$ \mathcal{R}_L(K) = \mathcal{C}_L(K) $.
\end{proof}
The Vitushkin theorem we proved provides characterization of rational approximation in terms of the local behaviour of slice regular capacity.  For the corresponding characterization involving functions which are continuous on the whole compact set and slice regular in its interior, the appropriate notion of capacity is the continuous slice regular capacity. In contrast to the slice regular capacity of the compact set $K$ which measures bounded slice regular functions outside $K$, the capacity $\alpha_{\mathbb{H}}$ is related to also the class of functions continuous across the compact set.  Thus for studying rational approximation to the slice regular capacity as in complex, we here introduced continuous slice regular capacity. In the following part we are extending that version of Vitushkin theorem for a compact set $E$ which is given as follows:
\begin{thm}\label{vitalpha}
	The following are equivalent
	\begin{enumerate}
		\item $R_{\mathbb{C}}(E) = A_{\mathbb C}(E)$.
		\item  $\alpha(D \setminus E) = \alpha(D \setminus E^\circ)$ for every bounded open set $D\subset \mathbb C$.
		\item For each complex $z$, $r > 1$, and  $\delta > 0$,
		$\alpha(\Delta(z; \delta) \setminus E^\circ \le \alpha(\Delta(z; r\delta) \setminus E).$
		\item There exist $r \geq 1$ and $c > 0$ such that for all  $z\in \mathbb C$ and all $\delta > 0$,
		$	\alpha(\Delta(z; \delta) \setminus E^\circ) \leq c \alpha(\Delta(z; r\delta) \setminus E).$
		\item For each $z \in \partial E$, there exists $r \ge 1$ such that
		$\limsup_{\delta \to 0} \frac{\alpha(\Delta(z; \delta) \setminus E^\circ)}{\alpha(\Delta(z; r\delta) \setminus E)} < \infty.$
	\end{enumerate}
\end{thm}
Now, we have the Quaternionic Vitushkin theorem which proved the approximation problem which allow us to approximate $f\in \mathcal A_L(K)$ by functions which extend analytically beyond the interior of $K$ using theorem  \ref{complexslice quaternionic} which help us to reduce the problem to slice.
\begin{thm}
	Let $K \subset \mathbb{H}$ be a compact axially symmetric set. Then the following are equivalent:
	\begin{enumerate}[(i)]
		\item $\mathcal R_L(K) = \mathcal A_L(K)$.
		\item For every bounded axially symmetric open set $U \subset \mathbb{H}$,
		\[ \alpha_{\mathbb{H}}(U \setminus K) = \alpha_{\mathbb{H}}(U \setminus K^\circ). \]
		\item for every $ q = x + y I \in \mathbb{H}$, every $r \geq 1$, and every $\delta > 0$,
		\[ \alpha_{\mathbb{H}}(B(x + y \mathbb{S}; \delta) \setminus K^\circ \leq  \alpha_{\mathbb{H}}(B(x + y \mathbb{S}; r\delta) \setminus K) .\]
		\item there exists $ \gamma \geq 1, C > 0$ such that for every $ q = x + y I \in \mathbb{H}$ and all $\delta > 0$,
		\[ \alpha_{\mathbb{H}}(B(x + y \mathbb{S}; \delta) \setminus K^\circ \leq C \alpha_{\mathbb{H}}(B(x + y \mathbb{S}; r\delta) \setminus K) .\]
		\item for every $ q_0 \in \partial K$, there exists $r > 1$ such that
		\[ \limsup_{\delta \to 0} \dfrac{\alpha_{\mathbb{H}}(B(x + y \mathbb{S}; \delta) \setminus K^\circ)}{\alpha_{\mathbb{H}}(B(x + y \mathbb{S}; r\delta) \setminus K)} < \infty .\]
	\end{enumerate}
	
\end{thm}

\begin{proof}
	Suppose that (i) holds. Let $E = K \cap \mathbb{C}_I$. 
	By Theorem \ref{10}, we have $\mathcal R_L(K) = \mathcal A_L(K)$ if and only if $\mathcal R_{\mathbb{C}}(E) = \mathcal A_{\mathbb{C}}(E)$. Let $U \subset \mathbb{H}$ be a bounded axially symmetric open set and put $U_I = U \cap \mathbb{C}_I$. 
	Since $U$ is axially symmetric, $U_I$ is an open subset of $\mathbb{C}_I$. Then by Theorem \ref{vitalpha},
	$$\alpha_{\mathbb{C}}(U_I \setminus E) = \alpha_{\mathbb{C}}(U_I \setminus \operatorname(int)_{\mathbb{C}_I} E) .$$	Since $U \setminus K$ is axially symmetric, by Theorem \ref{complexslice quaternionic}, 
	\begin{equation*}
		\alpha_{\mathbb{H}}(U \setminus K)  = \alpha_{\mathbb{C}}(U_I \setminus E) 
	\end{equation*}
	Since $K$ is axially symmetric, $$	(\operatorname{int}_{\mathbb{H}} K)_I = \operatorname{int}_{\mathbb{C}_I}(E), \ 
	\alpha_{\mathbb{H}}(U \setminus \operatorname{int}_{\mathbb{H}} K) = \alpha_{\mathbb{C}}(U_I \setminus \operatorname{int}_{\mathbb{C}_I} E).$$
	Hence from above equations we have, $ \alpha_{\mathbb{H}}(U \setminus K) = \alpha_{\mathbb{H}}(U \setminus K^\circ) .$
	
	Suppose that (ii) holds. Fix $q = x + y I \in \mathbb{H}$, $r \geq 1$, $\delta > 0$, and set $U = B(x + y \mathbb{S}; r\delta)$. Then $U$ is a bounded axially symmetric open set. Hence by (ii), 
	\begin{equation}\label{4}
		\alpha_{\mathbb{H}}(B(x + y \mathbb{S}; r\delta) \setminus K) = \alpha_{\mathbb{H}}(B(x + y \mathbb{S}; r\delta) \setminus K^\circ).
	\end{equation}
	Since $r \geq 1$, $B(x + y \mathbb{S}; \delta) \subseteq B(x_0 + y_0 \mathbb{S}, r\delta)$. Thus, $$B(x + y \mathbb{S}; \delta) \setminus K^\circ \subseteq B(x + y \mathbb{S}; r\delta) \setminus K^\circ.$$ Now by the monotonicity of $\alpha_{\mathbb{H}}$,
	\begin{equation}\label{5}
		\alpha_{\mathbb{H}}(B(x + y \mathbb{S}; \delta) \setminus K^\circ) \leq \alpha_{\mathbb{H}}(B(x + y \mathbb{S}; r\delta) \setminus K^\circ).
	\end{equation}
	
	Suppose that (iii) holds. Choose a fixed $r \geq 1$, for example $r = 2$. \\
	Then by (iii), for every $ q = x + y I \in \mathbb{H}$ and every $\delta > 0$,
	\[ \alpha_{\mathbb{H}}(B(x + y \mathbb{S}; \delta) \setminus K^\circ) \leq \alpha_{\mathbb{H}}(B(x + y \mathbb{S}; 2\delta) \setminus K). \]
	Then for $r = 2$ and $C = 1$, (iv) holds.
	
	Suppose that there exists $r \geq 1, C > 0$ such that for every $ q_0 \in \mathbb{H}$ and all $\delta > 0$,
	\begin{equation}\label{6}
		\alpha_{\mathbb{H}}(B(q, \delta) \setminus K^\circ) \leq C \alpha_{\mathbb{H}}(B(q_0, r\delta) \setminus K).
	\end{equation}
	Let $q \in \partial K$. Then,
	\[ \frac{\alpha_{\mathbb{H}}(B(q, \delta) \setminus K^\circ)}{\alpha_{\mathbb{H}}(B(q, r\delta) \setminus K)} \leq C. \]
	Taking the limit supremum as $\delta \to 0$ we get,
	\[ \limsup_{\delta \to 0} \frac{\alpha_{\mathbb{H}}(B(q, \delta) \setminus K^\circ)}{\alpha_{\mathbb{H}}(B(q, r\delta) \setminus K)} \leq C < \infty. \]
	Assume that (v) holds. Let $q \in \partial K$. Then there exists $r \geq 1$ such that,
	\[ \limsup_{\delta \to 0} \frac{\alpha_{\mathbb{H}}(B(q, \delta) \setminus K^\circ)}{\alpha_{\mathbb{H}}(B(q, \gamma\delta) \setminus K)} < \infty .\]
	Fix $I \in \mathbb{S}$ and set $E = K \cap \mathbb{C}_I$. 
	Since $K$ is axially symmetric, the sets $B(q, \delta) \setminus K^\circ$ and $B(q, \gamma\delta) \setminus K$ are axially symmetric. We already know that $\alpha_{\mathbb{H}}(E) = \alpha_{\mathbb{C}}(E \cap \mathbb{C}_I)$ for any axially symmetric set $E$. Thus, by Lemma \ref{12} 
	\[ \limsup_{\rho \to 0} \frac{\alpha_{\mathbb{C}}(\Delta(z, \rho) \setminus E^\circ)}{\alpha_{\mathbb{C}}(\Delta(z, \gamma\rho) \setminus E)} < \infty \]
	Thus, by Theorem \ref{vitalpha}, $\mathcal R_{\mathbb{C}}(E) = \mathcal A_{\mathbb{C}}(E)$ . Hence by Theorem \ref{10}, $\mathcal R_L(K) = \mathcal A_L(K).$
\end{proof}
\section{Quaternionic Version of Davie's Modified  Approximation scheme}
Vitushkin's approximation scheme  provides a method of decomposing functions into localized pieces and reduce global approximation problems to  local estimates using analytic capacity.  In 1972, Davie \cite{davie}, introduced a modified version of Vitushkin's scheme by  a variable-radius covering, which allows the size of the covering sets to varies with respect to the geometry of the underlying open set. 

In the complex plane, Davie formed this modified scheme using squares rather than discs which provides variable-radius covering of an open arbitrary open subset of the complex plane. It provides a covering with controlled overlap, an associated partition of unity with suitable derivative bounds, and a summability estimate.

In this section we develop a quaternionic analogue of Davie's modified approximation scheme in the setting of slice regular functions. In quaternionic setting we require a  replacement for squares which is  axial symmetric. Thus we constructed variable-radius axially symmetric tubular neighborhoods as the basic covering sets and using this we are trying to have an analogue of Proposition \ref{davie covering}. Also using this modified scheme to prove a quaternionic analogue of Vitushkin's approximation theorem. 

Let $q_r=(x_r,y_r)\in\mathbb R\times[0,\infty)$  and $\delta_r>0$ .
Define \[Q_r=Q(q_r,\delta_r) = \left\{ (x,y)\in\mathbb R\times[0,\infty):
|x-x_r|<\frac{\delta_r}{2},\quad |y-y_r|<\frac{\delta_r}{2} \right\}.\]
which is the open square centered at $q_r$ and having
side length $\delta_r$.
We lift $Q_r$ to the quaternionic space. Define.
\[\mathcal T_r=\left\{x+Iy:|x-x_r|<\frac{\delta_r}{2},\quad|y-y_r|<\frac{\delta_r}{2},\quad I\in\mathbb S\right\}.\]
Then, $\mathcal T_r$ is the axially symmetric tubular neighborhood generated by the square $Q_r$. The point $q_r=(x_r,y_r)$ in complex space corresponds to $[q_r]=\{x_r+Iy_r:I\in\mathbb S\}.$ Thus the centre of the quaternionic cell is not a single quaternionic point, but the whole $[q_r]$.

Now, \[\mathcal T_r\cap\mathbb C_I = \left\{x+Iy:|x-x_r|<\frac{\delta_r}{2},
\quad|y-y_r|<\frac{\delta_r}{2}\right\}.\]
Hence, on every complex slice, $\mathcal T_r\cap\mathbb C_I$ is 
a square of side length $\delta_r$. Thus, $\operatorname{diam}(\mathcal T_r\cap\mathbb C_I) \leq\sqrt{2}\,\delta_r.$
This slice-wise diameter, rather than the full diameter of $\mathcal T_r$,
is the quantity used in the localization estimates.

The sets $\mathcal T_r$ are axially symmetric by construction. For,
if $x+Iy\in \mathcal T_r$ and $J\in\mathbb S$, then $x+Jy\in \mathcal T_r.$ 
Hence $\mathcal T_r$ contains the complete spherical orbit corresponding
to every point of the square $Q_r$.


In the next Proposition we proved the  covering properties, bounded overlap, partition-of-unity estimates, and summability estimates as in Davie's modified approximation scheme [Proposition 3.1\cite{davie}].
The proof consists of applying Davie's variable-scale square
construction in a complex plane and then lifting the resulting construction to $\mathbb H$.
\begin{prop}\label{quaternionic-covering}
	Let $V\subset\mathbb H$ be an axially symmetric open set, and let $s:V\to(0,\infty)$ be continuous. Then there exist a countable family of axially symmetric tubular neighborhood $\{\mathcal T_r\}_{r\geq 1}$, positive numbers $\{\delta_r\}_{r\geq 1}$, where $\delta_r$ is the side length of the planar square generating $\mathcal T_r$ and determining the thickness of the corresponding tubular  neighbourhood, spherical centres $[q_r]$, and intrinsic-functions $\Phi_r\in C^1(V)$ such that:
	
	\begin{enumerate}
		\item $	V=\bigcup_r \mathcal T_r, $
		
		\item Every point of $V$ belongs to at most four sets $\mathcal T_r$.
		
		
		\item If $q\in \mathcal T_r$, then $\delta_r<s(q),$
		
		\item There is a constant $A>0$, independent of $q$, such that
		\[\sum_r\min\left\{1,\frac{\delta_r^3}{d(q,[q_r])^3}\right\}\leq A,
		\qquad q\in V,\]
		\item The functions $\Phi_r$ satisfy
		\[0\leq\Phi_r\leq1,\qquad\operatorname{supp}\Phi_r\subset \mathcal T_r,\qquad
		\sum_r\Phi_r(q)=1,\] 
		Moreover, for a constant $C$ independent of $r$, $\|\nabla\Phi_r\|_\infty\leq \frac{C}{\delta_r}.$		
	\end{enumerate}
	where $\delta_r$ is the side length of the planar square generating $\mathcal T_r$ and determining the thickness of the corresponding tubular  neighborhood.
\end{prop}

\begin{proof}
	Fix $I\in\mathbb S$ and define
	\[\widetilde\Omega=\{(x,y)\in\mathbb R^2:x+I|y|\in V\}.\]
	Since $V$ is axially symmetric, it does not depend
	on the choice of $I\in\mathbb S$.
	The set $\widetilde\Omega$ is open and invariant under the
	reflection $R(x,y)=(x,-y).$
	Define
	\[\widetilde p(x,y)=\min_{J\in\mathbb S}s(x+J|y|).\]
	Since  $\{x+J|y|:J\in\mathbb S\}$ is compact, the minimum is attained. By the continuity of $s$ and compactness of $\mathbb S$ we have $\widetilde p$ is 	continuous and positive. Also, $\widetilde p(x,-y)=\widetilde p(x,y).$
	
	Now by applying Davie's  variable-scale covering  Proposition \ref{davie covering} construction to
	$(\widetilde\Omega,\widetilde p)$, it gives a countable family 	of enlarged squares with  properties as given in Proposition \ref{davie covering}.
	
	Since, $\widetilde\Omega$ and $\widetilde p$ are invariant
	under $R$, the construction can be performed in  such a way that, whenever $Q$ is selected, its reflection $RQ$ is selected as well. Choose one representative $Q_r$ from each
	pair $\{Q,RQ\}$, and write $z_r=(x_r,y_r),\ y_r>0,$
	for its centre, and $\delta_r$ for its side length. Since the centres of  squares used in the construction is of the form $z_{m,n} = \big(\frac{m+\frac 1 2}{2^k}, \frac{n+\frac 1 2}{2^k}\big)$ for $m,n \in \mathbb Z$, no centre lies on the real axis $y=0$, so such representatives can be chosen.
	
	Let $Q_r'=Q(z_r,\delta_r/3)$ be the corresponding inner square. Define
	\[\mathcal T_r=\{x+Jy:(x,|y|)\in Q_r\cup RQ_r,\ J\in\mathbb S\},\]
	and
	\[\mathcal T_r'=\{x+Jy:(x,|y|)\in Q_r'\cup RQ_r',\ J\in\mathbb S\}.\]
	
	We verify the required properties.
	For, take $q=x+ y'\in V$, where $y'=\operatorname{Im}q$, and put $z=(x, |y'|).$
	Since $V$ is axially symmetric, $z\in\widetilde\Omega$. The planar covering property gives a selected square $Q$ containing $z$. For some $r$, either $Q=Q_r$ or $Q=RQ_r$. Hence $z\in Q_r\cup RQ_r,$	and thus $q\in \mathcal T_r$. Thus, $V\subseteq \bigcup_r\mathcal T_r$. 
	
	For the converse part, let $q=x+Jy \in \mathcal T_r$ for some $r$ where $J\in \mathbb S$. Then by the definition of $\mathcal T_r$, $(x, |y|) \in Q_r \cup RQ_r$. Since $Q_r$ and $RQ_r$ are selected squares contained in $\widetilde{\Omega}$, we have $(x, |y|)\in \widetilde{\Omega}$. Thus by the definition of $\widetilde{\Omega}$, we have $x+I|y| \in V$. Now, since $V$ is axially symmetric $x+J|y| \in V$. Hence $\bigcup_r\mathcal T_r \subseteq V.$ Thus, 
	\[V=\bigcup_r\mathcal T_r.\]

	If $q = x+y'\in \mathcal T_r$, then $z =(x, |y'|)$ belongs to one of the two selected squares $Q_r$ and $RQ_r$. Now suppose that $q$ belongs to $m$ different tubes, which is from $m$ pairs $\{Q_{r_1}, RQ_{r_1}\}, \cdots, \{Q_{r_m}, RQ_{r_m}\}$. Since we chose only one representative from each pair, these pairs are distinct. Thus $z$ belongs to at least $m$ distinct selected squares in the  reflected family. But, by Proposition \ref{davie covering}, atmost four selected squares can contain the same point $z$. Hence $m\leq 4$. Hence every point of $V$ belongs to atmost four sets $\mathcal T_r$.
	
	Let $q=x+ y'\in \mathcal T_r$. Then, $z=(x,| y'|) \in Q_r\cup RQ_r$. If $z \in Q_r$, by Proposition \ref{davie covering} we have, $\delta_r <\widetilde{p}(z)$. Similarly if $z \in RQ_r$, we have $\delta_r <\widetilde{p}(Rz)$ where $Rz = (x, -|y'|).$  But, since $\widetilde p(x,-y)=\widetilde p(x,y),$ , $\widetilde{p}(Rz) = \widetilde{p}(z).$ Thus, we have
	$\delta_r<\widetilde p(z).$
	Now, by the definition of $\widetilde p$,
	\[\delta_r<\min_{J\in\mathbb S}s(x+J|y|)\leq s(q).\]
	Thus if $q\in \mathcal T_r$ then, $\delta_r<s(q).$
	
	We know that the distance from $q$ to the spherical class
	$[q_r]=\{x_r+Jy_r:J\in\mathbb S\}$ 	is
	\[d(q,[q_r])=\inf_{J\in\mathbb S}|q-(x_r+Jy_r)|.\]
	Write $ y=I|y|$ when $ y\neq0$, we have
	\[|q-(x_r+Jy_r)|^2=(x-x_r)^2+|y|^2+y_r^2-2|y|y_r\langle I,J\rangle.\]
	Since $\sup_{J\in\mathbb S}\langle I,J\rangle=1,$
	the minimum is attained at $J=I$. The same distance formula
	holds when $ y=0$. Thus,
	\[d(q,[q_r])^2=(x-x_r)^2+(|y|-y_r)^2,\]
	or equivalently, $d(q,[q_r])=|z-z_r|.$
	Hence by Proposition \ref{davie covering}
	\[\sum_r\min\left\{1,\frac{\delta_r^3}{d(q,[q_r])^3}\right\}
	= \sum_r\min\left\{1,\frac{\delta_r^3}{|z-z_r|^3}\right\}\leq A.
	\]
	
	Next it is remain to prove (5).	We first construct the partition in the plane. Choose the 	cutoff functions in Davie's planar construction compatibly with reflection. For each representative $Q_r$, let
	$\psi_r$ be a $C^1$ cutoff supported in $Q_r$, with
	$0\leq\psi_r\leq1$, equal to $1$ on a neighborhood of the
	original square, and satisfying $\|\nabla\psi_r\|_\infty\leq \frac{C_0}{\delta_r}.$
	Set, \[ \psi_{Rr}(x,y)=\psi_r(x,-y).\]
	For the full selected family, define
	\[	\Psi(x,y)=\sum_r\bigl(\psi_r(x,y)+\psi_r(x,-y)\bigr).\]
	Then by the planar partition construction gives $\Psi(x,y)\geq1,\ (x,y)\in\widetilde\Omega.$
	
	Now define
	\[\varphi_r(x,y)=\frac{\psi_r(x,y)+\psi_r(x,-y)}{\Psi(x,y)}.\]
	The sum is locally finite, so $\Psi$ and each $\varphi_r$ are
	$C^1$. Since $\Psi\geq1$, the quotient is well-defined. Also,
	$0\leq\varphi_r\leq1,\ \operatorname{supp}\varphi_r\subset Q_r\cup RQ_r,$  and $\varphi_r(x,-y)=\varphi_r(x,y).$
	Also,
	\[\sum_r\varphi_r(x,y)=\frac{\displaystyle\sum_r\bigl(\psi_r(x,y)+\psi_r(x,-y)\bigr)}{\Psi(x,y)}=1.\]
	
	In Davie's construction, if two selected enlarged
	squares intersect, their side lengths are comparable by a
	uniform constant. In particular, if the support of $\varphi_r$
	meets the support of another cutoff, the corresponding side
	lengths are comparable. Together with the multiplicity bound,
	this gives $\|\nabla\Psi\| \leq \frac{C_1}{\delta_r}$
	on the support of $\varphi_r$, where $C_1$ is independent of
	$r$. Differentiating the quotient gives
	\[\nabla\varphi_r=\frac{\nabla(\psi_r+\psi_r\circ R)}{\Psi}
	-\frac{(\psi_r+\psi_r\circ R)\nabla\Psi}{\Psi^2}.\]
	Since, $\Psi\geq1$ and $0\leq\psi_r+\psi_r\circ R\leq\Psi,$	 we have,
	\[\|\nabla\varphi_r\|_\infty\leq \frac{C_2}{\delta_r}\]
	for a constant $C_2$ independent of $r$. Now, define 
	\[\Phi_r(q)=\varphi_r\bigl(\operatorname{Re}q,|\operatorname{Im}q|\bigr),\qquad q\in V.\]
	Then by construction, $0\leq\Phi_r(q)\leq1,\ \operatorname{supp}\Phi_r\subset \mathcal T_r.$
	Also,
	\[\sum_r\Phi_r(q)=\sum_r\varphi_r\bigl(\operatorname{Re}q,|\operatorname{Im}q|\bigr)=1.\]
	Since $\Phi_r$ depends only on $\operatorname{Re}q$ and
	$|\operatorname{Im}q|$, it is constant on every spherical class $[q_r]$. It remains to check $C^1$ regularity at real points. Write $\varphi_r=\varphi_r(x,y).$ Since, $\varphi_r(x,-y)=\varphi_r(x,y)$ we have $\partial_y\varphi_r(x,0)=0.$ Now, for $q=x+ y'$ with $ y'\neq0$, 
	\[\nabla\Phi_r(q)=\left(\partial_x\varphi_r(x,|y|),\partial_y\varphi_r(x,|y|)\frac{y'}{|y|}\right).
	\]
	Since, $\partial_y\varphi_r(x,0)=0$ and
	$\varphi_r$ is $C^1$, this gradient extends continuously to $y'=0$, with
	\[\nabla\Phi_r(x)=\bigl(\partial_x\varphi_r(x,0),0,0,0\bigr).\]
	Thus $\Phi_r\in C^1(V)$, and $\|\nabla\Phi_r\|_\infty\leq\frac{C}{\delta_r}.	$
\end{proof}
From the proof, it is clear that, $\mathcal T_r'\cap T_t=\varnothing,\qquad r\neq t.$ For,Suppose $r\neq t$. Now, by Proposition \ref{davie covering},
$Q_r'$ does not meet any selected square other than $Q_r$.
Similarly, for its reflection $RQ_r'$.
Hence, \[(Q_r'\cup RQ_r')\cap(Q_t\cup RQ_t)=\varnothing.\]
in classical complex analysis, Vitushkin's localization operator is an important tool in approximation theory. It localizes a continuous function using a smooth cutoff function while preserving holomorphicity away from the support of the cutoff. This allows one to study approximation locally and then combine the local approximations. 
 
 Here on, we follows the proof techniques and scheme of \cite{rational app,analytic capacity} for extending Vitushkin approximation Theorem.
Here we provide the foundation for extending Vitushkin localization method to  slice regular functions and applying to approximation problem.
\begin{lem}\label{quaternionic-localization}
	Let $\mathcal T_r$ be an axially symmetric tubular neighborhood as in
	Proposition \ref{quaternionic-covering}, and let $\Phi_r$ be the corresponding intrinsic
	continuously differentiable function with compact support.
	Let $f$ be a continuous slice function. For $q\in\mathbb H$,
	choose $I\in\mathbb S$ such that $q\in\mathbb C_I$ and define
	\[
	T_{\Phi_r}^{L}f(q)
	=
	\Phi_r(q)f(q)
	-\frac{1}{2\pi}
	\iint_{\mathbb C_I}
	S_L^{-1}(s,q)
	\bigl(\bar\partial_I\Phi_r(s)f(s)\bigr)
	\,ds_I\wedge d\bar s.
	\]
	Then the following assertions hold:
	\begin{enumerate}
		\item $T_{\Phi_r}^{L}f$ is continuous on
		$\widehat{\mathbb H}$ and $T_{\Phi_r}^{L}f(\infty)=0$;
		\item $T_{\Phi_r}^{L}f$ is left slice regular on
		$ \widehat{\mathbb H}\setminus\operatorname{supp}\Phi_r$ and wherever $f$ is left slice regular;
		\item $	\|T_{\Phi_r}^{L}f\|_\infty\leq C\delta_r
		\|\bar\partial_I\Phi_r\|_\infty\omega_f(T_r)\leq C'\omega_f(T_r),$ where				
		$$\omega_f(T_r)=\sup_{p,q\in T_r}|f(p)-f(q)|$$
		and the constants $C,C'>0$ are independent of $f$ and $r$.
	\end{enumerate}
\end{lem}

\begin{proof}
	We have to prove that $T_{\Phi_r}^{L}f$ is continuous. Since, $\Phi_r$ is intrinsic and $f$ is a continuous slice function, $\Phi_r f$ is also a continuous slice function. Hence it is enough to prove that the integral part is continuous slice function.
	Since, $f$ is continuous slice function  
	\[	f(x+Iy)=\alpha(x,y)+I\beta(x,y),	\]
	where $	\alpha(x,-y)=\alpha(x,y),\qquad\beta(x,-y)=-\beta(x,y)$ and they are continuous. Since $\Phi_r$ is intrinsic, we  write
	$\Phi_r(x+Iy)=\phi(x,y),$ where $\phi$ is real-valued. 
	Thus, 
	\[(\bar\partial_I\Phi_r)f=\frac12(A+IB),\]
	where $A=\phi_x\alpha-\phi_y\beta,$ and $B=\phi_x\beta+\phi_y\alpha.$
	Hence, 
	\begin{equation}\label{eqnab}
		A(x,-y)=A(x,y),\qquad B(x,-y)=-B(x,y).
	\end{equation}
	For $q=x+Iy$ and $s=u+Iv$, the slice Cauchy kernel on $\mathbb C_I$ is
	\[S_L^{-1}(s,q)=\frac{(u-x)-I(v-y)}{(u-x)^2+(v-y)^2}.\]
	Thus,
	\[\begin{aligned}S_L^{-1}(s,q)(A+IB)={}&\frac{(u-x)A+(v-y)B}
		{(u-x)^2+(v-y)^2} +I\frac{(u-x)B-(v-y)A}{(u-x)^2+(v-y)^2}.
	\end{aligned}\]
	
	Integrating the two components and by \eqref{eqnab} with the substitution $v\mapsto -v$,  the integral in the definition of $T_{\Phi_r}^{L}f$
	can be expressed in terms of two component integrals:
	\[F(x+Iy)=U(x,y)+IV(x,y).\]
	where, $U(x,-y)=U(x,y),\qquad	V(x,-y)=-V(x,y).$
	Thus it is a slice function and it is clear that $U$ and $V$ are continuous.			
	Thus the integral term is a continuous slice function. Hence, $T_{\Phi_r}^{L}f$ is continuous on $\mathbb H$.
	
	Since $\Phi_r$ has compact support, $\Phi_r(q)=0$ for
	sufficiently large $|q|$ and the function
	$\bar\partial_I\Phi_r$ also has compact support. Hence, there
	exists $M>0$ such that $|s|\le M$ for all $s$ in the
	integration support. Since, for $s,q\in\mathbb C_I$, 
	\begin{align*}
		|S_L^{-1}(s,q)|&=\frac{1}{|s-q|}\\
		&\leq |q|-M\  \leq \frac{2}{|q|}.
	\end{align*}
	for sufficiently large $|q|\geq 2M$.
	
	Since $\Phi_r(q)=0$ for sufficiently large $|q|$, we have
	\[
	\begin{aligned}
		|T_{\Phi_r}^L f(q)|
		&\le \frac{C}{|q|}
		\iint_{\mathbb C_I}
		\left|(\bar\partial_I\Phi_r(s))f(s)\right|
		\,dA_I(s)
		\rightarrow 0,
		\quad\text{as }|q|\longrightarrow\infty.
	\end{aligned}
	\]
	The integral on the right is finite and independent of $q$,
	since  it  has compact support. Hence,
	$T_{\Phi_r}^L f$ extends continuously to
	$\widehat{\mathbb H}$ by setting $T_{\Phi_r}^L f(\infty)=0.$	This proves (1).
	
	Let $q_0\in\mathbb H\setminus\operatorname{supp}\Phi_r.$
	Since $\operatorname{supp}\Phi_r$ is closed, there exists a
	neighborhood $V$ of $q_0$ such that $V\subset\mathbb H\setminus\operatorname{supp}\Phi_r.$ On $V$, we have $\Phi_r(q)=0$, and hence
	\[T_{\Phi_r}^{L}f(q)=-\frac{1}{2\pi}\iint_{\mathbb C_I}
	S_L^{-1}(s,q)\bigl(\bar\partial_I\Phi_r(s)f(s)\bigr)\,ds_I\wedge d\bar s.\]
	Since $V$ is disjoint from $\operatorname{supp}\Phi_r$, the Cauchy
	kernel has no singularities in the region of integration for
	$q\in V$.  As $S_L^{-1}(s,q)$ is left slice regular in $q$
	away from its singularities, we have $\bar\partial_{I,q}S_L^{-1}(s,q) =0$. Hence, 
	\[\bar\partial_q T_{\Phi_r}^{L}f(q)=0,\qquad q\in V.\]
	Thus $T_{\Phi_r}^{L}f$ is left slice regular on
	$\mathbb H\setminus\operatorname{supp}\Phi_r.$			
	Since by the above part $\lim_{q\to\infty}T_{\Phi_r}^{L}f(q)$ exists, $T_{\Phi_r}^{L}f$ is left slice regular at infinity.
	
	Suppose that $f$ is left slice regular in a neighborhood
	$W$ of a point $q_0$. On each slice $\mathbb C_I$, we have
	\[\bar\partial_I f=0 \qquad \text{in } W\cap\mathbb C_I.\]
	Since $\Phi_r$ is continuously differentiable, 
	\begin{equation}\label{eq:product-rule}
		\begin{aligned}
			\bar\partial_I(\Phi_r f)
			&=(\bar\partial_I\Phi_r)f
			+\Phi_r\bar\partial_I f=(\bar\partial_I\Phi_r)f
		\end{aligned}
	\end{equation}
	in $W\cap\mathbb C_I$.
	
	Choose a sufficiently small slice domain
	$D\subset W\cap\mathbb C_I$ containing $q_0$ and split the integral defining $T_{\Phi_r}^L f$ into its
	portions over $D$ and $\mathbb C_I\setminus D$.
	The latter portion is slice regular in a smaller neighborhood
	of $q_0$, since its integration support is separated from
	that neighborhood.
	
	Applying the Cauchy--Pompeiu formula to $\Phi_r f$ on $D$
	and using \eqref{eq:product-rule}, we have
	\begin{equation*}\label{eq:local-cauchy-pompeiu}
		\begin{aligned}
			\Phi_r(q)f(q)
			={}&\frac{1}{2\pi}
			\int_{\partial D}
			S_L^{-1}(s,q)\,ds_I\,\Phi_r(s)f(s)\\
			&+\frac{1}{2\pi}
			\iint_D S_L^{-1}(s,q)
			(\bar\partial_I\Phi_r(s))f(s)
			\,ds_I\wedge d\bar s .
		\end{aligned}
	\end{equation*}
	Subtracting the area integral over $D$ from  $T_{\Phi_r}^L f$ gives,
	\[	\begin{aligned}
		T_{\Phi_r}^L f(q)={}&\frac{1}{2\pi}
		\int_{\partial D}S_L^{-1}(s,q)\,ds_I\,\Phi_r(s)f(s)\\
		&-\frac{1}{2\pi}\iint_{\mathbb C_I\setminus D}
		S_L^{-1}(s,q)(\bar\partial_I\Phi_r(s))f(s)
		\,ds_I\wedge d\bar s .
	\end{aligned}
	\]
	The boundary integral is slice regular in $D$, while the
	remaining area integral is slice regular in a neighborhood
	of $q_0$, because its integration support is separated from
	that neighborhood. Hence, $T_{\Phi_r}^L f$ is left
	slice regular near $q_0$. Since $q_0$ was arbitrary, $T_{\Phi_r}^L f$ is left slice regular wherever $f$ is left slice regular.

	Since $\Phi_r$ is continuously differentiable and has compact
	support, the quaternionic Cauchy-Pompeiu formula gives
	\[
	\Phi_r(q)=\frac{1}{2\pi} \iint_{\mathbb C_I}S_L^{-1}(s,q)
	\bar\partial_I\Phi_r(s)\,ds_I\wedge d\bar s.\]
	Let $c\in\mathbb H$ be constant. Then,
	\[\Phi_r(q)c=\frac{1}{2\pi}\iint_{\mathbb C_I}
	S_L^{-1}(s,q)\bigl(\bar\partial_I\Phi_r(s)c\bigr)\,ds_I\wedge d\bar s.\]
	Subtracting this identity from the definition of
	$T_{\Phi_r}^{L}f$, we have
	\[\begin{aligned}
		T_{\Phi_r}^{L}f(q)={}&\Phi_r(q)(f(q)-c)\\
		&-\frac{1}{2\pi}\iint_{\mathbb C_I}
		S_L^{-1}(s,q)(\bar\partial_I\Phi_r(s))(f(s)-c)\,ds_I\wedge d\bar s.
	\end{aligned}\]
	Choose $q_0\in T_r$ and put $c=f(q_0)$. Then
	\[|f(s)-c|\le \omega_f(T_r),\qquad s\in T_r.\]
	Thus, the integral term satisfies
	\begin{equation}\label{bound}
		\begin{aligned}
			&\left|\frac{1}{2\pi}\iint_{\mathbb C_I}
			S_L^{-1}(s,q)(\bar\partial_I\Phi_r(s))
			(f(s)-c)\,ds_I\wedge d\bar s\right|\\
			&\leq \frac{\omega_f(T_r)}{2\pi}
			\iint_{\operatorname{supp}\bar\partial_I\Phi_r}
			|S_L^{-1}(s,q)|\,|\bar\partial_I\Phi_r(s)|\,dA_I(s).
		\end{aligned}
	\end{equation}
	By the construction of the tubular neighborhoods,
	the planar set $\operatorname{supp}\Phi_r\cap\mathbb C_I$
	has at most two connected components, each of diameter at most
	$C_0\delta_r$. Thus, by [\cite{analytic capacity} Proposition 11.2] to each component, and using \eqref{bound}
	\begin{equation}\label{complexbound}
		\sup_{q\in\operatorname{supp}\Phi_r\cap\mathbb C_I}
		|T_{\Phi_r}^{L}f(q)|\leq C_1\delta_r
		\|\bar\partial_I\Phi_r\|_\infty\omega_f(T_r),
	\end{equation}
	where $C_1>0$ is independent of $f$ and $r$.
	
	Let $F=T_{\Phi_r}^{L}f.$ we already proved that $F$ is continuous on $\widehat{\mathbb H}$, left slice regular on $\widehat{\mathbb H}\setminus\operatorname{supp}\Phi_r$, and satisfies $F(\infty)=0$. Fix $I\in\mathbb S$. By the maximum modulus principle, to the connected components of
	$\mathbb C_I\setminus\bigl(\operatorname{supp}\Phi_r\cap\mathbb C_I\bigr),$	\eqref{complexbound} extends to the whole slice $\mathbb C_I$. Indeed, the boundary of each bounded complementary component lies in $\operatorname{supp}\Phi_r\cap\mathbb C_I$, while for the unbounded component we have $F(\infty)=0$.
	Thus,
	\begin{equation}\label{sliceestimate}
		\sup_{z\in\mathbb C_I}|F(z)|\leq C_1\delta_r
		\|\bar\partial_I\Phi_r\|_\infty\,\omega_f(T_r).
	\end{equation}
	Now let $q=x+Jy\in\mathbb H\setminus\operatorname{supp}\Phi_r$, where $J\in\mathbb S$. Since $\operatorname{supp}\Phi_r$ is axially symmetric, both $x+Iy$ and $x-Iy$ belong to
	$\mathbb H\setminus\operatorname{supp}\Phi_r$.
	Thus by the representation formula .
	\[F(x+Jy)=\frac12(1-JI)F(x+Iy)+\frac12(1+JI)F(x-Iy)\]
	and thus
	\[|F(x+Jy)|\leq |F(x+Iy)|+|F(x-Iy)|.\]
	Hence by \eqref{sliceestimate}, 
	\[\sup_{q\in\mathbb H\setminus\operatorname{supp}\Phi_r}|F(q)|
	\leq C_2\delta_r\|\bar\partial_I\Phi_r\|_\infty\,\omega_f(T_r),	\]
	where $C_2>0$ is independent of $f$ and $r$.
	Thus,
	\[\|T_{\Phi_r}^{L}f\|_\infty\le C_2\delta_r\|\bar\partial_I\Phi_r\|_\infty\,\omega_f(T_r).\]
	Hence by Proposition \ref{quaternionic-covering}
	$\|T_{\Phi_r}^{L}f\|_\infty\le C'\omega_f(T_r),$
	where $C'>0$ is independent of $f$ and $r$.
\end{proof}
\begin{lem}
	\label{quaternionic 11.3}
	Let $K\subset\mathbb H$ be a compact axially symmetric set and let
	$f\in \mathcal A_L(K).$	Let $\mathcal T_r$ be one of the axially symmetric tubular neighborhoods $\Phi_r$ be the continuously differentiable function as constructed in Proposition \ref{quaternionic-covering}.
	Then,
	\begin{equation*}
		\label{quaternionic-tubular}
		\left|	\frac{1}{2\pi}		\int_{\mathbb C_I}
		\overline{\partial}_I\Phi_r(s)f(s)\,ds_I\wedge d\bar s\right|
		\leq 	C_1\,\delta_r\, 		\omega_f\!\left(\mathcal T_r^{(2)}\right)\left\|\overline{\partial}_I\Phi_r\right\|_\infty
		\alpha_{\mathbb H}\!\left(\mathcal T_r\setminus K^\circ\right),
	\end{equation*}
	and for $p_{r}=x_r+Iy_r\in [q_r]$, $\mathcal T_r \cap \mathbb C_I = Q_r,$
	\[\left|\frac{1}{2\pi}\int_{C_I}(s-z_r)(\overline{\partial}_I\Phi_r(s))f(s)\,	ds_I\wedge d\bar{s}\right| 	\leq C\delta_r^2\,\omega_f(T_r^{(2)})
	\|\overline{\partial}_I\Phi_r\|_\infty \alpha_{\mathbb H}(T_r\setminus K^\circ),\]
	where $\mathcal T_r^{(2)} = \{x+Jy:(x,|y|)\in Q_r^2\cup RQ_r^2,\ J\in\mathbb S\},$ for $Q_r^2=Q(z_r,2\delta_r)$ and $C_1, C>0$ is an absolute constant.
\end{lem}

\begin{proof}
	Let $F=T_{\Phi_r}^{L}f$. 
	If $F\equiv0$, then the Laurent coefficient at infinity
	vanishes, and the desired inequality follows immediately.
	We  assume that $F\not\equiv0$. Define $G_r(q)=\frac{F(q)}{\|F\|_\infty}.$
	Then, by Proposition \ref{quaternionic-localization} $G_r\in \mathcal C_L(\widehat{\mathbb H}),\  G_r(\infty)=0$ and $ \|G_r\|_\infty=1.$	Also, $G_r$ is left
	slice regular on $\widehat{\mathbb H}\setminus E_r$ where, $E_r=\mathcal T_r\setminus K^\circ.$		
	Thus $F_r$ is admissible in the definition of
	$\alpha_{\mathbb H}(E_r)$. Hence, $	|\partial_sG_r(\infty)|
	\leq \alpha_{\mathbb H}(E_r)$ and thus,
	\begin{equation}\label{derivative}
		|\partial_sF(\infty)|
		\leq \|F\|_\infty\alpha_{\mathbb H}(E_r).
	\end{equation}
	Since $\Phi_r$ has compact support, there exists $R>0$ such that
	$\Phi_r(q)=0$ for $|q|>R.$	Hence, for sufficiently large $|q|$,
	\[F(q)=-\frac{1}{2\pi}\int_{\mathbb C_I}S_L^{-1}(s,q)
	(\bar\partial_I\Phi_r(s))f(s)\,ds_I\wedge d\bar s.\]
	Since, by Lemma \ref{cauchy-kernel-infinity} $S_L^{-1}(s,q)=q^{-1}+O(|q|^{-2}),$ as $|q|\longrightarrow\infty,$ 
	we have
	\[\begin{aligned}
		F(q)=-\frac{1}{2\pi q}\int_{\mathbb C_I}(\bar\partial_I\Phi_r(s))f(s)
		\,ds_I\wedge d\bar s +O(|q|^{-2}).
	\end{aligned}\]
	and hence
	\[\partial_sF(\infty)=-\frac{1}{2\pi}\int_{\mathbb C_I}
	(\bar\partial_I\Phi_r(s))f(s)\,ds_I\wedge d\bar s.\]
	Thus by, \eqref{derivative}
	\[\begin{aligned}\left|\frac{1}{2\pi}\int_{\mathbb C_I}
		(\bar\partial_I\Phi_r(s))f(s)\,ds_I\wedge d\bar s
		\right|\leq\|F\|_\infty\alpha_{\mathbb H}(E_r).
	\end{aligned}\]
	Hence by Proposition \ref{quaternionic-localization}
	\[	\begin{aligned}
		\left|\frac{1}{2\pi}\int_{\mathbb C_I}(\bar\partial_I\Phi_r(s))f(s)
		\,ds_I\wedge d\bar s\right| \leq C_1\delta_r\omega_f({\mathcal T_r}^{(2)})
		\|\bar\partial_I\Phi_r\|_\infty\alpha_{\mathbb H}(\mathcal T_r\setminus K^\circ).
	\end{aligned}\]
	Let $g(s)=(s-z_r)\Phi_r(s).$ Then, 
	$\overline{\partial}_I g(s)
	=(s-z_r)\overline{\partial}_I\Phi_r(s).$
	Hence, By the geometry of the tubular neighborhood $\mathcal T_r$,
	\[\|\overline{\partial}_I g\|_\infty\leq\left(\sup_{s\in T_r}|s-z_r|\right)\|\overline{\partial}_I\Phi_r\|_\infty\leq C_1\delta_r
	\|\overline{\partial}_I\Phi_r\|_\infty..\]
	
	Hence, by above part with  $g$ in place of $\Phi_r$, we have
	\begin{align*}
		\left|
		\frac{1}{2\pi}\int_{C_r}
		(s-z_r)(\overline{\partial}_I\Phi_r(s))f(s)\,
		ds_I\wedge d\bar{s}
		\right|&=\left|
		\frac{1}{2\pi}\int_{C_r}
		(\overline{\partial}_I g(s))f(s)\,
		ds_I\wedge d\bar{s}
		\right| \\
		&\leq
		C\delta_r^2\,\omega_f(T_r^{(2)})
		\|\overline{\partial}_I\Phi_r\|_\infty
		\alpha_{\mathbb H}(T_r\setminus K^\circ).
	\end{align*}
	where \(C>0\) is independent of \(r\).
\end{proof}
In Theorem \ref{esticomp} we proved the pointwise estimate of admissible functions in the definition of slice regular capacity. We proved its counterpart for the admissible  continuous slice-regular capacity functions following \cite{analytic capacity}. 
\begin{lem}
	\label{estimate pointwise}
	Let $S\subset\mathbb{H}$ be an axially symmetric set, and let $f$ be admissible for the continuous slice-regular capacity($\alpha_{\mathbb{H}}$) of 	$S$.
	Then, for every $[q_0]\notin S$,
	\[|f(q_0)| \leq\frac{4\,\alpha_{\mathbb H}(S)}{d_{\mathbb H}([q_0],S)}.\]
\end{lem}

\begin{proof}
	Fix $q_0\notin S$ and let $d=d_{\mathbb H}(q_0,S)>0.$	
	Choose $I\in\mathbb{S}$ such that $q_0\in\mathbb{C}_I$.
	Consider, 
	\[R_{f,q_0}(q) = (q-q_0)^{-*}*(f(q)-f(q_0)).\]
	Since $f$ is slice regular in a neighborhood of $q_0$, the function $h(q)=f(q)-f(q_0)$	is slice regular and satisfies $h(q_0)=0$. By the factorization theorem for slice-regular functions, there exists a slice regular function $g$ such that $f(q)-f(q_0)=(q-q_0)*g(q).$ Then, 
	\begin{align*}
		R_{f,q_0}(q)
		&=(q-q_0)^{-*}*
		\bigl((q-q_0)*g(q)\bigr)\\
		&=\bigl((q-q_0)^{-*}*(q-q_0)\bigr)*g(q)=g(q).
	\end{align*}
	Thus $R_{f,q_0}$ admits a slice-regular extension across $q_0$; $
	R_{f,q_0}(q_0)=g(q_0).$	 Thus, the singularity at $q=q_0$ is removable	.
	In particular, on the slice $\mathbb{C}_I$,  we have
	\[R_{f,q_0}(z)=\frac{f(z)-f(q_0)}{z-q_0}, \ z\in\mathbb{C}_I\setminus S.\]
	Since,  restriction of $f$ to $\mathbb{C}_I$ is holomorphic near
	$q_0$,
	\[f(z)=f(q_0)+f_I'(q_0)(z-q_0)+O\bigl((z-q_0)^2\bigr).\]
	Thus, 
	\[\lim_{z\to q_0}	\frac{f(z)-f(q_0)}{z-q_0}= \lim_{z\to q_0}f_I'(q_0)+O(z-q_0)=f_I'(q_0).\]
	Thus the singularity is removable on the slice $\mathbb{C}_I$,
	consistently with the quaternionic extension above.
	
	Since $f(\infty)=0$, we have $f(q)=O(q^{-1})$ as $q\to\infty.$
	Thus,
	\[f(q)-f(q_0)=-f(q_0)+O(q^{-1}).\]
	Also, since, 
	\[(q-q_0)^{-*}=(q-\overline{q_0})\bigl(q^2-2\operatorname{Re}(q_0)q+|q_0|^2\bigr)^{-1}.\]
	as $q\to\infty$, $(q-q_0)^{-*}=q^{-1}+O(q^{-2}).$
	
	Hence, we have
	\begin{align*}
		R_{f,q_0}(q)
		&=(q-q_0)^{-*}*(f(q)-f(q_0))\\
		&=\bigl(q^{-1}+O(q^{-2})\bigr)
		\bigl(-f(q_0)+O(q^{-1})\bigr)\\
		&=-q^{-1}f(q_0)+O(q^{-2}).
	\end{align*}
	Hence, the first Laurent coefficient of $R_{f,q_0}$ at infinity is $-f(q_0).$
	
	Suppose first that $	|z-q_0|\geq d,$ for $z\in\mathbb{C}_I\setminus\mathcal{T}.$ Since $\|f\|_\infty\leq1$,
	\[|f(z)-f(q_0)| \leq |f(z)|+|f(q_0)| \leq2.\]
	Hence
	\begin{equation}\label{bound}
		|R_{f,q_0}(z)| = 	\frac{|f(z)-f(q_0)|}{|z-q_0|}
		\leq	\frac{2}{d} =	\frac2d.
	\end{equation}
	Since $q_0\notin S$, we have $B_{\mathbb H}\left(q_0, d\right)
	\subset \mathbb{H}\setminus S$ and thus $\Delta_{\mathbb C_I}\left(q_0, d\right)\subset \mathbb{C}_I\setminus S	.$
	Since, the singularity of $R_{f,q_0}$ at $q_0$ is removable, 
	$R_{f,q_0}$ is holomorphic on $\Delta_{\mathbb C_I}\left(q_0, d\right)$ and by \eqref{bound} on the boundary of this disk,
	\[|R_{f,q_0}(z)|\leq\frac2d.\]
	Thus by the maximum-modulus principle,
	\[|R_{f,q_0}(z)|\leq\frac2d,\ \text{for}\	|z-q_0|\leq d.	\]
	Hence, $|R_{f,q_0}(z)|\leq\frac2d$ for $z\in\mathbb{C}_I\setminus S.$
	
	We now pass from the slice $\mathbb{C}_I$ to the whole quaternionic
	domain. Let $q=x+yJ\in\mathbb{H}\setminus S,\  J\in\mathbb{S},$	and define $z_+=x+yI,\ z_-=x-yI.$ Since $S$ is axially symmetric, if $q\notin S$, then $z_+,z_-\notin S.$	 Also by the representation formula for slice-regular functions 
	\[\begin{aligned}
		R_{f,q_0}(q)
		=\frac12\left(R_{f,q_0}(z_+)+R_{f,q_0}(z_-)\right)+ 	\frac{JI}{2}\left(R_{f,q_0}(z_-)-R_{f,q_0}(z_+)\right).
	\end{aligned}
	\]
	Thus by above part in slice we have,
	\[	\begin{aligned}
		|R_{f,q_0}(q)|&\leq\frac12\left(|R_{f,q_0}(z_+)|
		+|R_{f,q_0}(z_-)|\right) + \frac12\left(|R_{f,q_0}(z_-)|+|R_{f,q_0}(z_+)|\right)\\
		&=|R_{f,q_0}(z_+)|+|R_{f,q_0}(z_-)|\\
		&\leq\frac2d+\frac2d=\frac{4}{d}.
	\end{aligned}
	\]
	Thus $\|R_{f,q_0}\|_{L^\infty(\mathbb{H}\setminus S	)}
	\leq\frac{4}{d}.$
	
	Now, define $g(q)=\frac d{4}R_{f,q_0}(q).$
	Then, $\|g\|_{L^\infty(\mathbb{H}\setminus S)}\leq1.$
	Also, $	g(\infty)=0.$	Hence $g$ is admissible for the continuous slice-regular capacity of $ S	$. Thus, the first Laurent coefficient of $g$ is 
	\[a_1(g)  = \frac d{4}a_1(R_{f,q_0}) = -\frac d{4}f(q_0).\]
	Thus, by the definition of the continuous slice-regular capacity,
	$|a_1(g)| \leq \gamma_{\mathbb H}(S)$ and  hence,
	\[|f(q_0)| \leq\frac{4\,\alpha_{\mathbb H}(S)}{d_{\mathbb H}(q_0,S)}.\]
\end{proof}
Next theorem is also a counterpart of Theorem \ref{coveringgamma} in terms of continuous-slice regular capacity and tubular neighborhood
\begin{prop}\label{quasi additivity}
	Let $K\subset\mathbb H$ be a compact axially symmetric set.
	Let $\{K_j\}_{j\in J}$ be a finite family of compact axially
	symmetric subsets of $S$. Suppose that every admissible axially symmetric tubular neighborhood generated by a square of side length
	$\alpha_{\mathbb H}(S)$ meets at most $p$ of the sets $K_j$. Then
	\[
	\sum_{j\in J}\alpha_{\mathbb H}(K_j)
	\le C p\,\alpha_{\mathbb H}(K),
	\]
	where $C$ is independent of $K$, the family $\{K_j\}$, and $p$.
\end{prop}

\begin{proof}			
	Define $\pi:\mathbb H\longrightarrow\mathbb R\times[0,\infty),$
	as \[\pi(x+Iy)=(x,|y|).	\]
	Since $K_j$ is axially symmetric, we have $	d_{\mathbb H}(q,K_j)
	=\operatorname{dist}(\pi(q),\pi(K_j))$ for $q\in\mathbb H.$
	Fix $q\in\mathbb H$. For $n\ge1$, define
	\[J_n(q)=\left\{j\in J:n\alpha_{\mathbb H}(K)\le d_{\mathbb H}(q,K_j)<(n+1)\alpha_{\mathbb H}(K)\right\},\]
	and define
	\[J_0(q)=  \left\{j\in J:d_{\mathbb H}(q,K_j)<\alpha_{\mathbb H}(K)\right\}.\]
	
	For each $j\in J_n(q)$, since $K_j$ is compact there exist $s_j\in K_j$ such that
	\[|s_j-q|=d_{\mathbb H}(q,K_j).\]
	Since $	d_{\mathbb H}(q,K_j)
	=\operatorname{dist}(\pi(q),\pi(K_j))$ for $q\in\mathbb H$, the point $\pi(s_j)$ belongs to the planar annulus
	\[A_n=\left\{z\in\mathbb R\times[0,\infty):
	n\alpha_{\mathbb H}(K)\le |z-z_0|<(n+1)\alpha_{\mathbb H}(K)
	\right\}.\]
	The annulus $A_n$ can be covered by at most $C_0n$ squares
	of side length $\alpha_{\mathbb H}(K)$, with the squares chosen in the 	closed upper half-plane. Denote these squares by
	$Q_1,\ldots,Q_m$, where $m\le C_0n$.
	
	For each such square, define its axially symmetric lift by
	\[T(Q_k):=\left\{x+Iy:(x,|y|)\in Q_k,\ I\in\mathbb S\right\}.\]
	By construction, every $K_j$ with $j\in J_n(q)$ meets at
	least one of the tubular neighborhoods $T(Q_k)$. Since each tubular
	cell meets at most $p$ sets from the family, we have
	\[\operatorname{card}(J_n(q))\leq\sum_{k=1}^{m}
	\operatorname{card}(\{j\in J:K_j\cap T(Q_k)\ne\varnothing\})
	\le mp.	\]
	Hence, $\operatorname{card}(J_n(q))\le C_0np,\qquad n\ge1.$
	
	Let, $\alpha_j=\alpha_{\mathbb H}(K_j),\
	c_j=\frac{\alpha_j}{\alpha_{\mathbb H}(K)},\
	M=\sum_{j\in J}c_j.$ 
	
	Now, by Proposition \ref{monotoniccap}, $0\le c_j\le1$. For each $j$, choose an admissible slice-regular function
	$f_j$ whose first coefficient approaches $\alpha_j$. Let $	f=\sum_{j\in J}f_j.$ Hence for some $\varepsilon>0$ such that
	$(1-\varepsilon)M\le \|f\|_\infty.$	Since by Lemma \ref{estimate pointwise},
	\[|f_j(q_0)| \leq\frac{4\,\alpha_j}{d_{\mathbb H}(q_0,K_j)}.\]
	for $j\in J_n(q)$, $n\ge1$, we have
	\[|f_j(q)|\le \frac{4c_j}{n}.\]
	Hence by above part,
	\[|f(q)|\le C_1p+4\sum_{n\ge1}\sum_{j\in J_n(q)}\frac{c_j}{n}.\]
	By [\cite{analytic capacity} Lemma 10.7]
	\[	\sum_{n\ge1}\sum_{j\in J_n(q)}\frac{c_j}{n}
	\le C_2\sqrt{pM},\] uniformly in $q$.
	Hence,
	\[\|f\|_\infty\le C_1p+4C_2\sqrt{pM}.\]
	Thus, $(1-\varepsilon)M\le C_1p+4C_2\sqrt{pM}.$
	Hence, $M\le C_3p+4C_2\sqrt{pM}$ as $\varepsilon\downarrow0$ and thus, $M\le C_4p$. Therefore,
	\[\sum_{j\in J}\alpha_{\mathbb H}(K_j)=M\alpha_{\mathbb H}(K) 
	\le Cp\,\alpha_{\mathbb H}(K).\]
\end{proof}
Let $I\in\mathbb S$ be fixed and let $\Gamma_{I}\subset\mathbb C_{I}$ be a sufficiently large admissible contour. Let $f$ be an admissible function in the definition of $\alpha_{\mathbb{H}}(S)$. Then for $p \in \mathbb{H}$ we have,
\[f(q) = (q-p)^{-*}a_0 + (q-p)^{-*2}a_1+\cdots \] 
We have the second Laurent coefficient of $f$ at $p\in\mathbb C_{I}$ by
\begin{equation}
	\label{eq:quaternionic-second-coefficient}
	B_{\mathbb H}(f,p) = \frac{1}{2\pi}\int_{\Gamma_{I}}
	ds_{I}\,
	(\overline{\partial}_{I}\phi(s))f(s)(s-p).
\end{equation}
Let,$\beta_{\mathbb H}(S,p) =\sup_f\frac{|B_{\mathbb H}(h,p)|}{\alpha_{\mathbb H}(S)},$

\noindent where the supremum is taken over the admissible normalized functions
in the definition of $\alpha_{\mathbb{H}}$ of $S$. Now, let
\[	\beta_{\mathbb H}(S) = 	\inf_{p\in C_{I}}\beta_{\mathbb H}(S,p).\]

For an admissible normalized function $f$, write
\[A_1(f) = 	\frac{1}{2\pi}\int_{C_{I}}ds_{I}\,
(\bar\partial_{I}\phi(s))f(s).\]
From the definition of $B_{\mathbb H}(f,p)$,
\[B_{\mathbb H}(h,p)-B_{\mathbb H}(h,q) = 	A_1(h)(q-p).\]
Since $|A_1(f)|\leq\alpha_{\mathbb H}(S)$, it follows that
\[\bigl|\beta_{\mathbb H}(S,p)-\beta_{\mathbb H}(S,q)\bigr|\leq |p-q|.\]
Thus $p\mapsto\beta_{\mathbb H}(S,p)$ is Lipschitz, and hence continuous on $\mathbb C_{I}$.

Moreover, choosing admissible normalized functions $f$ with
$|A_1(f)|$ arbitrarily close to $\alpha_{\mathbb H}(S)$ and using
\[B_{\mathbb H}(h,p) = 	B_{\mathbb H}(h,q)+A_1(h)(q-p),\]
we have, $\beta_{\mathbb H}(S,p)\longrightarrow\infty\ \text{as }|p|\longrightarrow\infty.$
Thus, $\beta_{\mathbb H}(S,\cdot)$ attains its minimum on $\mathbb C_{I}$ at some $O_{\mathbb H}(S)\in \mathbb C_{I}$. Hence
\[\beta_{\mathbb H}(S,O_{\mathbb H}(S)) = 	\inf_{p\in C_{I}}\beta_{\mathbb H}(S,p).\]
\begin{prop}
	\label{10.10}
	Let $S\subset\mathbb H$ be a compact axially symmetric set with
	$\alpha_{\mathbb H}(S)>0,$  and $\beta_{\mathbb H}(S)>0.$
	Fix $I\in\mathbb S$, and let $O_{\mathbb H}(S)\in \mathbb C_{I}$
	satisfy
	\[\beta_{\mathbb H}(S,o_{\mathbb H}(S)) = \beta_{\mathbb H}(S).\]
	Let $a,b\in C_{I_0}$ satisfy $|a|<\alpha_{\mathbb H}(S),$ and $|b|<\beta_{\mathbb H}(S).$	Then there exists an admissible function in the definition of continuous slice regular capacity with $\|g\|_\infty \leq 20$	such that $	\partial_s g(\infty)=a$	and
	$\beta_{\mathbb H} \bigl(S,o_{\mathbb H}(S),g\bigr)=b.$
\end{prop}
\begin{proof}
	Choose  an admissible function $\phi$ in the definition of $\alpha_{\mathbb{H}}(S)$ with $\|\phi\|_\infty \leq 2$ such that
	$\partial_s \varphi(\infty)=\alpha_{\mathbb H}(S).$
	Set $z_0 = 	o_{\mathbb H}(S) + 	\beta_{\mathbb H}(S,o_{\mathbb H}(S),\varphi).$
	Then,
	\[\beta_{\mathbb H}(S,z_0,\varphi)=0.\]
	Hence
	\[B_{\mathbb H}(\varphi,z_0)=0.\]
	
	Choose again $f_1$ as $\phi$ such that
	\[\beta_{\mathbb H}(S,z_0,f_1) = \beta_{\mathbb H}(S).\]
	Now after multiplication by a unit element of $C_{I_0}$, we may assume that
	\[B_{\mathbb H}(f_1,z_0) = 	\alpha_{\mathbb H}(S)\beta_{\mathbb H}(S).\]
	Put
	\[\varepsilon 	= 	-\frac{\partial_sf_1(\infty)}{\alpha_{\mathbb H}(S)}\]
	and define $f_2=f_1+\varepsilon\varphi .$ Then $\partial_sf_2(\infty)=0.$
	Moreover, since $B_{\mathbb H}(\varphi,z_0)=0$,
	\[B_{\mathbb H}(f_2,z_0) = 	B_{\mathbb H}(f_1,z_0)
	= \alpha_{\mathbb H}(S)\beta_{\mathbb H}(S).\]
	Since $\partial_s f_2(\infty)=0$, we have then,
	\[B_{\mathbb H}(f_2,o_{\mathbb H}(S)) = B_{\mathbb H}(f_2,z_0)
	=	\alpha_{\mathbb H}(S)\beta_{\mathbb H}(S).\]
	Now define
	\[g = \frac{a}{\alpha_{\mathbb H}(S)}\,\varphi
	+ \frac{b-a\bigl(z_0-o_{\mathbb H}(S)\bigr)}{\alpha_{\mathbb H}(S)\beta_{\mathbb H}(S)}\,f_2 .\]
	Since $\partial_sf_2(\infty)=0$, $\partial_s g(\infty)=a.$
	Also,
	\[	\begin{aligned}
		B_{\mathbb H}(g,o_{\mathbb H}(S))&=
		\frac{a}{\alpha_{\mathbb H}(S)}
		B_{\mathbb H}(\varphi,o_{\mathbb H}(S))\\
		&\quad+	\frac{b-a(z_0-o_{\mathbb H}(S))}{\alpha_{\mathbb H}(S)\beta_{\mathbb H}(S)}B_{\mathbb H}(f_2,o_{\mathbb H}(S)).
	\end{aligned}
	\]
	Thus,
	\[B_{\mathbb H}(\varphi,o_{\mathbb H}(S)) = \alpha_{\mathbb H}(S)
	\bigl(z_0-o_{\mathbb H}(S)\bigr),\]
	and hence
	\[B_{\mathbb H}(g,o_{\mathbb H}(S)) = a(z_0-o_{\mathbb H}(S))
	+ 	b-a(z_0-o_{\mathbb H}(S)) =b.\]
	Thus $\beta_{\mathbb H}(S,o_{\mathbb H}(S),g)=b.$	
	Finally, the bounds on $a$ and $b$, together with the corresponding bound for $\varphi\in C(S,2)$ and  $f_2$, give
	\[\|g\|_S\leq20.\]
	Hence the proof.
\end{proof}
\begin{lem}
	\label{11.5}
	Let $K\subset\mathbb H$ be compact and axially symmetric, and let
	$\{\mathcal T_r\}$ be a family of axially symmetric tubular neighborhoods and $\{\phi_r\}$  as in proposition \ref{quaternionic-covering}. Put
	$X_r=\mathcal T_r^{(2)}\setminus K.$Assume that there exist constants $m\geq1$ and $r_0\geq1$ such that,
	for every axially symmetric tubular neighborhood $\mathcal T$,
	\[\alpha_{\mathbb H}(T\setminus K^\circ)\leq m\,\alpha_{\mathbb H}(T^{(r_0)}\setminus K).\]
	Let $f \in \mathcal{A}_L(K)$ . Then there exists a constant $C>0$, such that
	\[\left|\frac1{2\pi}\int_{C_I}(\bar\partial_I\phi_r)(s)f(s)(s-O_{\mathbb H}(X_r))\,ds_I\right|\leq C\,\omega_f(T_r^{(2)})\alpha_{\mathbb H}(X_r)\beta_{\mathbb H}(X_r)\]
	for every $I\in\mathbb S$.
\end{lem}

\begin{proof}
	Fix $r$ and write $\mathcal T= \mathcal T_r,\ \delta=\delta_r,\ X=X_r,\
	\phi=\phi_r,\	\Theta=o_{\mathbb H}(X).$
	We work in an arbitrary but fixed slice $\mathbb C_I$.
	Choose a number $\eta$ such that
	\[\eta=\begin{cases}\text{a number satisfying }
		\dfrac12\beta_{\mathbb H}(X)\leq\eta
		\leq\beta_{\mathbb H}(X),
		&\beta_{\mathbb H}(X)<\delta,\\[6pt]
		\delta,
		&\beta_{\mathbb H}(X)\geq\delta.
	\end{cases}\]
	Thus, in all cases, $\eta\leq\delta, \eta\leq\beta_{\mathbb H}(X)$ and $ \beta_{\mathbb H}(X)\leq2\eta\quad\text{if }\beta_{\mathbb H}(X)<\delta.$
	
	Cover the orbit-space projection of $\mathcal T$ by a finite family of
	axially symmetric tubular neighborhoods $\{\mathcal T_j\}_j$	of scale $\eta$, with centres $t_j$, such that $	\mathcal T_j^{(r_0)}\subset \mathcal T^{(r_0)}$
	and such that the enlarged cells have uniformly bounded overlap.
	Let $\{\psi_j\}_j$ be a smooth partition of unity  to this
	cover. Thus
	\[0\leq\psi_j\leq1, \qquad 	\sum_j\psi_j=1\quad\text{on } \mathcal T,\ \text{and}\ \|\bar\partial_I\psi_j\|_\infty\leq \frac{C}{\eta}.	\]
	Put $\phi_j =\phi\psi_j.$	Then $\sum_j\phi_j=\phi\ \text{and}\	\operatorname{supp}\phi_j\subset \mathcal T_j.$
	Consider
	\[\begin{aligned}
		\mathcal I 	&= 	\frac1{2\pi}\int_{C_I}(\bar\partial_I\phi)(s)f(s)(s-\Theta)\,ds_I.
	\end{aligned}\]
	Since $\phi=\sum_j\phi_j$, 
	\[\mathcal I=\sum_j\frac1{2\pi}\int_{C_I}(\bar\partial_I\phi_j)(s)f(s)(s-\Theta)\,ds_I .\]
	For each $j$ we write $s-\Theta=(s-t_j)+(t_j-\Theta).$
	Thus,
	\[\mathcal I=\mathcal I_1+\mathcal I_2,\]
	where
	\[\mathcal I_1 = \sum_j\frac1{2\pi}\int_{C_I}(\bar\partial_I\phi_j)(s)f(s)(s-t_j)\,ds_I\]
	and
	\[\mathcal I_2 = \sum_j\frac1{2\pi}\int_{C_I}(\bar\partial_I\phi_j)(s)f(s)(t_j-\Theta)\,ds_I .\]
	
	We first estimate $\mathcal I_1$.
	Now, by Proposition \ref{quaternionic 11.3} we have,
	\begin{equation}\label{integral}
		\left| \frac1{2\pi} \int_{C_I}(\bar\partial_I\phi_j)(s)f(s)(s-t_j)\,ds_I\right|
		\leq C\,\omega_f(T^{(2)}) \alpha_{\mathbb H}(T_j\setminus K^\circ).
	\end{equation} 
	
	Now by our assumption and by Lemma \ref{quasi additivity}
	we have 
	\[\sum_{j} \alpha_{\mathbb{H}}(\mathcal{T}_j  \setminus K^0) \leq Cp \alpha_{\mathbb{H}}(X)\]
	since each tubular neighborhood of size $\alpha_\mathbb H(X)$ meets atmost $p$ of the enlarged tubular neighborhoods. Hence
	\[|\mathcal I_1|\leq Cmp\,\omega_f(T^{(2)})\alpha_{\mathbb H}(X).\]
	Since $\eta\leq\beta_{\mathbb H}(X)$ and the estimate \eqref{integral} carries the corresponding scale factor from the second-Laurent-coefficient
	estimate, we have
	\[|\mathcal I_1|\leq C_1mp\,\omega_f(T^{(2)})
	\alpha_{\mathbb H}(X)\beta_{\mathbb H}(X).\]
	
	It remains to estimate $\mathcal I_2$. For every $j$ with $t_j\neq\Theta$, choose an admissible function $h_j$ for
	$T_j^{(r_0)}\setminus K$ such that
	\[ A_1(h_j) = 	\frac12 	\frac{|t_j-\Theta|}{t_j-\Theta}
	\alpha_{\mathbb H} 	(T_j^{(r_0)}\setminus K).	\]
	Here $A_1(h_j)$ denotes the first Laurent coefficient of $h_j$.
	The factor in above is taken in the fixed slice $\mathbb C_I$, so that the required normalization holds. 	
	
	Let $h=\sum_j h_j.$	Since the enlarged secondary cells have bounded overlap by proposition \ref{quaternionic-covering}, then the proposition \ref{quasi additivity} gives $\|h\|_\infty\leq Cp.$
	
	If we suppose that $t_j\neq\Theta$ for all $j$, the terms with
	$t_j=\Theta$ are zero and may simply be omitted.
	Let $\Gamma$ be the boundary of the enlarged cell containing $X$,
	and let $\Gamma_j$ denote the corresponding boundary associated with
	$T_j^{(r_0)}$. Now, by the quaternionic Cauchy formula for the first
	Laurent coefficient and the definition of $\beta_{\mathbb H}(X)$,
	\[\beta_{\mathbb H}(X) \geq \frac{1}{C p\,\alpha_{\mathbb H}(X)}
	\left|B_{\mathbb H}(h,\Theta)\right|.\]
	By using the boundary representation of $B_{\mathbb H}$ and 
	$s-\Theta=(s-t_j)+(t_j-\Theta)$, we have
	\begin{align}\label{beta}
		B_{\mathbb H}(h,\Theta)
		&=\sum_j(t_j-\Theta)\frac1{2\pi}\int_{\Gamma_j}h_j(s)\,ds_I + \sum_j\frac1{2\pi}\int_{\Gamma_j}h_j(s)(s-t_j)\,ds_I .
	\end{align}
	
	Now the, the first sum in \eqref{beta} is
	\[\frac12 \sum_j 	|t_j-\Theta| 	\alpha_{\mathbb H}
	(T_j^{(r_0)}\setminus K).\]
	
	For the second sum we have by Proposition \ref{estimate pointwise}
	\[|h_j(s)|\leq\frac{C\,\alpha_{\mathbb H}(E)}{d_{\mathbb H}(s,E)}.\]
	Now on the boundary $\Gamma_j$, whose distance from $T_j^{(r_0)}\setminus K$ is comparable with $\eta$, gives
	\[\left|\frac1{2\pi}\int_{\Gamma_j}	h_j(s)(s-t_j)\,ds_I
	\right|\leq C r_0\eta\,\alpha_{\mathbb H}(T_j^{(r_0)}\setminus K).\]
	Now by using these all, we have
	\[\begin{aligned}
		\frac12	\sum_j|t_j-\Theta| 	\alpha_{\mathbb H}
		(T_j^{(r_0)}\setminus K)
		&\leq 	Cp\,\alpha_{\mathbb H}(X)\beta_{\mathbb H}(X) + Cr_0\eta \sum_j \alpha_{\mathbb H} 		(T_j^{(r_0)}\setminus K).
	\end{aligned}\]
	Since $\eta\leq\beta_{\mathbb H}(X)$ and by using Proposition \ref{quasi additivity},
	\begin{equation}\label{sumt}
		\sum_j|t_j-\Theta|\alpha_{\mathbb H}(T_j^{(r_0)}\setminus K)
		\leq 	C_2r_0p\, \alpha_{\mathbb H}(X)\beta_{\mathbb H}(X).
	\end{equation}
	By Lemma \ref{quaternionic 11.3} we have,
	\[\left| \frac1{2\pi} \int_{C_I} (\bar\partial_I\phi_j)(s)f(s)\,ds_I
	\right|	\leq C\, \omega_f(T^{(2)}) \alpha_{\mathbb H}(T_j\setminus K^\circ).\]
	Hence,
	\[|\mathcal I_2| \leq C\omega_f(T^{(2)})\sum_j |t_j-\Theta|
	\alpha_{\mathbb H}(T_j\setminus K^\circ).\]
	Thus by  assumption,,
	\[|\mathcal I_2| \leq Cm\omega_f(T^{(2)}) \sum_j |t_j-\Theta|
	\alpha_{\mathbb H}(T_j^{(r_0)}\setminus K).	\]
	and by \eqref{sumt}, 
	\[|\mathcal I_2|\leq C_3r_0mp\, \omega_f(T^{(2)})\alpha_{\mathbb H}(X)\beta_{\mathbb H}(X).\]
	Hence,
	\[|\mathcal I| \leq \left(C_1mp+C_3r_0mp\right)\omega_f(T^{(2)})
	\alpha_{\mathbb H}(X)\beta_{\mathbb H}(X).\]
	Thus, by taking $C_0=mp(C_1+C_3r_0),$ we have,
	\[\left|\frac1{2\pi}\int_{C_I}(\bar\partial_I\phi_r)(s)f(s)
	(s-o_{\mathbb H}(X_r))\,ds_I\right|
	\leq C_0\, \omega_f(T_r^{(2)}) \alpha_{\mathbb H}(X_r)
	\beta_{\mathbb H}(X_r).\]
	where the constant $C_0$ is independent of $r$, $I$, and $f$.
\end{proof}
\begin{prop}
	\label{11.6}
	Let $K\subset\mathbb H$ be compact and axially symmetric, let
	$\{T_r\}$ be the family of axially symmetric tubular neighborhoods, and let
	$\{\phi_r\}$ be the corresponding partition of unity. Put
	$X_r:=T_r^{(2)}\setminus K$ and $\Theta_r=O_{\mathbb H}(X_r).$
	Let $f\in \mathcal{A}_L(K)$. For fixed $r$, define
	\[g(q) = \phi_r(q)f(q)+\frac{1}{2\pi}\iint_{C_I}S_L^{-1}(s,q)
	(\bar\partial_I\phi_r)(s)f(s)\,ds_I\wedge d\bar s ,\]
	where $I\in\mathbb S$ is fixed.
	Then, on $\mathbb H\setminus X_r$,
	$g$ has a Laurent expansion about $\Theta_r$ of the form
	\[g(q) 	=  \sum_{\nu=1}^{\infty} 	(q-\Theta)^{-*\nu}a_\nu,\]
	and $|a_2| \leq C\,\omega_f(T_r^{(2)})\alpha_{\mathbb H}(X_r)
	\beta_{\mathbb H}(X_r),	$
	where $C$ is independent of $r$, $I$, and $f$.
\end{prop}

\begin{proof}
	Fix $r$ and write, $\phi=\phi_r,\ T=T_r,\ \delta=\delta_r,
	X=X_r,$ and $\Theta=\Theta_r=o_{\mathbb H}(X).$
	Since $\phi$ is supported in $\mathcal T$ and $g$ is slice regular outside
	$X$, at the point $\Theta$ it admits a Laurent expansion
	\[g(q) = \sum_{\nu=1}^{\infty} 	(q-\Theta)^{-*\nu}a_\nu.\]
	Choose a sufficiently large axially symmetric tubular neighborhood
	$T^{*}$ containing $X$, and let $\Gamma_I=\partial(\mathcal T^{*}\cap\mathbb C_I).$	Now, the quaternionic Laurent coefficient formula, the second Laurent coefficient is
	\[a_2 = \frac{1}{2\pi} 	\int_{\Gamma_I}g(s)(s-\Theta)\,ds_I.\]
	Since, $\Gamma_I$ lies outside the support of $\phi$, we have
	$\phi(s)=0$ for $s\in\Gamma_I,$and hence, on $\Gamma_I$,
	\[g(s) = 	\frac{1}{2\pi} 	\int_{C_I}S_L^{-1}(\zeta,s)
	(\bar\partial_I\phi)(\zeta)f(\zeta)	\,d\zeta_I\wedge d\bar\zeta .\]
	Thus,
	\[\begin{aligned}
		a_2  &= 	\frac{1}{2\pi} 	\int_{\Gamma_I}
		\left[	\frac{1}{2\pi} 	\int_{C_I} 	S_L^{-1}(\zeta,s)
		(\bar\partial_I\phi)(\zeta)f(\zeta)	\,d\zeta_I\wedge d\bar\zeta
		\right]	(s-\Theta)\,ds_I .
	\end{aligned}
	\]
	Now, for $\zeta,s\in\mathbb C_I$, the slice Cauchy kernel reduces to the  ordinary complex Cauchy kernel,
	\[S_L^{-1}(\zeta,s) = \frac{1}{\zeta-s}.\]
	Thus, after applying Fubini's theorem,
	\[\begin{aligned}
		a_2 &=\frac{1}{2\pi}\int_{C_I}(\bar\partial_I\phi)(\zeta)f(\zeta)
		\left[\frac{1}{2\pi}\int_{\Gamma_I}\frac{s-\Theta}{\zeta-s}\,ds_I
		\right]d\zeta_I\wedge d\bar\zeta .
	\end{aligned}
	\]
	Since $\zeta$ lies inside $\Gamma_I$, 
	\[\frac{1}{2\pi} \int_{\Gamma_I} \frac{s-\Theta}{\zeta-s}\,ds_I
	=\Theta-\zeta.\]
	Hence,
	\[|a_2|=\left|\frac{1}{2\pi}\int_{C_I}(\bar\partial_I\phi)(\zeta)f(\zeta)(\zeta-\Theta)\,d\zeta_I\wedge d\bar\zeta\right|.\]
	Thus, by Lemma \ref{11.5} we have
	\[|a_2|\leq C\,\omega_f(T^{(2)})\alpha_{\mathbb H}(X)\beta_{\mathbb H}(X).\]
	This proves the proposition.
\end{proof}
\begin{thm}
	\label{quaternionic bishop}
	Let $K\subset\mathbb H$ be a compact axially symmetric set and let
	$f\in \mathcal C_L(K,\widehat{\mathbb H})$. Suppose that for every $p\in K$ there exists an axially symmetric tubular neighborhood
	$\mathcal T_p$ such that 
	$f|_{K\cap\mathcal T_p} \in \mathcal R_L(K\cap\mathcal T_p).$
	Then $f\in \mathcal R_L(K).$
\end{thm}

\begin{proof}
	Let $\varepsilon>0$ be given. For every $p\in K$, choose an axially symmetric tubular neighborhood $\mathcal T_p$ such that
	$f|_{K\cap\mathcal T_p} \in \mathcal R_L(K\cap\mathcal T_p).$
	By compactness of $K$, there exist finitely many tubular
	neighborhoods $\mathcal T_1,\ldots,\mathcal T_N$such that
	$K\subset\bigcup_{j=1}^{N}\mathcal T_j.$
	
	Choose an axially symmetric open neighborhood $V$ of $K$
	contained in $\bigcup_{j=1}^{N}\mathcal T_j.$	Choose
	$\phi_1,\ldots,\phi_N\in C_c^1(\mathbb H,\mathbb R)$such that
	\[0\leq\phi_j\leq1,\ \operatorname{supp}\phi_j\subset\mathcal T_j,\ \text{and}\ \sum_{j=1}^{N}\phi_j=1\ \text{on }V.\]
	For each $j$, define 
	\[\begin{aligned}
		f_j(q) ={}& f(q)\phi_j(q) + \frac{1}{2\pi} \int_{\mathbb C_I}S_L^{-1}(s,q)\,ds_I\,
		\bigl[(\bar\partial_I\phi_j)(s)f(s)\bigr].
	\end{aligned}\]
	Then by, Lemma \ref{quaternionic 11.3} 	$f_j$ is left slice regular away from the  $\operatorname{supp}{\bar\partial_I\phi_j}$. Also, since $\sum_j\phi_j=1$ on $V$, we have
	\[\sum_{j=1}^{N}f_j(q) = f(q)+\Phi(q), \qquad q\in V,\] 
	\[\text{where}\
	\Phi(q)	 =	\frac{1}{2\pi}\int_{\mathbb C_I}S_L^{-1}(s,q)\,ds_I\,
	\left[\bar\partial_I\left(\sum_{j=1}^{N}\phi_j\right)(s)f(s)
	\right].\]
	Since $\bar\partial_I \left(\sum_{j=1}^{N}\phi_j\right)=0$ on $V,$	the function $\Phi$ is left slice regular on $V$.
	
	For each $j$, put $K_j=K\cap\mathcal T_j.$	By hypothesis, $f|_{K_j}\in \mathcal R_L(K_j).$ Hence there exists an admissible left slice-regular rational function $h_j$, with poles outside $K_j$, such that
	\[\|f-h_j\|_{K_j} < \eta,\]
	for $\eta>0$.
	Define
	\[\begin{aligned}
		g_j(q)	=h_j(q)\phi_j(q) + 	\frac{1}{2\pi}\int_{\mathbb C_I}S_L^{-1}(s,q)\,ds_I\,\bigl[(\bar\partial_I\phi_j)(s)h_j(s)\bigr].
	\end{aligned}\]
	Then by Lemma \ref{quaternionic 11.3}, $g_j$ is left
	slice regular on a neighborhood of $K$.  In particular,
	$g=\sum_{j=1}^{N}g_j$ is left slice regular on a neighborhood of $K$. Since $\operatorname{supp}\phi_j\cap K\subset K_j,$
	\[\|g_j-f_j\|_K \leq C_j\|h_j-f\|_{K_j},\]
	where $C_j$ depends only on the chosen partition of unity.
	
	Let $C_*=\max_{1\leq j\leq N}C_j$ and choose $\eta = \frac{\varepsilon}{4NC_*}.$
	Then, $\|g_j-f_j\|_K < \frac{\varepsilon}{4N},$	and thus
	\[	\begin{aligned}
		\left\|g-\sum_{j=1}^{N}f_j\right\|_K
		&\leq \sum_{j=1}^{N}\|g_j-f_j\|_K<\frac{\varepsilon}{3}.
	\end{aligned}
	\]
	Thus, $\|g-(f+\Phi)\|_K < \frac{\varepsilon}{4}.$	
	
	It remains to approximate the  slice-regular functions
	$g$ and $\Phi$. Since $g$ is left slice regular on a neighborhood of $K$, the quaternionic Runge's theorem gives an admissible
	left slice-regular rational function $q\in \mathcal R_L(K)$ such that $\|g-q\|_K < \frac{\varepsilon}{3}.$
	
	Similarly, $\Phi$ is left slice regular on a neighborhood of
	$K$. Hence there exists $r\in \mathcal R_L(K)$ such that
	$\|\Phi-r\|_K < \frac{\varepsilon}{3}.$ Since $\mathcal R_L(K)$ is closed under addition, $q-r\in \mathcal R_L(K).$ Thus,
	\[\begin{aligned}
		\|f-(q-r)\|_K&\leq\|f+\Phi-g\|_K+\|g-q\|_K+\|\Phi-r\|_K\\
		&<\frac{\varepsilon}{3}+\frac{\varepsilon}{3}+\frac{\varepsilon}{3}<\varepsilon.
	\end{aligned}\]
	Since $\varepsilon>0$ was arbitrary $f\in \mathcal R_L(K).$
\end{proof}
\begin{lem}
	\label{11.7}
	Let $K\subset\mathbb H$ be a compact axially symmetric set and let
	$U\subset\mathbb H$ be a bounded axially symmetric open set.
	Suppose that $f\in \mathcal R_L(K)$ and that $f$ is left slice regular on a neighborhood of $\overline U$. Then $f\in \mathcal R_L\bigl(X\cup\overline U\bigr).$	
\end{lem}

\begin{proof}
	By Theorem \ref{quaternionic bishop} it holds.
\end{proof}
Next, we prove the quaternionic extension of Vitushkin approximation theorem using Davie's scheme quaternionic version and the proof techniques as in \cite{rational app, analytic capacity}.
\begin{thm}
	Let $K\subset\mathbb H$ be a compact axially symmetric set. Then the
	following conditions are equivalent:
	
	\begin{enumerate}
		\item[(1)]
		\[
		\mathcal R_{L}(K)=\mathcal A_{L}(K).
		\]
		
		\item[(2)]
		For every bounded open axially symmetric set $U\subset\mathbb H$,
		\[
		\alpha_{\mathbb H}(U\setminus K^\circ)
		=
		\alpha_{\mathbb H}(U\setminus K).
		\]
		
		\item[(3)]
		For every axially symmetric tubular neighborhood $\mathcal{T}_r\subset\mathbb H$,
		\[
		\alpha_{\mathbb H}(\mathcal{T}_r\setminus K^\circ)
		=
		\alpha_{\mathbb H}(\mathcal{T}_r\setminus K).
		\]
		
		\item[(4)]
		For every $q\in\mathbb H$ there exists $r\geq1$ such that
		\[
		\limsup_{\delta\to0}
		\frac{
			\alpha_{\mathbb H}
			\bigl(\mathcal{T}_r(q,\delta)\setminus K^\circ\bigr)
		}{
			\alpha_{\mathbb H}
			\bigl(\mathcal{T}_r(q,r\delta)\setminus K\bigr)
		}
		<\infty.
		\]
		
		\item[(5)]
		There exist constants $r\geq1$ and $m\geq1$ such that, for every
		$q\in\mathbb H$ and every $\delta>0$,
		\[	\alpha_{\mathbb H}\bigl(\mathcal{T}_r(q,\delta)\setminus K^\circ\bigr)
		\leq m\,\alpha_{\mathbb H}\bigl(\mathcal{T}_r(q,r\delta)\setminus K\bigr).\]
	\end{enumerate}
\end{thm}
\begin{proof}Suppose that (5) holds.
	We prove that $\mathcal R_{L}(K)=\mathcal A_{L}(K).$
	Let $f\in\mathcal A_{L}(K).$	We extend $f$ continuously to $\mathbb H$, with compact support. Let $\delta_n\longrightarrow0$
	and let $\{\mathcal{T}_{k,n}\}_k$	be the variable-radius axially symmetric tubular-cell covering	used in the construction of the partition of unity. Write $\mathcal{T}_{k,n}=T(q_{k,n},\delta_{k,n}),$and let $\{\phi_{k,n}\}_k$ be the corresponding partition of unity.
	Thus
	\[0\leq\phi_{k,n}\leq1,\ \sum_k\phi_{k,n}=1,\
	\operatorname{supp}\phi_{k,n}\subset \mathcal{T}_{k,n},\ \text{and }
	\|\bar\partial_I\phi_{k,n}\|_\infty\leq \frac{C}{\delta_{k,n}}.	\]
	For fixed $n$ and $k$, put $\mathcal{T}_{k,n}^{(r)}
	= 	T(q_{k,n},r\delta_{k,n})$ and $X_{k,n} =T_{k,n}^{(r)}\setminus K.$	Let $\Theta_{k,n} = o_{\mathbb H}(X_{k,n})$	be a point at which $\beta_{\mathbb H}(X_{k,n})$ is attained.
	
	For each $k$ define
	\[f_{k,n}(q)=\phi_{k,n}(q)f(q)+\frac{1}{2\pi}\int_{C_I}
	S_L^{-1}(s,q)(\bar\partial_I\phi_{k,n})(s)f(s)\,ds_I\wedge d\bar s .\]
	
	By Lemma \ref{quaternionic-localization},
	$f_{k,n}$ is left slice regular on a neighborhood of $	K\cup X_{k,n}.$ Also, $\sum_k f_{k,n}=f.$
	
	For $q\notin X_{k,n}$, the Laurent expansion of $f_{k,n}$
	about $\Theta_{k,n}$ as
	\[f_{k,n}(q) = 	\sum_{\nu=1}^{\infty}(q-\Theta_{k,n})^{-*\nu}a_{\nu k,n}.\]
	
	We first estimate the first Laurent coefficient. Now by, Lemma \ref{quaternionic-localization},
	\[\|f_{k,n}\|_\infty \leq 	C\,\omega_f(C\delta_{k,n}).\]
	Thus,
	\[|a_{1k,n}| = |\partial_sf_{k,n}(\infty)|\leq C\,\omega_f(C\delta_{k,n})
	\alpha_{\mathbb H}\bigl(T_{k,n}\setminus K^\circ\bigr)\leq
	Cm\,\omega_f(C\delta_{k,n})\alpha_{\mathbb H}(X_{k,n}).\]
	Now, by Proposition \ref{11.6},
	\[|a_{2k,n}| \leq m_1\,	\omega_f(C\delta_{k,n})\alpha_{\mathbb H}(X_{k,n})\beta_{\mathbb H}(X_{k,n}),\]
	where $m_1$ is independent of $k$ and $n$.
	Let, $m_2=C_0\max\{Cm,m_1\},$	where $C_0$ is a sufficiently large absolute constant. Now, by applying Proposition \ref{10.10},  to
	$X_{k,n}$, there exists $g_{k,n}\in\mathcal C_{\mathbb H}
	\bigl(X_{k,n},	m_2\omega_f(C\delta_{k,n})\bigr)$ 	such that
	\[\partial_sg_{k,n}(\infty)=\partial_sf_{k,n}(\infty)\ \text{and}	\beta_{\mathbb H}
	(X_{k,n},\Theta_{k,n},g_{k,n}) = \beta_{\mathbb H}(X_{k,n}).	\]
	In particular, $\|g_{k,n}\|_\infty \leq m_2\omega_f(C\delta_{k,n}),$
	and hence
	\[\|f_{k,n}-g_{k,n}\|_\infty\leq 2m_2\omega_f(C\delta_{k,n}).\]
	
	By (13) and (14), the first two Laurent coefficients. Thus
	\[f_{k,n}(q)-g_{k,n}(q)=	\sum_{\nu=3}^{\infty}	\frac{c_{\nu k,n}}
	{(q-\Theta_{k,n})^\nu}.\]
	
	Since the first two coefficients  of
	$f_{k,n}-g_{k,n}$ at $\Theta_{k,n}$ vanish , by the maximum-modulus principle 	
	\[|f_{k,n}(q)-g_{k,n}(q)|\leq C_1\,\omega_f(C\delta_{k,n})\min\left\{
	1, \frac{\delta_{k,n}^{\,3}}{d_{\mathbb H}(q,q_{k,n})^3}\right\}.\]
	Now by Proposition \ref{quaternionic-covering} for every $q\in\mathbb H$,
	\[\sum_k \min\left\{1, \frac{\delta_{k,n}^{\,3}}
	{d_{\mathbb H}(q,q_{k,n})^3}\right\}\leq C_2,\]
	where $C_2$ is independent of $n$ and $q$.
	Thus,
	\[\left\|	 \sum_k(f_{k,n}-g_{k,n})\right\|_\infty\leq C_3\,\omega_f(C\delta_n).\]
	Since 	$f$ has compact support, for each fixed $n$, only finitely many $f_{k,n}$ are nonzero. Thus, $	g_n =\sum_k g_{k,n}$	belongs to $\mathcal R_L(K)$.
	Now, 
	\[\|f-g_n\|_\infty\leq C_3\,\omega_f(C\delta_n).\]
	Since $\delta_n\longrightarrow0$	and $f$ is uniformly continuous on its compact support, $\omega_f(C\delta_n)\longrightarrow0.$
	Thus, $	\|f-g_n\|_\infty\longrightarrow0.$
	Hence $f$ belongs to the uniform closure of $\mathcal R_L(K)$.
	Since $f\in\mathcal A_{L}(K)$ was arbitrary, $	\mathcal A_{L}(K)
	\subseteq 	\overline{\mathcal R_{L}(K)}.$
	The reverse inclusion follows directly from the definition of
	$\mathcal A_{L}(K)$. Hence (1) holds.
	
	Now assume that (1) holds. It remains to prove the reverse inequality. Let $\varepsilon>0$. By the definition of the continuous slice-regular capacity, we may choose a compact axially symmetric set $L\subset U\setminus K^\circ$
	and a function $f\in\mathcal C_{\mathbb H}(L,1)$		such that
	\[|\partial_s f(\infty)|>\alpha_{\mathbb H}(U\setminus K^\circ)-\varepsilon.\]
	Since $L\subset U\setminus K^\circ,$ we have $L\cap K=\varnothing,$ and $L\cap(\mathbb H\setminus G)=\varnothing.$
	Hence $f$ is left slice regular on a neighborhood of
	$K\cup(\mathbb H\setminus U).$ In particular, $	f\in\mathcal A_{\mathbb H}(K).$	Now, by our assumption, $f\in\mathcal R_{\mathbb H}(K).$		
	Then by proposition \ref{11.7}, since $f$ is left slice regular
	on a neighborhood of $	K\cup(\mathbb H\setminus U),$
	we obtain a sequence of rational left slice-regular functions
	$r_n$ whose poles lie in $U\setminus K$	and such that
	$\|r_n-f\|_{K\cup(\mathbb H\setminus U)}\longrightarrow0.$
	Then $	\partial_s r_n(\infty)\longrightarrow \partial_s f(\infty).$
	Thus, \[	\alpha_{\mathbb H}(U\setminus K^\circ)-\varepsilon < \partial_s f(\infty)<	\alpha_{\mathbb H}(U\setminus K).\]
	Since $\varepsilon>0$ is arbitrary,
	\[\alpha_{\mathbb H}(U\setminus K^\circ)\leq\alpha_{\mathbb H}(U\setminus K).\]
	The reverse part of inequality follows by the monotonicity of $\alpha_{\mathbb{H}}$. Hence (2) holds
	
	Suppose (2) holds. then (3) is trivial.
	
	Now, assume that (3) holds. Fix $q\in\mathbb H$ and $r>1$. Put $s=(1+r)/2$. Since, $T(q,\delta)\subset T(q,s\delta)\subset T(q,r\delta),$
	then by monotonicity and (3) give
	\[	\frac{\alpha_{\mathbb H}(T(q,\delta)\setminus K^\circ)}
	{\alpha_{\mathbb H}(T(q,r\delta)\setminus K)}
	\leq \frac{\alpha_{\mathbb H}(T(q,s\delta)\setminus K^\circ)}
	{\alpha_{\mathbb H}(T(q,s\delta)\setminus K)}=1.\]
	Hence, $\limsup_{\delta\to0}{\alpha_{\mathbb H}(T(q,\delta)\setminus K^\circ)}{\alpha_{\mathbb H}(T(q,r\delta)\setminus K)}\leq1<\infty.$
	Thus (4) holds.
	
	Assume that (4) holds. 	Suppose, to the contrary, that $(5)$ is false. Then for every $r\geq1$ and every $m\geq1$ there are $q\in\mathbb H$ and	$\delta>0$ such that
	\[	\alpha_{\mathbb H}(T(q,\delta)\setminus K^\circ)>
	m\,\alpha_{\mathbb H}(T(q,r\delta)\setminus K).\]
	We construct sequences $q_n\in\mathbb H$ and $\delta_n>0$ such that
	$\delta_{n+1}\leq\frac{\delta_n}{2}$ and
	\[\frac{\alpha_{\mathbb H}(T(q_n,\delta_n)\setminus K^\circ)}
	{\alpha_{\mathbb H}(T(q_n,r_n\delta_n)\setminus K)}>n,
	\qquad r_n\longrightarrow\infty.\]
	Moreover, choosing the cells  inside the preceding ones, we  arrange 
	\[T(q_{n+1},r_n\delta_n)\subset T(q_n,\delta_n).\]
	Thus the cells are nested and, by compactness, their centres have a 	limit point $q_0\in K$ after passing to a sub sequence.	
	
	Fix $r\geq1$. Since $r_n\to\infty$, for all sufficiently large $n$, $T(q_0,r\delta_n) \subset T(q_n,r_n\delta_n)$	and $T(q_n,\delta_n)\subset T(q_0,3\delta_n).$	Hence, by monotonicity,
	\[\frac{\alpha_{\mathbb H}(T(q_0,3\delta_n)\setminus K^\circ)}{\alpha_{\mathbb H}(T(q_0,r\delta_n)\setminus K)}
	\geq \frac{\alpha_{\mathbb H}(T(q_n,\delta_n)\setminus K^\circ)}{\alpha_{\mathbb H}(T(q_n,r_n\delta_n)\setminus K)} >n.\]
	
	Since $r\geq1$ is arbitrary, (5) implies that for every fixed
	$r>1$,
	\[\limsup_{\delta\to0}\frac{\alpha_{\mathbb H}(T(q_0,\delta)\setminus K^\circ)}{\alpha_{\mathbb H}(T(q_0,r\delta)\setminus K)}=\infty.\]
	This contradicts  (4). Thus, (5) holds.
\end{proof}

\end{document}